\documentclass[11pt]{article}

\usepackage[sort]{cite}       
\usepackage{amssymb}          
\usepackage{latexsym}         
\usepackage{amsmath}          
\usepackage{amscd}            
\usepackage[latin1]{inputenc} 
\usepackage{eucal}            
\usepackage{bbm}              
\usepackage{theorem}          
\usepackage[ps,dvips]{xypic}  
\usepackage{diagxy}           
\usepackage{paralist}         

\title{States and representations in deformation quantization}

\author{\textbf{Stefan Waldmann}\thanks{E-mail:
    Stefan.Waldmann@physik.uni-freiburg.de}
  \\[0.1cm]
  Fakult{\"a}t f{\"u}r Mathematik und Physik\\
  Albert-Ludwigs-Universit{\"a}t Freiburg\\
  Physikalisches Institut\\
  Hermann Herder Stra{\ss}e 3\\
  D 79104 Freiburg\\
  Germany}

\date{August 2004\\[0.5cm] FR-THEP 2004/15}

\renewcommand{\mathbb}[1]{\mathbbm{#1}} 

\newcounter{comment}

\newcommand{\cond}[1] {{\rm \textbf{(#1)}}}

\newcommand{\cc}[1]      {\overline{{#1}}}              
\newcommand{\id}         {\operatorname{\mathsf{id}}}

\newcommand{\tr}         {\operatorname{\mathsf{tr}}}    
\newcommand{\ad}         {\operatorname{\mathrm{ad}}}    
\newcommand{\Ad}         {\operatorname{\mathrm{Ad}}}    
\newcommand{\Hom}        {\operatorname{\mathsf{Hom}}}   
\newcommand{\End}        {\operatorname{\mathsf{End}}}   
\newcommand{\SP}[1]      {\left\langle{#1}\right\rangle} 
\newcommand{\ring}[1]    {\mathsf{#1}}                 
\newcommand{\Unit}       {\mathbb{1}}                  
\newcommand{\Aut}        {\operatorname{\mathsf{Aut}}} 
 
\newcommand{\cl}         {\mathrm{cl}}

\newcommand{\I}          {\mathrm{i}}
\newcommand{\E}          {\mathrm{e}}
\newcommand{\D}          {\operatorname{\mathrm{d}}}

\newcommand{\group}[1]   {\mathrm{#1}}

\newcommand{\HdR}        {\mathrm{H}_{\scriptscriptstyle{\mathrm{dR}}}}
 
\newcommand{\Exp}        {\operatorname{\mathrm{Exp}}}
\newcommand{\Def}        {\operatorname{\mathrm{Def}}}

\newcommand{\starweyl}  {\mathbin{\star_{\scriptscriptstyle\mathrm{Weyl}}}}
\newcommand{\starwick}  {\mathbin{\star_{\scriptscriptstyle\mathrm{Wick}}}}
\newcommand{\starstd}   {\mathbin{\star_{\scriptscriptstyle\mathrm{Std}}}}

\newcommand{\starBCH}   {\mathbin{\star_{\scriptscriptstyle\mathrm{BCH}}}}

\newcommand{\weylrep}  {\operatorname{\varrho_{\scriptscriptstyle\mathrm{Weyl}}}}
\newcommand{\stdrep}   {\operatorname{\varrho_{\scriptscriptstyle\mathrm{Std}}}}
\newcommand{\wickrep}  {\operatorname{\varrho_{\scriptscriptstyle\mathrm{Wick}}}}

\newcommand{\Bimod}[5] {\sideset{^{\scriptscriptstyle{#1}}_{\scriptscriptstyle{#2}}}{^{\scriptscriptstyle{#4}}_{\scriptscriptstyle{#5}}}{\operatorname{#3}}}

\newcommand{\EA}   {\Bimod{}{}{\mathcal{E}}{}{\mathcal{A}}}

\newcommand{\BEA}  {\Bimod{}{\mathcal{B}}{\mathcal{E}}{}{\mathcal{A}}}
\newcommand{\HD}   {\Bimod{}{}{\mathcal{H}}{}{\mathcal{D}}}
\newcommand{\HpD}  {\Bimod{}{}{\mathcal{H}}{\prime}{\mathcal{D}}}
\newcommand{\CFB}  {\Bimod{}{\mathcal{C}}{\mathcal{F}}{}{\mathcal{B}}}
\newcommand{\FB}   {\Bimod{}{}{\mathcal{F}}{}{\mathcal{B}}}

\newcommand{\DGDp} {\Bimod{}{\mathcal{D}}{\mathcal{G}}{}{\mathcal{D}^\prime}}
\newcommand{\AM}   {\Bimod{}{\mathcal{A}}{\mathcal{M}}{}{}}
\newcommand{\AEpB} {\Bimod{}{\mathcal{A}}{\mathcal{E}}{\prime}{\mathcal{B}}}
\newcommand{\AEppB}{\Bimod{}{\mathcal{A}}{\mathcal{E}}{\prime\prime}{\mathcal{B}}}
\newcommand{\AccEB}{\Bimod{}{\mathcal{A}}{\cc{\mathcal{E}}}{}{\mathcal{B}}}
\newcommand{\AAA}  {\Bimod{}{\mathcal{A}}{\mathcal{A}}{}{\mathcal{A}}}
\newcommand{\BBB}  {\Bimod{}{\mathcal{B}}{\mathcal{B}}{}{\mathcal{B}}}

\newcommand{\IP}[4]{{\,}_{\scriptscriptstyle{#2}\!\!}\left\langle{{#1}}\right\rangle^{\scriptscriptstyle{#3}}_{\scriptscriptstyle{#4}}}

\newcommand{\SPEA}[1]  {\IP{{#1}}{}{\mathcal{E}}{\mathcal{A}}}
\newcommand{\BSPE}[1]  {\IP{{#1}}{\mathcal{B}}{\mathcal{E}}{}}
\newcommand{\SPFB}[1]  {\IP{{#1}}{}{\mathcal{F}}{\mathcal{B}}}
\newcommand{\SPD}[1]   {\IP{{#1}}{}{}{\mathcal{D}}}

\newcommand{\SPFEA}[1] {\IP{{#1}}{}{\mathcal{F}\otimes\mathcal{E}}{\mathcal{A}}}
\newcommand{\SPA}[1]   {\IP{{#1}}{}{}{\mathcal{A}}}
\newcommand{\BSP}[1]   {\IP{{#1}}{\mathcal{B}}{}{}}
\newcommand{\MnASP}[1]   {\IP{{#1}}{M_n(\mathcal{A})}{}{}}
\newcommand{\ASPccE}[1]  {\IP{{#1}}{\mathcal{A}}{\cc{\mathcal{E}}}{}}
\newcommand{\SPccEB}[1]  {\IP{{#1}}{}{\cc{\mathcal{E}}}{\mathcal{B}}}

\newcommand{\tensor}[1][{}]{\mathbin{\otimes_{\scriptscriptstyle{#1}}}}
\newcommand{\tensM}[1][{}] {\mathbin{\widehat{\otimes}_{\scriptscriptstyle{#1}}}}
\newcommand{\tensB}[1][{}] {\mathbin{\widetilde{\otimes}_{\scriptscriptstyle{#1}}}}

\newcommand{\rep}[1][{}]  {\sideset{^*}{_{#1}}{\operatorname{\textrm{-}\mathsf{rep}}}}
\newcommand{\Rep}[1][{}]  {\sideset{^*}{_{#1}}{\operatorname{\textrm{-}\mathsf{Rep}}}}

\newcommand{\Mod}  {\operatorname{\textrm{-}\mathsf{Mod}}}

\newcommand{\Pic}     {\operatorname{\mathsf{Pic}}}
\newcommand{\StrPic}  {\operatorname{\mathsf{Pic}^{\mathrm{str}}}}
\newcommand{\starPic} {\operatorname{\mathsf{Pic}^*}}

\newtheorem{lemma} {Lemma} [section]
\newtheorem{proposition} [lemma] {Proposition}

\newtheorem{theorem} [lemma] {Theorem}
\newtheorem{corollary} [lemma] {Corollary}
\newtheorem{definition}[lemma] {Definition}

\theorembodyfont{\rmfamily}
\newtheorem{example}[lemma]{Example}
\newtheorem{remark}[lemma]{Remark}
\newtheorem{question}[lemma]{Question}

\newenvironment{proof}[1][{}]{
  \par\noindent
  \textsc{Proof{#1}:}
}
{
  \hspace*{\fill} $\blacksquare$\newline
}

\numberwithin{equation}{section}

\setdefaultenum{\itshape i.)}{a.)}{I.)}{A.)}

\begin{document}

\maketitle

\begin{abstract}
    In this review we discuss various aspects of representation theory
    in deformation quantization starting with a detailed introduction
    to the concepts of states as positive functionals and the GNS
    construction. But also Rieffel induction of representations as
    well as strong Morita equivalence, the Dirac monopole and the
    strong Picard groupoid are discussed.
\end{abstract}

\newpage

%
%

\tableofcontents

%
%

\section{Introduction}
\label{sec:Intro}

Based on works of Weyl, Groenewold, Moyal, Berezin and others
\cite{weyl:1931a, groenewold:1946a, moyal:1949a, berezin:1975a,
  berezin:1975b, berezin:1975c} on the physical side and on
Gerstenhaber's deformation theory of associative algebras on the
mathematical side \cite{gerstenhaber.schack:1988a, gerstenhaber:1974a,
  gerstenhaber:1968a, gerstenhaber:1966a, gerstenhaber:1964a,
  gerstenhaber:1963a}, in the 1970's Bayen, Flato, Fr{\o}nsdal,
Lichnerowicz and Sternheimer coined the notion of a star product and
laid the foundations of deformation quantization in their seminal work
\cite{bayen.et.al:1978a}, see \cite{dito.sternheimer:2002a,
  gutt:2000a, weinstein:1994a} for recent reviews. Since then
deformation quantization developed into one of the most attractive and
successful quantization theories, both from the mathematical and
physical point of view.

The principle idea is to quantize the classical observable algebra
which is modeled by the smooth complex-valued functions on a Poisson
manifold $M$, see e.g.~\cite{cannasdasilva.weinstein:1999a,
  weinstein:1998a}, by simply replacing the commutative product of
functions by some non-commutative new product, the star product, now
depending on Planck's constant to control the non-commutativity, but
keeping the underlying vector space of observables. Thereby the
interpretation of the quantum observables is trivial: they are the
same elements of the observable algebra as classically. It turns out
that many other well-known quantization schemes can actually be cast
into this form whence it is fair to say that deformation quantization
is more a theory of quantization itself rather than a particular
quantization scheme.

In formal deformation quantization one has very strong existence and
classification results for the star products which depend on $\hbar$
in the sense of a formal power series. For the symplectic case the
general existence was shown first by DeWilde and
Lecomte~\cite{dewilde.lecomte:1988a, dewilde.lecomte:1983b}, later by
Fedosov~\cite{fedosov:1985a, fedosov:1986a, fedosov:1989a,
  fedosov:1994a} and Omori, Maeda, and
Yoshioka~\cite{omori.maeda.yoshioka:1991a}. The case where the
classical phase space is a general Poisson manifold turned out to be
much more difficult and was finally solved by
Kontsevich~\cite{kontsevich:1997:pre}, see also
\cite{kontsevich:2001a, kontsevich:1999a} by proving his formality
conjecture \cite{kontsevich:1997a}. The classification of star
products was obtained again first for the symplectic case by Nest and
Tsygan~\cite{nest.tsygan:1995a, nest.tsygan:1995b}, Bertelson, Cahen
and Gutt~\cite{bertelson.cahen.gutt:1997a},
Deligne~\cite{deligne:1995a}, Weinstein and
Xu~\cite{weinstein.xu:1998a}. The classification in the Poisson case
follows also from the formality theorem of Kontsevich
\cite{kontsevich:1997:pre}. For an interpretation of Kontsevich's
formality in terms of the Poisson-sigma model as well as globalization
aspects see the work of Cattaneo, Felder and
Tomassini~\cite{cattaneo.felder:2000a,
  cattaneo.felder.tomassini:2002b, cattaneo.felder.tomassini:2002a}.

It should also be mentioned that star products find physical
applications far beyond the original quantization problem: recently
the most prominent applications come from non-commutative geometry
\cite{connes:1994a} and the non-commutative field theory models
arising from it. Here one endows the space-time manifold with the
non-commutative star product and studies field theories on this
non-commtutative space-time, see e.g.
\cite{bahns.doplicher.fredenhagen.piacitelli:2003a:pre,
  jurco.schupp.wess:2000a, jurco.schupp.wess:2001a,
  jurco.et.al.:2001a, seiberg.witten:1999a,
  connes.douglas.schwarz:1998a} and references therein.

For the quantization problem it is for physical reasons not enough to
consider only the space of observables and their quantization. One
also needs a notion for the states and their quantization. To give an
overview over the concepts of states in deformation quantization is
therefor the main topic of this review. It turns out that the question
of states is, as in any other quantum theory based on the notion of
observables, intimately linked to the question of representations of
the observable algebra. A systematic investigation of representations
of the deformed algebras started with the work
\cite{bordemann.waldmann:1998a} leading to the general representation
theory for algebras with involution defined over a ring $\ring{C} =
\ring{R}(\I)$ with an ordered ring $\ring{R}$ as developed in a series
of articles \cite{bordemann.waldmann:1997b,
  bordemann.roemer.waldmann:1998a, bordemann.neumaier.waldmann:1999a,
  bordemann.neumaier.pflaum.waldmann:2003a,
  bordemann.neumaier.waldmann:1998a, bursztyn.waldmann:2000a,
  bursztyn.waldmann:2001b, bursztyn.waldmann:2001a,
  bursztyn.waldmann:2003a:pre, bursztyn.waldmann:2004a,
  bursztyn.waldmann:2002a, waldmann:2000a, waldmann:2001a,
  waldmann:2002a, waldmann:2003b:pre, waldmann:2003a:pre}. The purpose
of this work is to give an introduction to these concepts and discuss
some of the basic results in representation theory of star product
algebras. Since the techniques are fairly general many of the results
will find applications also in other areas of mathematical physics.

The plan of this review is as follows: In Section~\ref{sec:Motivation}
we briefly remind on the basic notions of states as positive
functionals and representations of observable algebras and discuss the
necessity of studying them. Section~\ref{sec:Background} gives then
the algebraic background on ordered rings, $^*$-algebras and notions
of positivity which will be crucial throughout this work.
Section~\ref{sec:ExamplesPosFun} is devoted to examples of positive
functionals from deformation quantization. In
Section~\ref{sec:ClassLimPosFun} we discuss the deformation and the
classical limit of positive functionals and introduce the important
notion of a positive deformation. Section~\ref{sec:GNS} establishes
the relation between positive functionals and representations via the
GNS construction of representations.
Section~\ref{sec:GeneralRepTheory} starts with an introduction to more
advanced topics in representation theory like Rieffel induction and
related tensor product constructions. This will be used in
Section~\ref{sec:StrongMoritaEquivalence} to establish the notion of
strong Morita equivalence and the strong Picard groupoid which encodes
the whole Morita theory. In Section~\ref{sec:SMEStarProducts} we
discuss the Morita theory for deformed and more specifically for star
product algebras yielding a new look at Dirac's monopole and the
corresponding charge quantization. Finally, Section~\ref{sec:Outlook}
contains several open questions and further ideas related to
representation theory in deformation quantization. We have included an
extensive though by no means complete bibliography.  For more details
and references one should also consult the Deformation Quantization
Homepage.

\noindent
\textbf{Acknowledgments:} It is a pleasure to thank Didier Arnal,
Giuseppe Dito and Paco Turrubiates for many discussions during my stay
in Dijon as well as for encouraging me to write these notes.  I also
would like to thank the University of Dijon for the warm hospitality.

%
%

\section{Motivation: Why states and representations?}
\label{sec:Motivation}

In this section we shall give some well-known remarks on the
quantization problem and the general approach to quantum theory based
on the notion of an observable algebra and specialize this to
deformation quantization.

%
%

\subsection{Observables and states}
\label{subsec:ObservablesStates}

When we want to learn something about the relation between the
classical and quantum description of a physical system we should first
discuss the similarities and differences as detailed as possible. Here
we follow the idea that the \emph{observables}, i.e. the possible
measurements one can perform on the system, characterize the system
itself. Moreover, the algebraic structure of the observables
determines what the possible \emph{states} of the system can be. From
this point of view classical and quantum theory behave quite similar.
We illustrate this for a system with finitely many degrees of freedom
though the main results will easily generalize to field theories or
thermodynamical systems.

For the classical side, we model the algebra of observables by the
\emph{Poisson $^*$-algebra} of complex-valued smooth functions on a
manifold $M$, the \emph{phase space} of the system. The
$^*$-involution will always be the pointwise complex conjugation. An
element in the observable algebra is called \emph{observable} if it is
a real-valued function $f = \cc{f}$.  The structure of a Poisson
bracket $\{\cdot, \cdot\}$ for the smooth functions is equivalent to a
Poisson structure on the manifold, i.e. a smooth anti-symmetric
$2$-tensor field $\pi \in \Gamma^\infty(\Lambda^2 TM)$ with $[\![\pi,
\pi]\!] = 0$, where $[\![\cdot, \cdot]\!]$ denotes the Schouten
bracket, and the relation is $\{f, g\} = \pi(\D f, \D g)$, see
e.g.~\cite{cannasdasilva.weinstein:1999a, weinstein:1998a,
  landsman:1998a} for more details and references on Poisson geometry.
Of course there are situations where the class of functions describing
the system most adequately may be a different one. The \emph{pure}
states are then the points of the phase space while the \emph{mixed}
states correspond to more general positive Borel measures on $M$. The
\emph{physical spectrum} of an observable $f$, i.e. the possible
values of $f$ in a measurement, coincides with its mathematical
spectrum, namely the set of values of the function. Finally, the
\emph{expectation value} of an observable $f$ in a pure state $x \in
M$ is given by the evaluation $E_x(f) = f(x)$ while in a mixed state
$\mu$ the expectation value is $E_\mu(f) = \int_M f(x) \D\mu(x)$. The
crucial feature of a classical observable algebra is its commutativity
which allows to have sharp measurement of \emph{all} observables in a
pure state.

In quantum theory, the observables are usually modeled by the
\emph{$^*$-algebra} of bounded operators $\mathfrak{B}(\mathfrak{H})$
on a complex Hilbert space $\mathfrak{H}$, or, more general, some
$^*$-algebra of densely defined and possibly unbounded operators on
$\mathfrak{H}$. It is clear that one has to specify in which sense one
wants to understand `$^*$-algebra' if the operators are only densely
defined, but these technical aspects can be made precise in a
completely satisfactory way, see e.g.~\cite{schmuedgen:1990a}. The
observable elements in the observable algebra are then the
self-adjoint operators.  The pure states are now complex rays in the
Hilbert space. Usually, not all rays have physical relevance as the
vectors defining them have to be in the domain of the observables of
interest which may only be a dense subspace of $\mathfrak{H}$. From
this point of view only a pre-Hilbert space is needed to describe the
physically relevant states while the Hilbert space in the background
is needed to have a `good' spectral calculus. More generally, mixed
states are described by density matrices $\varrho$, i.e. positive
trace class operators with trace $1$. Indeed, the pure states are just
the rank-one projection operators from this point of view. The
spectrum is now the spectrum in the sense of self-adjoint operators.
Finally, the expectation value of an observable $A$ in the pure state
defined by $\phi \in \mathfrak{H}$ is $E_\phi(A) = \frac{\SP{\phi,
    A\phi}}{\SP{\phi,\phi}}$ and in a mixed state $\varrho$ it is
given by $E_\varrho(A) = \tr(\varrho A)$.

Up to now this is the standard description of classical and quantum
theory as it can be found in text books. However, the way we presented
it allows for a uniform framework for both theories which is better
suited concerning questions of quantization and classical limit.
Indeed, the structure of the observables is in both cases encoded in a
(unital) complex $^*$-algebra $\mathcal{A}$. The difference between
the classical and quantum side is that $\mathcal{A}$ is
non-commutative for the quantum theory as we have to incorporate
\emph{uncertainty relations} while in the classical situation
$\mathcal{A}$ is commutative but has an additional structure, the
Poisson bracket.  The states are now identified with the expectation
value functionals and can thus be described by \emph{normalized
  positive linear functionals} $\omega: \mathcal{A} \longrightarrow
\mathbb{C}$, i.e.  linear functionals such that $\omega(a^*a) \ge 0$
for all $a \in \mathcal{A}$ and $\omega(\Unit) = 1$. The question
whether a state is pure or mixed becomes now the question whether
$\omega$ can be decomposed in a non-trivial way into a convex
combination $\omega = c_1 \omega_1 + c_2 \omega_2$ of two other states
$\omega_1$, $\omega_2$. The expectation value of an observable $a$ in
a state $\omega$ is then $E_\omega(a) = \omega(a)$. Clearly, all our
above examples of `states' fit into this framework.  Unfortunately,
one has to leave this purely algebraic framework as soon as one wants
to have a reasonable notion for `spectrum'. Here one has to impose
some analytical conditions on the $^*$-algebra in question in order to
get physically acceptable answers. In the above examples this
corresponds to the choice of an appropriate class of functions on the
phase space on the classical side and the questions about
self-adjointness on the quantum side. Typically, some $C^*$-algebraic
structures behind $\mathcal{A}$ will be responsible for a good
spectral calculus.

Except for this last difficulty the \emph{problem of quantization} can
now be seen as the task to construct the quantum observable algebra
out of the knowledge of the classical observable algebra. The above
formulation then gives automagically a construction of the states as
well, since the states of a $^*$-algebra are defined in a uniform way,
whether the algebra is commutative or not. In this sense the algebraic
structure of the observables \emph{determines} the possible states
whence for quantization it is sufficient to find the observable
algebra.

Before we discuss one of the approaches to quantization in more detail
let us state clearly that from a physicists point of view the whole
question about quantization is in some sense completely artificial: by
our present knowledge the world \emph{is} already quantum whence there
is nothing left to be quantized. The true physical problem is the
inverse question: why and how does a classical world emerge out of
this quantum world, at least for certain scales of energy, momentum,
length, time, etc.? Nevertheless, as physicist one is still interested
in quantization since up to now we have not developed a sufficient
intuition which would allow us to formulate quantum theories \emph{a
  priori} without the usage of classical counterparts, except for very
few cases. We shall not speculate too much on the more philosophical
question why this is (still) the case. Instead, we consider
`quantization' as a pragmatic approach to find relevant quantum
descriptions for physical systems we are interested in.

%
%

\subsection{Super-positions and super-selection rules}
\label{subsec:SuperpositionSelections}

Once we succeeded in finding the quantum algebra of observables and
having determined the states as its positive functionals do we then
have a complete quantum theory? The answer is no, there is still one
important piece missing, most crucial for quantum physics: We still
need the \emph{super-position principle} for the (pure) states. In the
usual Hilbert space formulation of quantum physics one simply takes
complex linear combinations of the state \emph{vectors} in order to
encode the super-position of the states they represent. Note that this
can not be done so simply in our more advanced formulation where
states are identified with their expectation value functionals. Of
course we can take \emph{convex combinations} of positive functionals
and get again positive functionals but this does \emph{not} correspond
to the super-position of the states but to a mixed state. Thus we need
this additional linear structure of the Hilbert space which is
precisely the reason why one has to represent the algebra of
observables on a Hilbert space such that the positive functionals
become the expectation value functionals for vector states. We want
to be able to write
\begin{equation}
    \label{eq:ExpectationVectorState}
    \omega(a) = \frac{\SP{\phi, \pi(a) \phi}}{\SP{\phi, \phi}}
\end{equation}
for some $^*$-representation $\pi$ on some Hilbert space
$\mathfrak{H}$ and some vector $\phi \in \mathfrak{H}$. Clearly, at
this stage we only need a pre-Hilbert space structure.

But which $^*$-representation shall we choose? In particular, for two
given positive functionals $\omega_1$ and $\omega_2$ can we always
find a $^*$-representation $(\mathfrak{H}, \pi)$ such that both states
$\omega_1$, $\omega_2$ can be written in the form
\eqref{eq:ExpectationVectorState} with some $\phi_1, \phi_2 \in
\mathfrak{H}$ in order to form their super-positions? If we require in
addition that the algebra acts irreducibly then the answer, in
general, is no. Dropping the irreducibility gives an easy answer
provided we can find $^*$-representations $(\mathfrak{H}_i, \pi_i)$
for each $\omega_i$ separately, since then we simply can take the
direct orthogonal sum of the $^*$-representations. Then however,
super-positions of the vectors $\phi_1$, $\phi_2$ will not produce any
interesting interference cross terms. There will be no transitions
between these two states. This phenomena is called a
\emph{super-selection rule}: one can not superpose the two states in a
non-trivial way.  It is well-known from quantum field theory that this
may happen indeed, see e.g. the discussion in \cite{haag:1993a}. The
presence of super-selection rules is usually interpreted as the
existence of non-trivial charges.  Mathematically speaking it
corresponds to the existence of \emph{inequivalent} (faithful)
irreducible $^*$-representations.  Note that in order to `see' the
super-selection rules it is \emph{not} enough to choose one particular
$^*$-representation from the beginning. So the problem of choosing a
$^*$-representation `is not a bug, it is a feature'.

In usual quantum mechanics of an uncharged particle moving in
Euklidian $\mathbb{R}^n$ super-selection rules are absent: this is the
statement of a classical theorem of von Neumann. Note however that
this statement is only true after some effort involving the completion
of the observable algebra to a $C^*$-algebra, the Weyl algebra, see
e.g. the discussion in \cite[Sect.~I.1]{haag:1993a}.

One should take this non-trivial result also as a \emph{warning}:
since the absence of super-selection rules in the above case is the
consequence of some rather strong topological context one has to
expect that one might see `super-selection rules' which are artifacts
of a possible `non-completeness' of the observable algebra in the
sense that they vanish immediately after one passes to an appropriate
completion.

As conclusion we see that we have to understand the whole
$^*$-representation theory of the observable algebra and determine the
`hard' super-selection rules, i.e. those which survive some
(physically motivated) completions. Of course, this has to be made
more specific in the sequel. In fact, in formal deformation
quantization this turns out to be a highly non-trivial issue.

%
%

\subsection{Deformation quantization}
\label{subsec:DeformationQuantization}

The main idea of \emph{deformation quantization}, as formulated in
\cite{bayen.et.al:1978a} is to construct the quantum observable
algebra as a \emph{non-commutative associative deformation} of the
classical observable algebra in the sense of Gerstenhaber's
deformation theory \cite{hazewinkel.gerstenhaber:1988a,
  gerstenhaber:1974a, gerstenhaber:1968a, gerstenhaber:1966a,
  gerstenhaber:1964a} where the first non-trivial term in the
commutator of the deformation is the classical Poisson bracket and the
deformation parameter is Planck's constant $\hbar$.  Thus the
classical Poisson bracket is the `shadow' of the non-commutativity of
quantum theory. Roughly speaking, deformation means that we endow the
same underlying vector space of the classical observable algebra with
a family of new products $\star_\hbar$, called \emph{star products},
which depend on the deformation parameter $\hbar$ in such a way that
for $\hbar = 0$ we recover the classical commutative product
structure.

There are at least two flavours of this quantization scheme:
\emph{strict deformations} and \emph{formal deformations}.  While in
strict deformation quantization, see e.g.~\cite{rieffel:1993a,
  rieffel:1989a, landsman:1998a}, one wants the products $\star_\hbar$
to depend in a continuous way on $\hbar$, usually within a
$C^*$-algebraic framework, in formal deformation quantization the
dependence is in the sense of formal power series. At least in some
good cases the formal deformations can be seen as an asymptotic
expansion of the strict ones. Of course, the deformation parameter,
being identified with Planck's constant $\hbar$, should not be a
formal parameter but a physical quantity. Thus starting with formal
deformations one should be able to establish at some point the
convergence of the formal series. In general this turns out to be a
rather delicate problem usually depending in a very specific way on
the particular example one considers, see e.g.
\cite{cahen.gutt.rawnsley:1990a, gutt:2000a,
  bordemann.brischle.emmrich.waldmann:1996b} and references therein.
So, unfortunately, not very much can be said about this point in
general.  On the other hand, the advantage of formal deformations is
that we can decide at which point we want to impose the convergence
conditions.  This gives usually more freedom in the beginning. In the
following we shall always consider the formal framework.

After these general remarks we can now state the definition of a star
product according to \cite{bayen.et.al:1978a} and recall some of the
basic results:
\begin{definition}[Star Product]
    \label{definition:StarProduct}
    Let $(M, \pi)$ be a Poisson manifold. Then a formal star product
    $\star$ for $(M, \pi)$ is an associative
    $\mathbb{C}[[\lambda]]$-bilinear multiplication for
    $C^\infty(M)[[\lambda]]$, written as
    \begin{equation}
        \label{eq:StarProduct}
        f \star g = \sum_{r=0}^\infty \lambda^r C_r(f,g)
    \end{equation}
    for $f, g \in C^\infty(M)[[\lambda]]$, such that
    \begin{compactenum}
    \item $C_0(f, g) = fg$ and $C_1(f, g) - C_1(g, f) = \I\{f,g\}$,
    \item $1 \star f = f = f \star 1$,
    \item $C_r$ is a bidifferential operator,
    \item $\cc{f \star g} = \cc{g} \star \cc{f}$.
    \end{compactenum}
\end{definition}
The first condition reflects the \emph{correspondance principle} while
the last condition is sometimes omitted. To stress the last condition
we shall call star products satisfying this condition also
\emph{Hermitian star products}. In the following we shall mainly be
interested in Hermitian star products.

If $S = \id + \sum_{r=1}^\infty \lambda^r S_r$ is a formal series of
differential operators $S_r: C^\infty(M) \longrightarrow C^\infty(M)$
with the property that $S_r$ vanishes on constants then for a given
star product $\star$ the definition
\begin{equation}
    \label{eq:EquivStarProd}
    f \star' g = S^{-1} (Sf \star Sg)
\end{equation}
defines again a star product deforming the same Poisson bracket as
$\star$. This is an immediate computation. If in addition $S \cc{f} =
\cc{Sf}$ then $\star'$ is a Hermitian star product if $\star$ is
Hermitian. If two star products are related by such an operator they
are called \emph{equivalent} or \emph{$^*$-equivalent}, respectively.

The existence of star products as well as their classification up to
equivalence is now well-understood:
\begin{theorem}[Existence of Star Products]
    \label{theorem:ExistenceStarProducts}
    On any Poisson manifold there exists a Hermitian star product.
\end{theorem}
The first proofs of this theorem for the case of symplectic manifolds
were obtained by DeWilde and Lecomte \cite{dewilde.lecomte:1983b,
  dewilde.lecomte:1988a} and independently by Fedosov
\cite{fedosov:1994a, fedosov:1989a, fedosov:1986a, fedosov:1985a} and
Omori, Maeda and Yoshioka \cite{omori.maeda.yoshioka:1991a}. The much
more involved existence in the Poisson case is a consequence of
Kontsevich's formality theorem \cite{kontsevich:1999a,
  kontsevich:1997:pre}.

The classification up to equivalence was first obtained for the
symplectic case by Nest and Tsygan \cite{nest.tsygan:1995a,
  nest.tsygan:1995b}, Bertelson, Cahen and Gutt
\cite{bertelson.cahen.gutt:1997a}, Deligne \cite{deligne:1995a} (see
also \cite{gutt.rawnsley:1999a, neumaier:2002a}), and Weinstein and Xu
\cite{weinstein.xu:1998a}. Here the equivalence classes are shown to
be in canonical bijection with formal series in the second (complex)
de Rham cohomology: One has a \emph{characteristic class}
\begin{equation}
    \label{eq:CharacteristicClass}
    c: \star \; \longmapsto \; c(\star) 
    \in \frac{[\omega]}{\I\lambda} + \HdR^2(M, \mathbb{C})[[\lambda]],
\end{equation}
where the origin of the above affine space is chosen by convention and
two symplectic star products are equivalent if and only if their
characteristic classes coincide. Moreover, a star product is
equivalent to a Hermitian star product if and only if its
characteristic class is imaginary \cite{neumaier:2002a}.  The above
classification is a particular case for the classification in the
Poisson case which is also obtained from Kontsevich's formality. In
general, the equivalence classes of star products are in bijection
with the formal deformations of the Poisson bivector modulo formal
diffeomorphisms \cite{kontsevich:1997:pre}.  Finally, Hermitian star
products are equivalent if and only if they are $^*$-equivalent, see
e.g. \cite{bursztyn.waldmann:2002a}.

Having understood this, the next problem is to define positive
functionals and $^*$-representations in this context and determine (as
far as possible) the representation theory of the star product
algebras. The first attempt of considering $\mathbb{C}$-linear
positive functionals $\omega: C^\infty(M)[[\lambda]] \longrightarrow
\mathbb{C}$ with $\omega(\cc{f} \star f) \ge 0$ turns out to be too
naive: One is faced immediately with convergence problems or one has
to ignore higher orders in $\lambda$ at some point. Both problems
limit this attempt too much. It is simply the wrong category and we
should better take the formal power series serious. Thus the better
choice is to look for \emph{$\mathbb{C}[[\lambda]]$-linear}
functionals
\begin{equation}
    \label{eq:POsFunInDefQuant}
    \omega: C^\infty(M)[[\lambda]] 
    \longrightarrow \mathbb{C}[[\lambda]].
\end{equation}
Then, of course, we have to \emph{define} what we mean by
\emph{positivity}.  Here we can use the following simple but crucial
fact that the real formal power series $\mathbb{R}[[\lambda]]$ are in
a natural way an \emph{ordered ring}: one defines for $a \in
\mathbb{R}[[\lambda]]$
\begin{equation}
    \label{eq:FormalSeriesPositive}
    a = \sum_{r=r_0}^\infty \lambda^r a_r > 0
    \quad
    \textrm{if and only if}
    \quad
    a_{r_0} > 0.
\end{equation}

This allows us to speak of positive linear functionals in a meaningful
way and follow the above program once we have adapted the concepts of
$^*$-representations etc. to $^*$-algebras defined over such an
ordered ring. To provide such a framework, extending the usual
framework of $^*$-algebras over $\mathbb{C}$, is the main objective of
this paper.

Let us conclude this section with moderate warnings on what we have to
expect to get from this approach. Clearly, we have to expect artifacts
like many inequivalent $^*$-representations which will disappear in a
convergent and more topological context. In some sense the situation
might turn out to be even more involved than for $^*$-algebras over
$\mathbb{C}$. Thus it will be a difficult task to detect the `hard'
super-selection rules in this framework. On the other hand, the
framework will hopefully be wide enough to contain all physically
interesting $^*$-representations. Obstructions found in this general
framework will be certainly difficult to overcome in even more strict
frameworks.

%
%

\section{Algebraic background: $^*$-Algebras over ordered rings}
\label{sec:Background}

In this section we set up the basic theory of $^*$-algebras over
ordered rings in order to have a unified approach for $^*$-algebras
over $\mathbb{C}$ like e.g. $C^*$-algebras or more general
$O^*$-algebras and the star product algebras being $^*$-algebras over
$\mathbb{C}[[\lambda]]$. The well-known case of operator algebras
(bounded or unbounded, see e.g.~\cite{sakai:1971a,
  kadison.ringrose:1997a, bratteli.robinson:1987a,
  kadison.ringrose:1997b, schmuedgen:1990a}) will be the motivation
and guideline throughout this section.

%
%

\subsection{Ordered rings}
\label{subsec:OrderedRings}

Let us first recall the definition and some basic properties of
ordered rings being a slight generalization of ordered fields, see
e.g.~\cite[Sect.~5.1]{jacobson:1985a}.
\begin{definition}[Ordered Ring]
    \label{definition:OrderedRing}
    An ordered ring $(\ring{R}, \ring{P})$ is a commutative,
    associative unital ring $\ring{R}$ together with a subset
    $\ring{P} \subset \ring{R}$, the positive elements, such that
    \begin{compactenum}
    \item $\ring{R} = -\ring{P} \mathbin{\dot{\cup}} \{0\}
        \mathbin{\dot{\cup}} \ring{R}$ (disjoint union),
    \item $\ring{P} \cdot \ring{P} \subseteq \ring{P}$ and $\ring{P} +
        \ring{P} \subseteq \ring{P}$.
    \end{compactenum}
\end{definition}
The subset $\ring{P}$ induces an \emph{ordering} defined by $a < b$ if
$b-a \in \ring{P}$. The symbols $\le$, $\ge$ $>$ will then be used in
the usual way.  In the following, we will fix an ordered ring
$\ring{R}$ and denote by
\begin{equation}
    \label{eq:CdefRi}
    \ring{C} = \ring{R}(\I) = \ring{R} \oplus \I \ring{R},
    \quad
    \textrm{where}
    \quad
    \I^2 = -1,
\end{equation}
the ring-extension by a square-root of $-1$. In $\ring{C}$ we have the
usual \emph{complex conjugation}
\begin{equation}
    \label{eq:ccDef}
    z = a + \I b \; \longmapsto \; \cc{z} = a - \I b,
\end{equation}
where $a,b \in \ring{R}$ and $\ring{R}$ is considered as a sub-ring of
$\ring{C}$ in the usual way.
\begin{remark}[Characteristics and Quotient Fields of Ordered Rings]
    \label{remark:OrderedRing}
    ~
    \begin{compactenum}
    \item If $a \ni \ring{R}$, $a \ne 0$ then $a^2 > 0$. Hence $1 =
        1^2 > 0$ and thus $1 + \cdots + 1 = n > 0$ for all $n \in
        \mathbb{N}$. Thus it follows that $\mathbb{Z} \subseteq
        \ring{R}$ whence $\ring{R}$ has \emph{characteristic zero}.
    \item Moreover, $z\cc{z} = a^2 + b^2 > 0$ for $z = a + \I b \in
        \ring{C}$, $z \ne 0$.  It also follows that the characteristic
        of $\ring{C}$ is zero, too.
    \item If $a,b \ne 0$ in $\mathbb{R}$ then we have four cases $a >
        0$ and $b>0$, $a>0$ and $b<0$, $a<0$ and $b>0$, $a<0$ and
        $b<0$. In each case we obtain $ab \ne 0$ whence $\ring{R}$ has
        no \emph{zero-divisors}. The same holds for $\ring{C}$. Hence
        we can pass to the \emph{quotient fields} $\hat{\ring{R}}$ and
        $\hat{\ring{C}}$, respectively. The field $\hat{\ring{R}}$ is
        canonically ordered and $\ring{R} \hookrightarrow
        \hat{\ring{R}}$ is order preserving. Finally, $\hat{\ring{C}}
        = \hat{\ring{R}}(\I)$.
    \end{compactenum}
\end{remark}
\begin{definition}[Archimedean Ordering]
    \label{definition:Archimedean}
    An ordered ring $\ring{R}$ is called Archimedean if for $a,b > 0$
    there is a $n \in \mathbb{N}$ with $na > b$. Otherwise $\ring{R}$
    is called non-Archimedean.
\end{definition}
\begin{example}[Ordered Rings]
    \label{example:OrderedRings}
    ~
    \begin{compactenum}
    \item $\mathbb{Z}$ is the smallest ordered ring and contained in
        any other. Clearly $\mathbb{Z}$ is Archimedean.
    \item $\mathbb{Q}$ and $\mathbb{R}$ are Archimedean ordered rings,
        even ordered fields.
    \item $\mathbb{R}[[\lambda]]$ is non-Archimedean as $n\lambda < 1$
        but $\lambda > 0$. The quotient field is the field of the
        formal Laurent series $\mathbb{R}(\!(\lambda)\!)$.
    \item More generally, if $\ring{R}$ is an ordered ring then
        $\ring{R}[[\lambda]]$ is canonically ordered again by the
        analogous definition as in \eqref{eq:FormalSeriesPositive} and
        it is always non-Archimedean. This already indicates that
        ordered rings and formal deformations fit together nicely, for
        the price of non-Archimedean orderings.
    \end{compactenum}
\end{example}

%
%

\subsection{Pre-Hilbert spaces}
\label{subsec:PreHilbertSpaces}

Having an ordered ring we have the necessary notion of positivity in
order to define pre-Hilbert spaces generalizing the usual complex
case.
\begin{definition}[Pre-Hilbert Space]
    \label{definition:PreHilbertSpace}
    A $\ring{C}$-module $\mathcal{H}$ with a map $\SP{\cdot,\cdot}:
    \mathcal{H} \times \mathcal{H} \longrightarrow \ring{C}$ is called
    pre-Hilbert space over $\ring{C}$ if
    \begin{compactenum}
    \item $\SP{\cdot,\cdot}$ is $\ring{C}$-linear in the second
        argument,
    \item $\SP{\phi, \psi} = \cc{\SP{\psi, \phi}}$ for all $\phi, \psi
        \in \mathcal{H}$,
    \item $\SP{\phi, \phi} > 0$ for $\phi \ne 0$.
    \end{compactenum}
\end{definition}
A map $A: \mathcal{H}_1 \longrightarrow \mathcal{H}_2$ is called
\emph{adjointable} if there exists a (necessarily unique) map $A^*:
\mathcal{H}_2 \longrightarrow \mathcal{H}_1$ with
\begin{equation}
    \label{eq:Adjointable}
    \SP{\phi, A \psi}_2 = \SP{A^*\phi, \psi}_1
\end{equation}
for all $\phi \in \mathcal{H}_2$, $\psi \in \mathcal{H}_1$. Clearly,
adjointable maps are $\ring{C}$-linear and we have the usual rules for
adjoints, i.e.
\begin{equation}
    \label{eq:AdjointStarInvolution}
    (zA + wB)^* = \cc{z} A^* + \cc{w} B^*,
    \quad
    (AB)^* = B^* A^*,
    \quad
    \textrm{and}
    \quad
    (A^*)^* = A,
\end{equation}
where the existence of the adjoints on the left is implied. The set of
all adjointable maps is denoted by
\begin{equation}
    \label{eq:BHeinsHzwei}
    \mathfrak{B}(\mathcal{H}_1, \mathcal{H}_2)
    =
    \{A: \mathcal{H}_1 \longrightarrow \mathcal{H}_2
    \; | \; A \; \textrm{is adjointable} \}
\end{equation}
\begin{equation}
    \label{eq:BH}
    \mathfrak{B}(\mathcal{H}) = \mathfrak{B}(\mathcal{H},\mathcal{H}).
\end{equation}
Clearly, $\mathfrak{B}(\mathcal{H})$ is a unital sub-algebra of all
$\ring{C}$-linear endomorphisms of $\mathcal{H}$.

Particular examples of adjointable operators are the rank-one and
finite rank operators. For $\phi \in \mathcal{H}_1$ and $\psi \in
\mathcal{H}_2$ we define the \emph{rank-one operator}
\begin{equation}
    \label{eq:Thetaphipsi}
    \Theta_{\psi, \phi}: \mathcal{H}_1 \ni \chi 
    \; \longmapsto \; \psi \SP{\phi, \chi} \in \mathcal{H}_2,
\end{equation}
which is clearly adjointable with adjoint $\Theta_{\psi, \phi}^* =
\Theta_{\phi, \psi}$. Moreover, we define the \emph{finite-rank
  operators}
\begin{equation}
    \label{eq:FiniteRankOperators}
    \mathfrak{F}(\mathcal{H}_1, \mathcal{H}_2)
    =
    \ring{C}\textrm{-}\mathrm{span}
    \{ \Theta_{\psi, \phi} \; | \; 
    \phi \in \mathcal{H}_1, \psi \in \mathcal{H}_2 
    \}
\end{equation}
and set
\begin{equation}
    \label{eq:FiniteRankH}
    \mathfrak{F}(\mathcal{H}) 
    = \mathfrak{F}(\mathcal{H}, \mathcal{H}).
\end{equation}
We have $\mathfrak{F}(\mathcal{H}_1, \mathcal{H}_2) \subseteq
\mathfrak{B}(\mathcal{H}_1, \mathcal{H}_2)$ and
$\mathfrak{F}(\mathcal{H}) \subseteq \mathfrak{B}(\mathcal{H})$. In
general, these inclusions are proper:
\begin{example}[Standard Pre-Hilbert Space]
    \label{example:BHFH}
    Let $\Lambda$ be a set and consider the free $\ring{C}$-module
    generated by $\Lambda$, i.e. $\mathcal{H} = \ring{C}^{(\Lambda)} =
    \oplus_{\lambda \in \Lambda} \ring{C}_\lambda$ with
    $\ring{C}_\lambda = \ring{C}$ for all $\lambda$. Then
    $\mathcal{H}$ becomes a pre-Hilbert space by
    \begin{equation}
        \label{eq:SpStandard}
        \SP{(x_\lambda), (y_\lambda)}
        = \sum_{\lambda \in \Lambda} \cc{x}_\lambda y_\lambda.
    \end{equation}
    In general $\mathfrak{F}(\mathcal{H}) \subsetneq
    \mathfrak{B}(\mathcal{H})$ unless $\#\Lambda = n < \infty$. In this
    case $\mathfrak{F}(\mathcal{H}) = \mathfrak{B}(\mathcal{H}) \cong
    M_n(\ring{C})$.
\end{example}

%
%

\subsection{$^*$-Algebras}
\label{subsec:StarAlgebras}

The sub-algebra $\mathfrak{B}(\mathcal{H}) \subseteq
\End(\mathcal{H})$ will be the motivating example of a $^*$-algebra:
\begin{definition}[$^*$-Algebra]
    \label{definition:StarAlgebra}
    An associative algebra $\mathcal{A}$ over $\ring{C}$ together with
    a $\ring{C}$-anti-linear, involutive anti-automorphism $^*:
    \mathcal{A} \longrightarrow \mathcal{A}$ is called a
    $^*$-algebra and $^*$ is called the $^*$-involution of
    $\mathcal{A}$. A morphism of $^*$-algebras is a morphism $\phi:
    \mathcal{A} \longrightarrow \mathcal{B}$ of associative
    $\ring{C}$-algebras with $\phi(a^*) = \phi(a)^*$.
\end{definition}
\begin{example}[$^*$-Algebras]
    \label{example:StarAlgebras}
    ~
    \begin{compactenum}
    \item Hermitian star products on Poisson manifolds give
        $^*$-algebras $(C^\infty(M)[[\lambda]], \star, \cc{})$ over
        $\ring{C} = \mathbb{C}[[\lambda]]$, where the $^*$-involution
        is the complex conjugation.
    \item For any pre-Hilbert space $\mathcal{H}$ the algebra
        $\mathfrak{B}(\mathcal{H})$ is a $^*$-algebra with
        $^*$-involution given by the adjoint. Moreover,
        $\mathfrak{F}(\mathcal{H}) \subseteq
        \mathfrak{B}(\mathcal{H})$ is a $^*$-ideal. 
    \item In particular, $M_n(\ring{C})$ is a $^*$-algebra.
    \item If $\mathcal{A}$, $\mathcal{B}$ are $^*$-algebras then
        $\mathcal{A} \otimes \mathcal{B}$ is again a $^*$-algebra with
        the obvious $^*$-involution.
    \item In particular, $M_n(\mathcal{A}) = \mathcal{A} \otimes
        M_n(\ring{C})$ is a $^*$-algebra.
    \end{compactenum}
\end{example}
Having a $^*$-algebra we can adapt the definitions of positive
functionals, positive algebra elements and positive maps from the
well-known theory of $^*$-algebras over $\mathbb{C}$, see e.g.
\cite{schmuedgen:1990a}, immediately to our algebraic context. This
motivates the following definitions \cite{bordemann.waldmann:1998a,
  bursztyn.waldmann:2003a:pre, bursztyn.waldmann:2001a}:
\begin{definition}[Positivity]
    \label{definition:PositiveStuff}
    Let $\mathcal{A}$ be a $^*$-algebra over $\ring{C}$.
    \begin{compactenum}
    \item A $\ring{C}$-linear functional $\omega: \mathcal{A}
        \longrightarrow \ring{C}$ is called positive if
        \begin{equation}
            \label{eq:PositiveFunctional}
            \omega(a^*a) \ge 0.
        \end{equation}
        If $\mathcal{A}$ is unital then $\omega$ is called a state if
        in addition $\omega(\Unit) = 1$.
    \item $a \in \mathcal{A}$ is called positive if $\omega(a) \ge 0$
        for all positive functionals $\omega$ of $\mathcal{A}$. We set
        \begin{equation}
            \label{eq:Aplus}
            \mathcal{A}^+ = \{a \in \mathcal{A} \; | \; a \;
            \textrm{is positive} \}
        \end{equation}
        \begin{equation}
            \label{eq:Aplusplus}
            \mathcal{A}^{++}
            = \left\{
                a \in \mathcal{A} 
                \; \left| \; 
                    a = \sum\nolimits_{i=1}^n \beta_i b_i^*b_i,
                    \; \textrm{with} \;
                    0 < \beta_i \in \ring{R}, b_i \in \mathcal{A}
                \right.
            \right\}.
        \end{equation}
    \item A linear map $\phi: \mathcal{A} \longrightarrow \mathcal{B}$
        into another $^*$-algebra $\mathcal{B}$ is called positive if
        $\phi(\mathcal{A}^+) \subseteq \mathcal{B}^+$. Moreover,
        $\phi$ is called completely positive if the componentwise
        extension $\phi: M_n(\mathcal{A}) \longrightarrow
        M_n(\mathcal{B})$ is positive for all $n \in \mathbb{N}$.
    \end{compactenum}
\end{definition}
\begin{remark}[Positive Elements and Maps]
    \label{remark:PositiveStuff}
    ~
    \begin{compactenum}
    \item Clearly, we have $\mathcal{A}^{++} \subseteq \mathcal{A}^+$
        but in general $\mathcal{A}^{++} \ne \mathcal{A}^+$.
    \item For $C^*$-algebras we have $\mathcal{A}^{++} =
        \mathcal{A}^+$ and, moreover, any positive element has a
        unique positive square root $a = (\sqrt{a})^2$. This follows
        from the spectral calculus.
    \item Any $^*$-homomorphism is a completely positive map.
    \item A linear map $\phi: \mathcal{A} \longrightarrow \mathcal{B}$
        is positive if and only if for any positive linear functional
        $\omega: \mathcal{B} \longrightarrow \ring{C}$ the pull-back
        $\phi^*\omega = \omega \circ \phi$ is a positive functional of
        $\mathcal{A}$. This is the case if and only if
        $\phi(\mathcal{A}^{++}) \subseteq \mathcal{B}^+$.
    \item A positive functional $\omega: \mathcal{A} \longrightarrow
        \ring{C}$ is a completely positive map. However, there are
        simple counter-examples which show that not every positive map
        is completely positive, even in the case $\ring{C} =
        \mathbb{C}$, see e.g.
        \cite[Exercise~11.\,5.\,15]{kadison.ringrose:1997b}.
    \end{compactenum}
\end{remark}
The following standard examples will be used later:
\begin{example}[Positive Maps]
    \label{example:PositiveElementsAndMaps}
    ~
    \begin{compactenum}
    \item The trace functional $\tr: M_n(\mathcal{A}) \longrightarrow
        \mathcal{A}$ is completely positive.
    \item The map $\tau: M_n(\mathcal{A}) \longrightarrow \mathcal{A}$
        defined by
        \begin{equation}
            \label{eq:tauComPos}
            \tau((a_{ij})) = \sum_{i,j=1}^n a_{ij}
        \end{equation}
        is completely positive.
    \item For $a \in \mathcal{A}^+$ and $b \in \mathcal{B}^+$ we have
        $a \otimes b \in (\mathcal{A} \otimes \mathcal{B})^+$.
        Indeed, for $b \in \mathcal{B}$ the map $a \mapsto a \otimes
        b^*b$ is clearly a positive map. Thus for $a \in
        \mathcal{A}^+$ the element $a \otimes b^*b$ is positive for
        all $b$. Hence also the map $b \mapsto a \otimes b$ is a
        positive map for positive $a$. Thus the claim follows. It also
        follows that the tensor product $\phi \otimes \psi$ of
        positive maps $\phi: \mathcal{A} \longrightarrow \mathcal{B}$,
        $\psi: \mathcal{C} \longrightarrow \mathcal{D}$ is positive.
        In particular, the tensor product of positive linear
        functionals is again a positive linear functional.
    \end{compactenum}
\end{example}

Though the above definition of positive functionals and elements is in
some sense the canonical one there are other concepts for positivity.
Indeed, in the theory of $O^*$-algebras the above definition does not
give the most useful concept, see the discussion in
\cite{schmuedgen:1990a}. In general, one defines a
\emph{$m$-admissible wedge} $\mathcal{K} \subseteq \mathcal{A}$ of a
unital $^*$-algebra to be a subset of Hermitian elements such that
$\mathcal{K}$ is closed under convex combinations, $a^*\mathcal{K}a
\subseteq \mathcal{K}$ for all $a \in \mathcal{A}$ and
$\mathcal{A}^{++} \subseteq \mathcal{K}$. Then the elements in
$\mathcal{K}$ are a replacement for the positive elements
$\mathcal{A}^+$. Also one can define a linear functional $\omega:
\mathcal{A} \longrightarrow \ring{C}$ to be \emph{positive with
  respect to} $\mathcal{K}$ if $\omega(\mathcal{K}) \ge 0$. In
particular, $\omega$ is positive in the usual sense, but not all
positive functionals will be positive with respect to $\mathcal{K}$.
Similarly, this gives a refined definition of (completely) positive
maps, leading to the notion of \emph{strong positivity} in the case of
$O^*$-algebras.

In the following we shall stick to the
Definition~\ref{definition:PositiveStuff} since it seems that for
deformation quantization this is the `correct' choice. At least in the
classical limit this definition produces the correct positive elements
in $C^\infty(M)$, see e.g. the discussion in
\cite[App.~B]{bursztyn.waldmann:2001a} and
\cite[Sect.~3]{waldmann:2003c:pre}.

%
%

\section{Examples of positive functionals in deformation quantization}
\label{sec:ExamplesPosFun}

We shall now discuss three basic examples of positive functionals in
deformation quantization: the $\delta$-functionals, the Schrödinger
functional and the positive traces and KMS functionals, see
\cite{bordemann.neumaier.waldmann:1999a,
  bordemann.neumaier.waldmann:1998a, waldmann:2000a,
  bordemann.waldmann:1998a} for these examples.

%
%

\subsection{The $\delta$-functional for the Weyl and Wick star
  product}
\label{subsec:DeltaFunctional}

We consider the most simple classical phase space $\mathbb{R}^{2n}$
with its standard symplectic structure and Poisson bracket. For this
example one knows several explicit formulas for star products
quantizing the canonical Poisson bracket. The most prominent one is
the \emph{Weyl-Moyal star product}
\begin{equation}
    \label{eq:WeylMoyalFormula}
    f \starweyl g = 
    \mu \circ \E^{\frac{\I\lambda}{2}\sum_{r=1}^n 
      \left(\frac{\partial}{\partial q^r} \otimes
          \frac{\partial}{\partial p_r}
          -
          \frac{\partial}{\partial p_r}
          \otimes \frac{\partial}{\partial q^r}
      \right)
    }
    f \otimes g,
\end{equation}
where $f, g \in C^\infty(\mathbb{R}^{2n})[[\lambda]]$ and $\mu(f
\otimes g) = fg$ is the pointwise (undeformed) product.

Consider the Hamiltonian of the isotropic harmonic oscillator $H(q,p)
= \frac{1}{2}(p^2 + q^2)$, where we put $m = \omega = 1$ for
simplicity.  Then we have
\begin{equation}
    \label{eq:HstarweylH}
    \cc{H} \starweyl H = H^2 - \frac{\lambda^2}{4}
\end{equation}
whence
\begin{equation}
    \label{eq:deltaNotPositive}
    \delta_0(\cc{H} \starweyl H) = - \frac{\lambda^2}{4} < 0.
\end{equation}
Thus the $\delta$-functional at $0$ (and similarly at any other point)
can \emph{not} be a positive functional for the Weyl-Moyal star
product, while classically all $\delta$-functionals are of course
positive. This has a very simple physical interpretation, namely that
points in phase space are no longer valid states in quantum theory: we
can not localize both space \emph{and} momentum coordinates because of
the uncertainly relations.

More interesting and in some sense surprising is the behaviour of the
\emph{Wick star product} (or normal ordered star product) which is
defined by
\begin{equation}
    \label{eq:WickStarproduct}
    f \starwick g = \mu \circ 
    \E^{2\lambda \sum_{r=1}^n
      \frac{\partial}{\partial z^r} 
      \otimes \frac{\partial}{\partial \cc{z}^r}
    }
    f \otimes g,
\end{equation}
where $z^r = q^r + \I p_r$ and $\cc{z}^r = q^r - \I p_r$. First recall
that $\starwick$ is equivalent to $\starweyl$ by the equivalence
transformation
\begin{equation}
    \label{eq:WeylToWick}
    f \starwick g = S\left( S^{-1} f \starweyl S^{-1} g\right),
\end{equation}
where
\begin{equation}
    \label{eq:SDef}
    S = \E^{\lambda \Delta}
    \quad
    \textrm{and}
    \quad
    \Delta = \sum_{r=1}^n \frac{\partial^2}{\partial z^r \partial \cc{z}^r}.
\end{equation}
The operator $\Delta$ is, up to a constant multiple, the Laplace
operator of the Euklidian metric on the phase space $\mathbb{R}^{2n}$.
With the explicit formula \eqref{eq:WickStarproduct} we find
\begin{equation}
    \label{eq:ccfstarwickf}
    \cc{f} \starwick f 
    = \sum_{r=0}^\infty \frac{(2\lambda)^r}{r!}
    \sum_{i_1, \ldots, i_r = 1}^n
    \left|
        \frac{\partial^r f}
        {\partial \cc{z}^{i_1} \cdots \partial \cc{z}^{i_r}}
    \right|^2
\end{equation}
for $f \in C^\infty(\mathbb{R}^{2n})[[\lambda]]$. Thus \emph{any}
classically positive functional of $C^\infty(\mathbb{R}^{2n})$ is also
positive with respect to the Wick star product. In particular, all the
$\delta$-functionals are positive. In some sense they have to be
interpreted as \emph{coherent states} in this context.

It should be remarked that this simple observation has quite drastic
consequences as we shall see in
Section~\ref{subsec:HermStarProductsComPos} which are far from being
obvious. We also remark that this result still holds for Wick-type
star products on arbitrary Kähler manifolds \cite{karabegov:1996a,
  bordemann.waldmann:1997a, bordemann.waldmann:1998a}.

%
%

\subsection{The Schrödinger functional}
\label{subsec:SchroedingerFunctional}

Consider again the Weyl-Moyal star product $\starweyl$ on
$\mathbb{R}^{2n}$ which we interprete now as the cotangent bundle
$\pi: T^*\mathbb{R}^n \longrightarrow \mathbb{R}^n$ of the
configuration space $\mathbb{R}^n$. Denote by
\begin{equation}
    \label{eq:ZeroSection}
    \iota: \mathbb{R}^n 
    \hookrightarrow T^*\mathbb{R}^n
\end{equation}
the zero section. We consider the following \emph{Schrödinger
  functional} $\omega: C^\infty_0(T^*\mathbb{R}^n)[[\lambda]]
\longrightarrow \mathbb{C}[[\lambda]]$ defined by
\begin{equation}
    \label{eq:Schroedingerfunctional}
    \omega(f) = \int_{\mathbb{R}^n} \iota^* f \D^n q.
\end{equation}
Thanks to the restriction to formal series of functions with compact
support the integration with respect to the usual Lebesgues measure
$\D^n q$ is well-defined. One defines the operator
\begin{equation}
    \label{eq:Neumaier}
    N = \E^{\frac{\lambda}{2\I} \Delta},
    \quad
    \textrm{where}
    \quad
    \Delta = \sum_{k=1}^n 
    \frac{\partial^2}{\partial p_k \partial q^k}
\end{equation}
is now the Laplacian (`d'Alembertian') with respect to the maximally
indefinite metric obtained by pairing the configuration space
variables with the momentum variables. In fact, $N$ is the equivalence
transformation between the Weyl-Moyal star product and the
standard-ordered star product as we shall see later. By some
successive integration by parts one finds
\begin{equation}
    \label{eq:SchroedingerOnProducts}
    \omega(\cc{f} \starweyl g) 
    = \int_{\mathbb{R}^n} \cc{(\iota^*Nf)} (\iota^*Ng) \D^nq,
\end{equation}
whence immediately
\begin{equation}
    \label{eq:SchroedingerPositive}
    \omega(\cc{f} \starweyl f) \ge 0.
\end{equation}
Thus the Schrödinger functional is a positive functional with respect
to the Weyl-Moyal star product.

In fact, there is a geometric generalization for any cotangent bundle
$\pi: T^*Q \longrightarrow Q$ of this construction,
see~\cite{bordemann.neumaier.pflaum.waldmann:2003a,
  bordemann.neumaier.waldmann:1999a,
  bordemann.neumaier.waldmann:1998a, pflaum:2000a, pflaum:1998c}: For
a given torsion-free connection $\nabla$ on the configuration space
$Q$ and a given positive smooth density $\mu$ on $Q$ one can construct
rather explicitly a star product $\starweyl$ which is the direct
analog of the usual Weyl-Moyal star product in flat space. Moreover,
using the connection $\nabla$ one obtains a maximally indefinite
pseudo-Riemannian metric on $T^*Q$ coming from the natural pairing of
the vertical spaces with the horizontal space. The Laplacian $\Delta$
(d'Alembertian) of this indefinite metric is in a bundle chart locally
given by the explicit formula
\begin{equation}
    \label{eq:DeltaGeometric}
    \Delta = 
    \sum_{k=1}^n \frac{\partial^2}{\partial q^k \partial p_k}
    +
    \sum_{k, \ell, m=1}^n p_k \pi^*\Gamma^k_{\ell m} 
    \frac{\partial^2}{\partial p_\ell \partial p_m}
    +
    \sum_{k, \ell=1}^n 
    \pi^*\Gamma^k_{k\ell} \frac{\partial}{\partial p_\ell},
\end{equation}
generalizing \eqref{eq:Neumaier} to the general curved framework. Here
$\Gamma^k_{\ell m}$ denote the Christoffel symbols of the connection
$\nabla$. This gives a geometric version of the operator $N$
(Neumaier's operator)
\begin{equation}
    \label{eq:NeumaierGeometrisch}
    N = \E^{\frac{\lambda}{2\I}(\Delta + \mathsf{F}(\alpha))},
\end{equation}
where $\alpha \in \Gamma^\infty(T^*Q)$ is the one-form determined by
$\nabla_X\mu = \alpha(X)\mu$ and $\mathsf{F}(\alpha)$ is the
differential operator 
\begin{equation}
    \label{eq:Falpha}
    (\mathsf{F}(\alpha) f)(\alpha_q) = \frac{\D}{\D t}\Big|_{t=0}
    f(\alpha_q + t\alpha(q)),
\end{equation}
where $q \in Q$ and $\alpha_q \in T^*_qQ$.  In particular, one has
$\alpha = 0$ if the density $\mu$ is covariantly constant. This is the
case if we choose $\nabla$ to be the Levi-Civita connection of a
Riemannian metric and $\mu = \mu_g$ to be the corresponding Riemannian
volume density.  Note that in a typical Hamiltonian system on $T^*Q$
we have a kinetic energy in the Hamiltonian which is nothing else than
a Riemannian metric. Thus there is a prefered choice in this
situation.

The Schrödinger functional in this context is simply given by the
integration with respect to the a priori chosen density $\mu$
\begin{equation}
    \label{eq:SchroedingerGeometric}
    \omega(f) = \int_Q \iota^*f \; \mu,
\end{equation}
where again we restrict to $f \in C^\infty_0(T^*Q)[[\lambda]]$ to have
a well-defined integration. Now the non-trivial result is that the
above formulas still hold in this general situation. We have
\begin{equation}
    \label{eq:SchroedingerGeometricfg}
    \omega(\cc{f} \starweyl g) 
    = \int_Q \cc{(\iota^*Nf)} (\iota^*Ng) \; \mu
\end{equation}
whence the Schrödinger functional is positive
\begin{equation}
    \label{eq:SchroedingerPositiveGeometric}
    \omega(\cc{f} \starweyl f) 
    = \int_Q \cc{(\iota^*Nf)} (\iota^*Nf) \; \mu \ge 0.
\end{equation}
The proof consists again in a successive integration by parts which is
now much more involved due to the curvature terms coming from
$\nabla$, see \cite{bordemann.neumaier.waldmann:1998a, neumaier:1998a,
  bordemann.neumaier.waldmann:1999a} for details.

%
%

\subsection{Positive traces and KMS functionals}
\label{subsec:PositiveTraces}

We consider a connected symplectic manifold $(M, \omega)$ with a
Hermitian star product $\star$. Then it is well-known that there
exists a unique \emph{trace functional} up to normalization and even
the normalization can be chosen in a canonical way
\cite{nest.tsygan:1995a, gutt.rawnsley:2002a, fedosov:2002a,
  karabegov:1998b}. Here a trace functional means a
$\mathbb{C}[[\lambda]]$-linear functional
\begin{equation}
    \label{eq:TraceFunctional}
    \tr: C^\infty_0(M)[[\lambda]] 
    \longrightarrow \mathbb{C}[[\lambda]]
\end{equation}
such that
\begin{equation}
    \label{eq:TraceProperty}
    \tr(f \star g) = \tr(g \star f).
\end{equation}
Furthermore, it is known that $\tr$ is of the form
\begin{equation}
    \label{eq:TraceLowestOrder}
    \tr(f) = c_0 \int_M f \; \Omega 
    + \; \textrm{higher orders in} \; \lambda,
\end{equation}
where $\Omega$ is the Liouville form on $M$ and $c_0$ is a
normalization constant. If one does not need the higher order
corrections then the star product is called \emph{strongly closed}
\cite{connes.flato.sternheimer:1992a}.

Since $\star$ is a Hermitian star product, the functional $f \mapsto
\cc{\tr(\cc{f})}$ is still a trace whence we can assume that the trace
we started with is already a \emph{real} trace. In particular, for
this choice $c_0 = \cc{c_0}$ is real. Passing to $-\tr$ if necessary
we can assume that $c_0 > 0$. Then
\begin{equation}
    \label{eq:TracePositive}
    \tr(\cc{f} \starweyl f)
    = c_0 \int_M \cc{f}f  \; \Omega
    + \; \textrm{higher orders in} \; \lambda.
\end{equation}
Hence, if $f \ne 0$, already the zeroth order in $\tr(\cc{f} \star f)$
is nonzero and clearly positive. Thus by definition of the ordering of
$\mathbb{R}[[\lambda]]$ we see that $\tr$ is a positive functional.
Note however the difference in the argument compared to the
$\delta$-functional.

More generally, we can consider \emph{thermodynamical states}, i.e.
\emph{KMS functionals}. Here we fix a Hamiltonian $H = \cc{H} \in
C^\infty(M)[[\lambda]]$ and an `inverse temperatur' $\beta > 0$.
Then the \emph{star exponential} $\Exp(-\beta H) \in
C^\infty(M)[[\lambda]]$ is well-defined as the solution of the
differential equation
\begin{equation}
    \label{eq:StarExpDef}
    \frac{\D}{\D \beta} \Exp(-\beta H) = - H \star \Exp(- \beta H)
\end{equation}
with initial condition $\Exp(0) = 1$. The star exponential has all
desired functional properties of a `exponential function', see e.g.
the discussion in \cite{bordemann.roemer.waldmann:1998a}. Using this,
the \emph{KMS functional} corresponding to this data is defined by
\begin{equation}
    \label{eq:KMSFunctional}
    \omega_{H, \beta}(f) = \tr(\Exp(-\beta H) \star f)
\end{equation}
for $f \in C^\infty_0(M)[[\lambda]]$. The positivity of $\tr$ and the
existence of a square root $\Exp(-\frac{\beta}{2}H)$ of $\Exp(-\beta
H)$ shows that the KMS functional is indeed a positive functional
again.
\begin{remark}
    \label{remark:KMS}
    Originally, KMS functionals are characterized by the so-called KMS
    condition \cite{kubo:1957a, martin.schwinger:1959a,
      haag.hugenholtz.winnink:1967a} in a more operator-algebraic
    approach which was transfered to deformation quantization in
    \cite{basart.lichnerowicz:1985a,
      basart.flato.lichnerowicz.sternheimer:1984a}. However, the
    existence of a unique trace functional allows to classify the KMS
    functionals completely yielding the above characterization
    \cite{bordemann.roemer.waldmann:1998a}. Note that this only holds
    in the (connected) symplectic framework as in the general Poisson
    framework traces are no longer unique, see
    e.g.~\cite{felder.shoikhet:2000a,
      bieliavsky.bordemann.gutt.waldmann:2003a}, so the arguments of
    \cite{bordemann.roemer.waldmann:1998a} do no longer apply. Thus it
    would be very interesting to get some more insight in the nature
    of KMS functionals in the general Poisson case.
\end{remark}

%
%

\section{Deformation and classical limit of positive functionals}
\label{sec:ClassLimPosFun}

The interpretation of states as positive linear functionals allows for
a simple definition of a classical limit of a state of a Hermitian
deformation of a $^*$-algebra. The converse question whether any
`classical' state can be deformed is much more delicate and shines
some new light on the relevance of the Wick star product.

%
%

\subsection{Completely positive deformations}
\label{subsec:CompletelyPositiveDeformations}

Let $\mathcal{A}$ be a $^*$-algebra over $\ring{C}$. In the spirit of
star products we consider a \emph{Hermitian deformation} of
$\mathcal{A}$, i.e. an associative $\ring{C}[[\lambda]]$-linear
multiplication $\star$ for $\mathcal{A}[[\lambda]]$ making
$\boldsymbol{\mathcal{A}} = (\mathcal{A}[[\lambda]], \star, {}^*)$ a
$^*$-algebra over $\ring{C}[[\lambda]]$. The
$\ring{C}[[\lambda]]$-bilinearity of $\star$ implies that
\begin{equation}
    \label{eq:astarb}
    a \star b = \sum_{r=0}^\infty \lambda^r C_r(a, b)
\end{equation}
with $\ring{C}$-bilinear maps $C_r: \mathcal{A} \otimes \mathcal{A}
\longrightarrow \mathcal{A}$, extended to $\boldsymbol{\mathcal{A}}$
by the usual $\ring{C}[[\lambda]]$-bilinearity. As usual, the
deformation aspect is encoded in the condition $C_0(a, b) = ab$. Note
that we do not want to deform the $^*$-involution though in principle
this can also be taken into account, see e.g. the discussion in
\cite{bursztyn.waldmann:2001a, bursztyn.waldmann:2000a}. As we
mentioned already in Example~\ref{example:OrderedRings} we are still
in the framework of $^*$-algebras over ordered rings as
$\ring{R}[[\lambda]]$ is canonically ordered.

Now assume $\boldsymbol{\omega}: \boldsymbol{\mathcal{A}}
\longrightarrow \ring{C}[[\lambda]]$ is a $\ring{C}[[\lambda]]$-linear
positive functional. Then the $\ring{C}[[\lambda]]$-linearity implies
that $\boldsymbol{\omega}$ is actually of the form
$\boldsymbol{\omega} = \sum_{r=0}^\infty \lambda^r \omega_r$ with
$\ring{C}$-linear functionals $\omega_r: \mathcal{A} \longrightarrow
\ring{C}$, the later being canonically extended to
$\boldsymbol{\mathcal{A}}$ by $\ring{C}[[\lambda]]$-linearity. Since
$\star$ deforms the given multiplication of $\mathcal{A}$ we obtain
from the positivity of $\boldsymbol{\omega}$
\begin{equation}
    \label{eq:BoldOmegaPos}
    0 \le 
    \boldsymbol{\omega}(a^* \star a)
    = \omega_0(a^*a) 
    + \lambda \left(
        \omega_0(C_1(a^*, a)) + \omega_1(a^*a)
    \right)
    + \textrm{higher order terms}.
\end{equation}
Thus it follows immediately from the ordering of $\ring{R}[[\lambda]]$
that $\omega_0$ has to be a \emph{positive} linear functional of
$\mathcal{A}$. In this sense, \emph{the classical limit of a quantum
  state is a classical state}. Note that this statement is
non-trivial, though physically of course more than plausible.

This observation immediately raises the question whether the converse
is true as well: can we always deform states? We know already from the
example of the $\delta$-functional and the Weyl-Moyal star product
that in general some quantum corrections $\omega_1$, $\omega_2$,
\ldots are unavoidable. The reason can easily be seen from the
expansion \eqref{eq:BoldOmegaPos}: If $\omega_0(a^*a) = 0$ then the
positivity of $\boldsymbol{\omega}$ is decided in the next order. But
this involves now the higher order terms $C_1(a^*, a)$, etc. of the
deformed product and these terms usually do not have any reasonable
positivity properties. Thus the terms $\omega_1$ etc. have to be
chosen well in order to compensate this. We state the following
definition \cite{bursztyn.waldmann:2000a,bursztyn.waldmann:2003a:pre}:
\begin{definition}[Positive Deformations]
    \label{definition:PositiveDeformation}
    A Hermitian deformation $\boldsymbol{\mathcal{A}} =
    (\mathcal{A}[[\lambda]], \star, {}^*)$ of a $^*$-algebra
    $\mathcal{A}$ over $\ring{C}$ is called a positive deformation if
    for any positive $\ring{C}$-linear functional there exists a
    deformation $\boldsymbol{\omega} = \sum_{r=0}^\infty \lambda^r
    \omega_r: \boldsymbol{\mathcal{A}} \longrightarrow
    \ring{C}[[\lambda]]$ into a $\ring{C}[[\lambda]]$-linear positive
    functional with respect to $\star$. Furthermore
    $\boldsymbol{\mathcal{A}}$ is called a completely positive
    deformation if $M_n(\boldsymbol{\mathcal{A}})$ is a positive
    deformation of $M_n(\mathcal{A})$ for all $n \in \mathbb{N}$.
\end{definition}
\begin{example}[A non-positive deformation]
    \label{example:NonPositiveDeformation}
    Let $\mathcal{A}$ be a $^*$-algebra over $\ring{C}$ with
    multiplication $\mu: \mathcal{A} \otimes \mathcal{A}
    \longrightarrow \mathcal{A}$. Then $\boldsymbol{\mathcal{A}} =
    (\mathcal{A}[[\lambda]], \lambda\mu, {}^*)$ is a Hermitian
    deformation of the trivial $^*$-algebra $\mathcal{A}_0$ being
    $\mathcal{A}$ as $\ring{C}$-module and equipped with the zero
    multiplication. Since for this trivial $^*$-algebra any (real)
    functional is positive, we can not expect to have a positive
    deformation in general since positivity with respect to $\star
    =\lambda \mu$ is a non-trivial condition. Thus not any Hermitian
    deformation is a positive deformation showing the non-triviality
    of the definition.
\end{example}

%
%

\subsection{Complete positivity of Hermitian star products}
\label{subsec:HermStarProductsComPos}

Recall that for the Wick star product $\starwick$ on $\mathbb{R}^{2n}
= \mathbb{C}^n$ we do not need to deform classically positive linear
functionals at all: They are positive with respect to $\starwick$ for
free.

Thus if $\star$ is any Hermitian symplectic star product on
$\mathbb{R}^{2n}$ it is equivalent and hence $^*$-equivalent to
$\starwick$ by some $^*$-equivalence transformation $T = \id +
\sum_{r=1}^\infty \lambda^r T_r$. Hence for any classically positive
linear functional $\omega_0$ we see that
\begin{equation}
    \label{eq:DeformClassPosFun}
    \omega = \omega_0 \circ T 
    = \omega_0 + \lambda \omega_0 \circ T_1 +
    \textrm{higher order terms}
\end{equation}
gives a positive linear functional with respect to $\star$ deforming
$\omega_0$.
\begin{example}
    \label{example:PositiveDelta}
    The functional $f \mapsto \delta \circ \E^{\lambda \Delta} f$ is a
    deformation of the $\delta$-functional into a positive linear
    functional with respect to the Weyl-Moyal star product.
\end{example}

Now let $\star$ be a Hermitian star product on a symplectic manifold
$M$. Using a quadratic partition of unity $\sum_\alpha
\cc{\chi}_\alpha \chi_\alpha = 1$ subordinate to some Darboux atlas of
$M$ we can localize a given classically positive linear functional
$\omega_0$ by writing
\begin{equation}
    \label{eq:OmegaLocalChialpha}
    \omega(f) = \sum_{\alpha} \omega_0 
    \left( 
        \cc{\chi}_\alpha \star f \star \chi_\alpha
    \right)
\end{equation}
such that each $\omega_0(\cc{\chi}_\alpha \star \: \cdot \: \star
\chi_\alpha)$ has now support in one Darboux chart. Hence we can
replace $\omega_0$ by some deformation (depending on $\alpha$) in
order to make it positive with respect to $\star$ and glue things
together in the end. The finaly result will then be a deformation of
$\omega_0$ which is now positive with respect to $\star$. This shows
\cite{bursztyn.waldmann:2000a}:
\begin{theorem}
    \label{theorem:SymplPosDef}
    Any symplectic Hermitian star product is a positive deformation.
\end{theorem}
It is easy to see that the same argument applies for
$M_n(\mathbb{C})$-valued functions whence any symplectic Hermitian
star product is even a \emph{completely} positive deformation.

The Poisson case proves to be more involved. First we note that we
again have to solve only the local problem of showing that Hermitian
star product in $\mathbb{R}^n$ are positive since the same glueing as
in the symplectic case can be applied. Hence we consider the local
case $\mathbb{R}^n$ with coordinates $q^1, \ldots, q^n$ equipped with
a Hermitian star product $\star$ deforming some Poisson structure.  We
define
\begin{equation}
    \label{eq:FormalWeylAlg}
    \mathcal{W}_0 = C^\infty(\mathbb{R}^n)[[p_1, \ldots, p_n]]
    \quad
    \textrm{and}
    \quad
    \mathcal{W} = \mathcal{W}_0[[\lambda]],
\end{equation}
whence we can equipp $\mathcal{W}$ with the Weyl-Moyal star product
$\starweyl$.  Clearly the formula \eqref{eq:WeylMoyalFormula} extends
to the `formal' momentum variables. We can think of $\mathcal{W}_0$ as
the functions on a `formal cotangent bundle' of $\mathbb{R}^n$. Thus
we have also the two canonical maps $\iota^*: \mathcal{W}_0
\longrightarrow C^\infty(\mathbb{R}^n)$ and $\pi^*:
C^\infty(\mathbb{R}^n) \hookrightarrow \mathcal{W}_0$ which are
algebra homomorphisms with respect to the undeformed products.

Now the idea is to deform $\pi^*$ into a $^*$-algebra homomorphism
\begin{equation}
    \label{eq:TauQFSR}
    \tau = \sum_{r=0}^\infty \tau_k: 
    \left(C^\infty(\mathbb{R}^n)[[\lambda]], \star\right)
    \longrightarrow
    \left(\mathcal{W}, \starweyl\right),
\end{equation}
where $\tau_0 = \pi^*$ and $\tau_k$ is homogeneous of degree $k$ with
respect to the grading induced by the $\lambda$-Euler derivation
(homogeneity operator)
\begin{equation}
    \label{eq:Homogen}
    \mathsf{H} = \lambda \frac{\partial}{\partial \lambda} 
    + \sum_{k=1}^n p_k \frac{\partial}{\partial p_k}.
\end{equation}
This can actually be done by an inductive construction of the $\tau_k$
using the vanishing of a certain Hochschild cohomology (essentially
$\mathcal{W}$ as a $C^\infty(\mathbb{R}^n)$-bimodule with bimodule
multiplication given by $\starweyl$), see
\cite{bordemann.neumaier.nowak.waldmann:1997a:misc, nowak:1997a} as
well as \cite{bursztyn.waldmann:2003a:pre} for details.
\begin{remark}[Quantized Formal Symplectic Realization]
    \label{remark:DefQuantSymplReal}
    Note that by setting $\lambda = 0$ in $\tau$ one obtains a
    `formal' symplectic realization of the Poisson manifold
    $\mathbb{R}^n$. Here formal is understood in the sense that the
    dependence on the momentum variables is formal. Thus $\tau$ can be
    seen as a \emph{quantized formal symplectic realization}. Note
    also that the same construction can be done for any Poisson
    manifold $Q$ if one replaces $\mathcal{W}_0$ by the formal
    functions on the cotangent bundle $T^*Q$ and $\starweyl$ by the
    homogeneous star product of Weyl type for $T^*Q$ constructed out
    of a connection on $Q$ as in
    \cite{bordemann.neumaier.waldmann:1999a,
      bordemann.neumaier.waldmann:1998a}.
\end{remark}

Having $\tau$ it is very easy to obtain a deformation of a given
positive functional $\omega_0: C^\infty(\mathbb{R}^n) \longrightarrow
\mathbb{C}$. The key observation is that $\omega_0 \circ \iota^*$ is a
positive linear functional of $\mathcal{W}_0$ and thus, using the Wick
star product again, $\omega_0 \circ \iota^* \circ S$ is a positive
$\mathbb{C}[[\lambda]]$-linear functional of $\mathcal{W}$, equipped
with the Weyl-Moyal star product $\starweyl$, where $S$ is defined as
in \eqref{eq:SDef} only using the formal momentum variables for the
definition of the partial derivatives $\frac{\partial}{\partial z^k}$
and $\frac{\partial}{\partial \cc{z}^k}$. Then clearly $\omega_0 \circ
\iota^* \circ S \circ \tau$ gives the deformation of $\omega_0$ into a
positive functional with respect to $\star$. It is also clear that the
matrix-valued case causes no further problems whence we have the
following result \cite{bursztyn.waldmann:2003a:pre}, answering thereby
a question raised in \cite{bieliavsky.bordemann.gutt.waldmann:2003a}:
\begin{theorem}
    \label{theorem:ComPosDef}
    Every Hermitian star product on a Poisson manifold is a completely
    positive deformation.
\end{theorem}

We discuss now some easy consequences of this result:
\begin{corollary}
    \label{corollary:SuffManyPosFun}
    For Hermitian star products one has sufficiently many positive
    linear functionals in the sense that for $0 \ne H = \cc{H} \in
    C^\infty(M)[[\lambda]]$ one finds a positive linear functional
    $\omega: C^\infty(M)[[\lambda]] \longrightarrow
    \mathbb{C}[[\lambda]]$ with $\omega(H) \ne 0$.
\end{corollary}
Indeed, classically this is certainly true and by
Theorem~\ref{theorem:ComPosDef} we just have to deform an appropriate
classically positive functional.
\begin{corollary}
    \label{corollary:ClassLimPosEl}
    If $H \in (C^\infty(M)[[\lambda]], \star)^+$ then for the
    classical limit we have $H_0 \in C^\infty(M)^+$ as well.
\end{corollary}
In general, it is rather difficult to characterize the positive
algebra elements in a star product algebra beyond this zeroth order.

A nice application is obtained for the following situation: consider a
Lie algebra $\mathfrak{g}$ of a compact Lie group and let $\starBCH$
be the Baker-Campbell-Hausdorff star product on $\mathfrak{g}^*$, as
constructed in \cite{gutt:1983a}. Then any $\mathfrak{g}$-invariant
functional on $C^\infty(\mathfrak{g}^*)$ defines a trace with respect
to $\starBCH$, see \cite{bieliavsky.bordemann.gutt.waldmann:2003a}.
The question is whether one can deform a classically positive trace
into a positive trace with respect to $\starBCH$. In
\cite{bieliavsky.bordemann.gutt.waldmann:2003a} this was obtained for
very particular trace functionals by some BRST like quantization of a
phase space reduction. However, as already indicated in
\cite[Sect.~8]{bieliavsky.bordemann.gutt.waldmann:2003a}, we can just
deform the trace in some positive functional thanks to
Theorem~\ref{theorem:ComPosDef}, loosing probably the trace property.
But averaging over the \emph{compact} group corresponding to
$\mathfrak{g}$ gives again a $\mathfrak{g}$-invariant functional,
hence a trace, now \emph{without} loosing the positivity. Thus we have
the result:
\begin{corollary}
    \label{corollary:BCHStarProduct}
    Any $\mathfrak{g}$-invariant functional on
    $C^\infty(\mathfrak{g}^*)$ can be deformed into a positive trace
    functional with respect to the Baker-Campbell-Hausdorff star
    product $\starBCH$.
\end{corollary}




%
%

\section{The GNS construction and examples}
\label{sec:GNS}

As we have argued in Section~\ref{subsec:SuperpositionSelections} we
have to construct $^*$-representations of the observable algebra in
order to implement the super-position principle for states. The GNS
construction, which is well-known for $^*$-algebras over $\mathbb{C}$,
see e.g.~\cite{kadison.ringrose:1997a, kadison.ringrose:1997b,
  sakai:1971a, bratteli.robinson:1987a, schmuedgen:1990a}, provides a
canonical way to construct such $^*$-representations out of a given
positive linear functional.

%
%

\subsection{$^*$-Representation theory}
\label{subsec:StarRepTheory}

We start with some general remarks on $^*$-representations of
$^*$-algebras transfering the usual notions to the algebraic framework
for $^*$-algebras over $\ring{C}$.
\begin{definition}[$^*$-Representations]
    \label{definition:StarRepresentations}
    Let $\mathcal{A}$ be a $^*$-algebra over $\ring{C}$.
    \begin{compactenum}
    \item A $^*$-representation $\pi$ of $\mathcal{A}$ on a
        pre-Hilbert space $\mathcal{H}$ is a $^*$-homomorphism
        \begin{equation}
            \label{eq:StarRep}
            \pi: \mathcal{A} 
            \longrightarrow \mathfrak{B}(\mathcal{H}).
        \end{equation}
    \item An intertwiner $T$ between two $^*$-representations
        $(\mathcal{H}_1, \pi_1)$ and $(\mathcal{H}_2, \pi_2)$ is a map
        $T \in \mathfrak{B}(\mathcal{H}_1, \mathcal{H}_2)$ with
        \begin{equation}
            \label{eq:Intertwiner}
            T \pi_1(a) = \pi_2(a) T
        \end{equation}
        for all $a \in \mathcal{A}$.
    \item Two $^*$-representations are called unitarily equivalent if
        there exists a unitary intertwiner between them.
    \item A $^*$-representation $(\mathcal{H}, \pi)$ is called
        strongly non-degenerate if $\pi(\mathcal{A}) \mathcal{H} =
        \mathcal{H}$.
    \item A $^*$-representation $(\mathcal{H}, \pi)$ is called cyclic
        with cyclic vector $\Omega \in \mathcal{H}$ if
        $\pi(\mathcal{A})\Omega = \mathcal{H}$.
    \end{compactenum}
\end{definition}
\begin{remark}[$^*$-Representations]
    \label{remark:StarRepresentations}
    ~
    \begin{compactenum}
    \item For unital $^*$-algebras we only consider unital
        $^*$-homomorphisms by convention. Thus in the unital case,
        $^*$-representations are always strongly non-degenerate by
        convention. This is reasonable since if $\pi$ is a
        $^*$-representation of a unital $^*$-algebra then $\pi(\Unit)
        = P$ is a projection $P = P^* = P^2$ and thus we can split the
        representation space $\mathcal{H} = P\mathcal{H} \oplus
        (\id-P)\mathcal{H}$ into an orthogonal direct sum. Then the
        $^*$-representation $\pi$ is easily seen to be block-diagonal
        with respect to this decomposition and $\pi$ is identically
        zero on $(\id-P)\mathcal{H}$. Thus the only `interesting' part
        is $P\mathcal{H}$ which is strongly non-degenerate.
    \item The space of intertwiners is a $\ring{C}$-module and clearly
        the composition of intertwiners gives again an
        intertwiner.
    \end{compactenum}
\end{remark}
This last observation allows to state the following definition of the
category of $^*$-representations:
\begin{definition}[$^*$-Representation theory]
    \label{definition:RepTheory}
    The $^*$-representation theory of $\mathcal{A}$ is the category of
    $^*$-representations with intertwiners as morphisms. It is denoted
    by $\rep(\mathcal{A})$, and $\Rep(\mathcal{A})$ denotes the
    subcategory of strongly non-degenerate $^*$-representations.
\end{definition}
Thus the final goal is to understand $\Rep(\mathcal{A})$ for a given
$^*$-algebra $\mathcal{A}$ like e.g.  a star product algebra
$\mathcal{A} = (C^\infty(M)[[\lambda]], \star, \cc{\phantom{a}})$.
Since direct orthogonal sums of $^*$-representations give again
$^*$-representations one would like to understand if and how a given
$^*$-representation can be decomposed into a direct sum of
non-decomposable (irreducible) $^*$-representations. In practice this
will be only achievable for very particular and simple examples. In
general, it is rather hopeless to obtain such a complete picture of
$\Rep(\mathcal{A})$. Thus other strategies have to be developed, like
e.g. finding `interesting' subclasses of $^*$-rep\-re\-sen\-ta\-tions.

Nevertheless, the GNS construction will allow to construct at least a
large class of $^*$-re\-pre\-sen\-ta\-tions, namely out of given
positive functionals. Hence one can discuss those GNS representations
which come from positive functionals that are of particular interest.

%
%

\subsection{The general GNS construction}
\label{subsec:GeneralGNS}

The whole GNS construction is a consequence of the Cauchy-Schwarz
inequality which itself is obtained from the following simple but
crucial lemma. Of course this is well-known for the case $\ring{C} =
\mathbb{C}$.
\begin{lemma}
    \label{lemma:Babylonian}
    Let $p(z, w) = a \cc{z}z + b z \cc{w} + b' \cc{z} w + c \cc{w}w
    \ge 0$ for all $z, w \in \ring{C}$, where $a, b, b', c \in
    \ring{C}$. Then
    \begin{equation}
        \label{eq:Babylonian}
        a \ge 0,\quad
        c \ge 0, \quad
        \cc{b} = b'
        \quad
        \textrm{and}
        \quad
        b\cc{b} \le ac.
    \end{equation}
\end{lemma}
\begin{proof}
    Taking $z = 0$ gives $c \ge 0$ and similarly $a \ge 0$ follows.
    Taking $z=1, \I$ and $w = 1$ implies $\cc{b} = b'$.  Now we first
    consider the case $a = 0 = c$.  Then taking $w = \cc{b}$ gives $b
    \cc{b} (z + \cc{z}) \ge 0$ for all $z \in \ring{C}$, whence $z =
    -1$ gives $b = 0$.  Thus we can assume, say, $a > 0$. Taking $z =
    \cc{b}$ and $w = -a$ gives $a (\cc{b}b - \cc{b} b - \cc{b} b + ac)
    \ge 0$ whence the inequality $\cc{b} b \le ac$ holds, too.
\end{proof}
\begin{corollary}[Cauchy-Schwarz inequality]
    \label{corollary:CauchySchwarz}
    Let $\omega: \mathcal{A} \longrightarrow \ring{C}$ be a positive
    linear functional. Then
    \begin{equation}
        \label{eq:omegaReal}
        \omega(a^*b) = \cc{\omega(b^*a)}
    \end{equation}
    and
    \begin{equation}
        \label{eq:CSU}
        \omega(a^*b) \cc{\omega(a^*b)} 
        \le \omega(a^*a) \omega(b^*b)
    \end{equation}
    for all $a, b \in \mathcal{A}$.
\end{corollary}
For the proof one considers $p(z, w) = \omega((za+wb)^*(za+wb)) \ge
0$.  In particular, if $\mathcal{A}$ is unital then we have
\begin{equation}
    \label{eq:OmegaReal}
    \omega(a^*) = \cc{\omega(a)}
\end{equation}
and $\omega(\Unit) = 0$ implies already $\omega =0$. Thus (by passing
to the quotient field if necessary) we can replace positive linear
functionals by normalized ones, i.e. by states $\omega(\Unit) = 1$.

Now we consider the following subset
\begin{equation}
    \label{eq:GelfandIdeal}
    \mathcal{J}_\omega 
    = \{a \in \mathcal{A} \; | \; \omega(a^*a) = 0\}
    \subseteq \mathcal{A}.
\end{equation}
Using the Cauchy-Schwarz inequality one obtains immediately that
$\mathcal{J}_\omega$ is a left ideal of $\mathcal{A}$, the so-called
\emph{Gel'fand ideal} of $\omega$.  Thus we can form the quotient
\begin{equation}
    \label{eq:Homega}
    \mathcal{H}_\omega = \mathcal{A} \big/ \mathcal{J}_\omega,
\end{equation}
which is a left $\mathcal{A}$-module in the usual way. We denote by
$\psi_b \in \mathcal{H}_\omega$ the equivalence class of $b \in
\mathcal{A}$. Then the left module structure can be written as
\begin{equation}
    \label{eq:GNSpiDef}
    \pi_\omega(a) \psi_b = \psi_{ab},
\end{equation}
for $a \in \mathcal{A}$ and $\psi_b \in
\mathcal{H}_\omega$. Furthermore, $\mathcal{H}_\omega$ becomes a
pre-Hilbert space by setting
\begin{equation}
    \label{eq:GNSInnerProduct}
    \SP{\psi_b, \psi_c}_\omega = \omega(b^*c).
\end{equation}
Indeed, this is well-defined thanks to the Cauchy-Schwarz inequality.
Moreover, $\SP{\cdot,\cdot}_\omega$ is positive definite since we
divided precisely by the `null-vectors' in $\mathcal{A}$. Finally, we
have
\begin{equation}
    \label{eq:GNSStarRep}
    \SP{\psi_b, \pi_\omega(a)\psi_c}_\omega
    = \omega(b^*ac)
    = \omega((a^*b)^*c)
    = \SP{\pi_\omega(a^*)\psi_b, \psi_c}_\omega,
\end{equation}
whence $\pi_\omega$ is a $^*$-representation of $\mathcal{A}$ on
$\mathcal{H}_\omega$. 
\begin{definition}[GNS representation]
    \label{definition:GNS}
    For a positive linear functional $\omega: \mathcal{A}
    \longrightarrow \ring{C}$ the $^*$-representation
    $(\mathcal{H}_\omega, \pi_\omega)$ is called the GNS
    representation of $\omega$.
\end{definition}

If $\mathcal{A}$ is unital, we can recover the functional $\omega$ as
`vacuum expectation value' with respect to the `vacuum vector'
$\psi_\Unit$ as follows
\begin{equation}
    \label{eq:VacuumGNS}
    \omega(a) = \SP{\psi_\Unit, \pi_\omega(a) \psi_\Unit}_\omega,
\end{equation}
and $\psi_\Unit$ is obviously a cyclic vector for the GNS
representation since
\begin{equation}
    \label{eq:psiUnitCyclic}
    \psi_b = \pi_\omega(b) \psi_\Unit
\end{equation}
for all $\psi_b \in \mathcal{H}_\omega$. It turns out that these
proporties already characterize the $^*$-representation
$(\mathcal{H}_\omega, \pi_\omega, \psi_\Unit)$ up to unitary
equivalence:
\begin{theorem}[GNS Representation]
    \label{theorem:GNSUnique}
    Let $\mathcal{A}$ be unital and let $\omega: \mathcal{A}
    \longrightarrow \ring{C}$ be a positive linear functional. If
    $(\mathcal{H}, \pi, \Omega)$ is a cyclic $^*$-representation such
    that
    \begin{equation}
        \label{eq:HisGNS}
        \omega(a) = \SP{\Omega, \pi(a)\Omega}
    \end{equation}
    for all $a \in \mathcal{A}$ then $(\mathcal{H}, \pi, \Omega)$ is
    unitarily equivalent to the GNS representation
    $(\mathcal{H}_\omega, \pi_\omega, \psi_\Unit)$ via the unitary
    intertwiner
    \begin{equation}
        \label{eq:GNSTootherGNSlike}
        U: \mathcal{H}_\omega \ni \psi_b 
        \; \longmapsto \;
        U \psi_b = \pi(b) \Omega \in \mathcal{H}.
    \end{equation}
\end{theorem}
The proof consists essentially in showing that $U$ is well-defined at
all. Then the remaining properties are immediate.
\begin{example}
    \label{example:BHonH}
    Let $\mathcal{H}$ be a pre-Hilbert space and $\phi \in
    \mathcal{H}$ a unit vector, $\SP{\phi, \phi} = 1$. Then for
    $\mathcal{A} = \mathfrak{B}(\mathcal{H})$ and for
    \begin{equation}
        \label{eq:ExpectationValueGNS}
        \omega(A) = \SP{\phi, A\phi}
    \end{equation}
    one recovers the defining $^*$-representation of $\mathcal{A}$ on
    $\mathcal{H}$ as the GNS representation corresponding to $\omega$.
    Note that one can replace $\mathfrak{B}(\mathcal{H})$ by
    $\mathfrak{F}(\mathcal{H})$ as well.
\end{example}

A slight generalization is obtained for the following situation: Let
$\mathcal{B} \subseteq \mathcal{A}$ be a $^*$-ideal and let $\omega:
\mathcal{B} \longrightarrow \ring{C}$ be a positive linear functional
which does not necessarily extend to $\mathcal{A}$. Let
$\mathcal{J}_\omega \subseteq \mathcal{B}$ be its Gel'fand ideal and
let $(\mathcal{H}_\omega = \mathcal{B} \big/ \mathcal{J}_\omega,
\pi_\omega)$ be the corresponding GNS representation of $\mathcal{B}$.
In this situation we have:
\begin{lemma}
    \label{lemma:BinAGNS}
    $\mathcal{J}_\omega \subseteq \mathcal{A}$ is a left ideal in
    $\mathcal{A}$ as well whence the GNS representation $\pi_\omega$
    of $\mathcal{B}$ extends canonically to $\mathcal{A}$ by the
    definition $\pi_\omega(a) \psi_b = \psi_{ab}$ and yields a
    $^*$-representation of $\mathcal{A}$ on $\mathcal{H}_\omega$.
\end{lemma}
The proof is again a consequence of the Cauchy-Schwarz inequality, see
\cite[Cor.~1]{bordemann.waldmann:1998a}. Nevertheless, this will be
very useful in the examples later.

%
%

\subsection{The case of $\delta$, Schrödinger and trace functionals}
\label{subsec:deltaSchroedingerTrace}

Let us now come back to deformation quantization and the examples of
positive functionals as in Section~\ref{sec:ExamplesPosFun}. We want
to determine their GNS representations explicitly.

\subsubsection*{The $\delta$-functional and the Wick star product}
\label{subsubsec:deltaGNS}

From the explicit formula for the Wick star product
\eqref{eq:ccfstarwickf} we see that the Gel'fand ideal of the
$\delta$-functional is simply given by
\begin{equation}
    \label{eq:GelfandDelta}
    \mathcal{J}_\delta =
    \left\{ f \in C^\infty(M)[[\lambda]]
        \; \left| \; 
            \frac{\partial^r f}
            {\partial \cc{z}^{i_1} \cdots \partial \cc{z}^{i_r}}(0) =
            0
            \;
            \textrm{for all}
            \;
            r \in \mathbb{N}_0,
            i_1, \ldots, i_r = 1, \ldots, n
        \right. 
    \right\}.
\end{equation}
In order to obtain explicit formulas for the GNS representation we
consider the $\mathbb{C}[[\lambda]]$-module
\begin{equation}
    \label{eq:deltaHilbert}
    \mathcal{H} = \mathbb{C}[[\cc{y}^1, \ldots, \cc{y}^n]][[\lambda]]
\end{equation}
which we endow with the inner product
\begin{equation}
    \label{eq:DeltaInnerProd}
    \SP{\phi, \psi} 
    = \sum_{r=0}^\infty \frac{(2\lambda)^r}{r!}
    \sum_{i_1, \ldots, i_r=1}^n
    \cc{
      \frac{\partial^r \phi}
      {\partial \cc{y}^{i_1} \cdots \partial \cc{y}^{i_r}}(0)
    }
    \;
    \frac{\partial^r \psi}{\partial \cc{y}^{i_1} \cdots \partial \cc{y}^{i_r}}(0).
\end{equation}
Clearly, this is well-defined as formal power series in $\lambda$ and
turns $\mathcal{H}$ into a pre-Hilbert space over
$\mathbb{C}[[\lambda]]$. Then we have the following characterization
of the GNS representation corresponding to the $\delta$-functional
\cite{bordemann.waldmann:1998a}:
\begin{theorem}[Formal Bargmann-Fock Representation]
    \label{theorem:GNSWick}
    The GNS pre-Hilbert space $\mathcal{H}_\delta$ is isometrically
    isomorphic to $\mathcal{H}$ via the unitary map
    \begin{equation}
        \label{eq:UnitaryGNSWick}
        U: \mathcal{H}_\delta \ni \psi_f 
        \; \longmapsto \;
        \sum_{r=0}^\infty \frac{1}{r!} 
        \sum_{i_1, \ldots, i_r=1}^n
        \frac{\partial^r f}
        {\partial \cc{z}^{i_1} \cdots \partial\cc{z}^{i_r}}(0)
        \;
        \cc{y}^{i_1} \cdots \cc{y}^{i_r}
        \in \mathcal{H},
    \end{equation}
    i.e. the formal $\cc{z}$-Taylor expansion around $0$. This way the
    GNS representation $\pi_\delta$ on $\mathcal{H}_\delta$ becomes
    the formal Bargmann-Fock representation in Wick ordering
    \begin{equation}
        \label{eq:FormalBargmannFock}
        \wickrep(f) = \sum_{r,s=0}^\infty \frac{(2\lambda)^r}{r!s!}
        \sum_{i_1, \ldots, i_r, j_1, \ldots, j_s=1}^n
        \frac{\partial^{r+s}f}
        {\partial \cc{z}^{j_1} \cdots \partial \cc{z}^{j_1} 
          \partial z^{i_1} \cdots \partial z^{i_r}}
        (0)
        \;
        \cc{y}^{j_1}\cdots \cc{y}^{j_s}
        \frac{\partial^r}{\partial \cc{y}^{i_1} \cdots \partial \cc{y}^{i_r}},
    \end{equation}
    i.e. $\wickrep(f) = U \pi_\delta(f) U^{-1}$.
\end{theorem}
With the explicit formulas for $\mathcal{J}_\delta$, $\starwick$ and
$U$ the proof is a straightforward verification. In particular,
$\wickrep$ is indeed the Bargmann-Fock representation in normal
ordering (Wick ordering): Specializing \eqref{eq:FormalBargmannFock}
we have
\begin{equation}
    \label{eq:NormalOrdering}
    \wickrep(z^i) = 2\lambda \frac{\partial}{\partial \cc{y}^i}
    \quad
    \textrm{and}
    \quad
    \wickrep(\cc{z}^i) = \cc{y}^i,
\end{equation}
together with normal ordering for higher polynomials in $z^i$ and
$\cc{z}^i$. Thus we obtain exactly the creation and annihilation
operators.

In fact, one can make the relation to the well-known `convergent'
Bargmann-Fock representation even more transparent. Recall that the
Bargmann-Fock Hilbert space is given by the \emph{anti-holomorphic}
functions
\begin{equation}
    \label{eq:BFHilbertSpace}
    \mathfrak{H}_{\mathrm{\scriptscriptstyle BF}}
    =
    \left\{f \in \cc{\mathcal{O}}(\mathbb{C}^n)
        \; \left| \; 
            \frac{1}{(2\pi\hbar)^n} \int |f(\cc{z})|^2
            \E^{-\frac{\cc{z}z}{2\hbar}} \D z \D \cc{z}
            < \infty
        \right.
    \right\},
\end{equation}
which are square-integrable with respect to the above Gaussian
measure, see \cite{bargmann:1967a, bargmann:1961a} as well as
\cite{landsman:1998a}. Here $\hbar > 0$ is Planck's constant. The
inner product is then the corresponding $L^2$-inner product and it is
well-known that those anti-holomorphic functions are a closed
sub-space of all square integrable functions $L^2(\mathbb{C}^n,
\E^{-\frac{\cc{z}z}{2\hbar}} \D z \D \cc{z})$. The quantization of
polynomials in $\cc{z}$ and $z$ is given by
\begin{equation}
    \label{eq:piBFannihilator}
    \pi_{\mathrm{\scriptscriptstyle BF}}(z^k) 
    = 2\hbar \frac{\partial}{\partial \cc{z}^k}
    = a_k
\end{equation}
\begin{equation}
    \label{eq:piBFcreator}
    \pi_{\mathrm{\scriptscriptstyle BF}}(\cc{z}^k) = \cc{z}^k
    = a_k^\dag
\end{equation}
plus \emph{normal ordering} for the higher monomials. Here
\eqref{eq:piBFannihilator} and \eqref{eq:piBFcreator} are densely
defined operators on $\mathfrak{H}_{\mathrm{\scriptscriptstyle BF}}$
which turn out to be mutual adjoints when the domains are chosen
appropriately.  Then the formal Bargmann-Fock space $\mathcal{H}$
together with the inner product \eqref{eq:DeltaInnerProd} can be seen
as asymptotic expansion of $\mathfrak{H}_{\mathrm{\scriptscriptstyle
    BF}}$ for $\hbar \longrightarrow 0$ and similarly $\wickrep$
arises as asymptotic expansion of $\pi_{\mathrm{\scriptscriptstyle
    BF}}$. Note that for a wide class of elements in $\mathcal{H}$ and
a large class of observables like the polynomials the asymptotic
expansion is already the \emph{exact} result.

In both cases the `vacuum vector' is just the constant function $1$
out of which all anti-holomorphic `functions' are obtained by
successively applying
the creation operators $\wickrep(\cc{z}^k)$ or
$\pi_{\mathrm{\scriptscriptstyle BF}}(\cc{z}^k)$, respectively.

\begin{remark}
    \label{remark:KaehlerFramework}
    In the formal framework, similar results can be obtained easily
    for any Kähler manifold being equipped with the Fedosov star
    product of Wick type, i.e. the star product with separation of
    variables according to the Kähler polarization
    \cite{bordemann.waldmann:1997a, bordemann.waldmann:1998a,
      karabegov:1996a}. In particular, all $\delta$-functionals are
    still positive linear functionals and essentially all formulas are
    still valid if one replaces the formal $\cc{z}$-Taylor expansion
    with the `anti-holomorphic part' of the Fedosov-Taylor series, see
    \cite{bordemann.waldmann:1998a}.  Note however, that the
    convergent analog is much more delicate and requires additional
    assumptions on the topology in the compact case. There is an
    extensive literatur on this topic, see e.g.
    \cite{cahen.gutt.rawnsley:1995a, cahen.gutt.rawnsley:1994a,
      cahen.gutt.rawnsley:1993a, cahen.gutt.rawnsley:1990a,
      landsman:1998a, berezin:1975a, berezin:1975b, berezin:1975c,
      karabegov.schlichenmaier:2001a, schlichenmaier:2000a,
      bordemann.meinrenken.schlichenmaier:1991a}. It would be
    interesting to understand the relations between both situations
    better, in particular concerning the representation point of view.
\end{remark}

\subsubsection*{The Schrödinger representation}
\label{subsubsec:SchroedingerRep}

To obtain the most important representation for mathematical physics,
the Schrödinger representation, we consider again the Weyl-Moyal star
product $\starweyl$ on $T^*\mathbb{R}^{n}$ together with the
Schrödinger functional $\omega$ as in
\eqref{eq:Schroedingerfunctional}. The functional is positive and
defined on the $^*$-ideal $C^\infty_0(T^*\mathbb{R}^n)[[\lambda]]$ of
$C^\infty(T^*\mathbb{R}^n)[[\lambda]]$. Thus we are in the situation
of Lemma~\ref{lemma:BinAGNS} since $\omega$ does not have an extension
to the whole $^*$-algebra. Nevertheless, the GNS representation
extends canonically to all observables. From
\eqref{eq:SchroedingerOnProducts} we immediately obtain
\begin{equation}
    \label{eq:SchroedingerGelfand}
    \mathcal{J}_\omega 
    = \left\{ f \ni C^\infty_0(T^*\mathbb{R}^n)[[\lambda]]
    \; \big| \; \iota^*Nf = 0 \right\}.
\end{equation}
This allows to identify the GNS representation of the Schrödinger
functional explicitly. It is the formal Schrödinger representation on
`formal' wave functions \cite{bordemann.waldmann:1998a}:
\begin{theorem}[Formal Schrödinger Representation]
    \label{theorem:SchroedingerRep}
    The GNS pre-Hilbert space $\mathcal{H}_\omega$ of the Schrödinger
    functional is isometrically isomorphic to the `formal wave
    functions' $\mathcal{H} = C^\infty_0(\mathbb{R}^n)[[\lambda]]$
    with inner product
    \begin{equation}
        \label{eq:FWFInnerProduct}
        \SP{\phi, \psi} = \int_{\mathbb{R}^n} \cc{\phi(q)} \psi(q) \D^nq
    \end{equation}
    via the unitary map
    \begin{equation}
        \label{eq:UnitarySchroedingerGNS}
        U: \mathcal{H}_\omega \ni \psi_f
        \; \longmapsto \;
        \iota^*Nf \in C^\infty_0(\mathbb{R}^n)[[\lambda]].
    \end{equation}
    The GNS representation $\pi_\omega$ becomes the formal Schrödinger
    representation (in Weyl ordering)
    \begin{equation}
        \label{eq:SchroedingerRepWeylOrder}
        \weylrep(f)
        =
        \sum_{r=0}^\infty \frac{1}{r!}
        \left(\frac{\lambda}{\I}\right)^r
        \sum_{i_1, \ldots, i_r = 1}^n
        \left.
            \frac{\partial^r (Nf)}{\partial p_{i_1} \cdots \partial p_{i_r}}
        \right|_{p=0}
        \frac{\partial^r}{\partial q^{i_1} \cdots \partial q^{i_r}},
    \end{equation}
    i.e. $\weylrep(f) = U \pi_\omega(f) U^{-1}$.
\end{theorem}
The proof consists again in a simple verification that $U$ is
well-defined and has the desired properties. Then
\eqref{eq:SchroedingerRepWeylOrder} is a straightforward computation.
Note that $\weylrep$ is indeed the usual Schrödinger representation
with $\hbar$ replaced by the formal parameter $\lambda$. In
particular,
\begin{equation}
    \label{eq:Schroedingerqp}
    \weylrep(q^k) = q^k
    \quad
    \textrm{and}
    \quad
    \weylrep(p_\ell) = - \I\lambda \frac{\partial}{\partial q^\ell}
\end{equation}
together with Weyl's \emph{total symmetrization rule} for the higher
monomoials. Here the correct combinatorics is due to the operator $N$.
Without $N$ in \eqref{eq:SchroedingerRepWeylOrder} one obtains the
Schrödinger representation in \emph{standard ordering} (all `$p_\ell$'
to the right). Using the \emph{standard ordered star product}
\begin{equation}
    \label{eq:starstd}
    f \starstd g = 
    \mu \circ 
    \E^{\frac{\lambda}{\I} \sum_{k=1}^n
      \frac{\partial}{\partial p_k} 
      \otimes \frac{\partial}{\partial q^k}}
    f \otimes g,
\end{equation}
which is equivalent to the Weyl-Moyal star product via the equivalence
transformation $N$, we can write
\begin{equation}
    \label{eq:weylrepNiceFormula}
    \weylrep(f)\psi = \iota^*(Nf \starstd \pi^*\psi)
\end{equation}
for $f \in C^\infty(T^*\mathbb{R}^n)[[\lambda]]$ and $\psi \in
C^\infty_0(\mathbb{R}^n)[[\lambda]]$. The corresponding standard
ordered representation is then simply given by
\begin{equation}
    \label{eq:stdrep}
    \stdrep(f)\psi = \iota^*(f \starstd \pi^*\psi).
\end{equation}
This is precisely the usual symbol calculus for differential operators
(in standard ordering) when we replace $\lambda$ by $\hbar$ and
restrict to polynomial functions in the momenta.
\begin{remark}[Formal Schrödinger Representation]
    \label{remark:SchroedingerGNSGlobal}
    ~
    \begin{compactenum}
    \item The `formal' Schrödinger representation can again be
        obtained from integral formulas for the Weyl-ordered symbol
        calculus by asymtptic expansions for $\hbar \longrightarrow
        0$. The asymtotic formulas are already exact for functions
        which are polynomial in the momenta.
    \item There are geometric generalizations not only for
        $\starweyl$, $\starstd$ and $N$ as discussed in
        Section~\ref{subsec:SchroedingerFunctional} but also the whole
        GNS construction can be translated to the geometric framework
        of cotangent bundles. Even the formulas for the
        representations $\weylrep$ and $\stdrep$ are still valid if
        partial derivatives are replaced by appropriate covariant
        derivatives, see \cite{bordemann.neumaier.waldmann:1999a,
          bordemann.neumaier.waldmann:1998a,
          bordemann.neumaier.pflaum.waldmann:2003a}.
    \end{compactenum}
\end{remark}

\subsubsection*{GNS representation of traces and KMS functionals}
\label{subsubsec:GNSforTracesKMS}

We consider again a connected symplectic manifold $M$ with a Hermitian
star product $\star$ and its positive trace functional $\tr:
C^\infty_0(M)[[\lambda]] \longrightarrow \mathbb{C}[[\lambda]]$ as in
Section~\ref{subsec:PositiveTraces}. Here $\tr$ is defined on a
$^*$-ideal whence we can again apply Lemma~\ref{lemma:BinAGNS} to
extend the GNS representation of $\tr$ to the whole algebra of
observables.

From the earlier investigation of the trace functional $\tr$ in
\eqref{eq:TraceLowestOrder} we see that the Gel'fand ideal of $\tr$ is
trivial,
\begin{equation}
    \label{eq:GelfandTrace}
    \mathcal{J}_{\tr} = \{0\}.
\end{equation}
Such positive functionals are called \emph{faithful}. Thus the GNS
pre-Hilbert space is simply $\mathcal{H}_{\tr} =
C^\infty_0(M)[[\lambda]]$ together with the inner product
\begin{equation}
    \label{eq:SPTrace}
    \SP{f, g}_{\tr} = \tr(\cc{f} \star g).
\end{equation}
The GNS representation is then the left regular representation
\begin{equation}
    \label{eq:GNSTrace}
    \pi_{\tr}(f)g = f \star g
\end{equation}
of $C^\infty(M)[[\lambda]]$ on the $^*$-ideal
$C^\infty_0(M)[[\lambda]]$.

It is clear that the analogous result still holds for the KMS states as
the star exponential function $\Exp(-\beta H)$ is invertible. Thus the
module structure is again the left regular one. Only the inner product
changed and is now twisted by the additional factor $\Exp(-\beta H)$
inside \eqref{eq:SPTrace}.

The characteristic property of the GNS representation of the trace
functional and the KMS functionals is that the \emph{commutant}
\begin{equation}
    \label{eq:Commutant}
    \pi(\mathcal{A})' 
    = \{ A \in \mathfrak{B}(\mathcal{H}) \; | \; 
    [\pi(a), A] = 0
    \; 
    \textrm{for all}
    \; a \in \mathcal{A} \}
\end{equation}
of the representation is as big as the algebra of observables
itself. In fact, it is anti-isomorphic
\begin{equation}
    \label{eq:CommutantAntiIso}
    \pi(\mathcal{A})' \cong \mathcal{A}^{\mathrm{\scriptscriptstyle opp}},
\end{equation}
since the commutant is given by all \emph{right-multiplications}. This
is in some sense the beginning of an algebraic `baby-version' of the
Tomita-Takesaki theory as it is well-known for operator algebras, see
\cite{waldmann:2000a} for more details.

%
%

\subsection{Deformation and classical limit of GNS representations}
\label{subsec:DefClassLimGNS}

Since the above examples prove that the GNS construction gives
physically meaningful representations also in formal deformation
quantization we shall now discuss the classical limit of GNS
representations and the corresponding deformation problem.

Let $\mathcal{A}$ be a $^*$-algebra over $\ring{C}$ and let $\star$ be
a Hermitian deformation. Thus $\boldsymbol{\mathcal{A}} =
(\mathcal{A}[[\lambda]], \star)$ is a $^*$-algebra over
$\ring{C}[[\lambda]]$. Then we want to understand how one can
construct a $^*$-representation $\pi$ out of a $^*$-representation
$\boldsymbol{\pi}$ of $\boldsymbol{\mathcal{A}}$ on some pre-Hilbert
space $\boldsymbol{\mathcal{H}}$ over $\ring{C}[[\lambda]]$. It turns
out that we can always take the \emph{classical limit} of a
$^*$-representation.

First we consider a pre-Hilbert space $\boldsymbol{\mathcal{H}}$ over
$\ring{C}[[\lambda]]$ and define
\begin{equation}
    \label{eq:Hnull}
    \boldsymbol{\mathcal{H}}_0
    =
    \left\{\boldsymbol{\phi} \in \boldsymbol{\mathcal{H}} \; \big| \;
    \SP{\boldsymbol{\phi},
      \boldsymbol{\phi}}_{\boldsymbol{\mathcal{H}}}
    \big|_{\lambda = 0} = 0\right\}.
\end{equation}
By the Cauchy-Schwarz inequality for the inner product
\begin{equation}
    \label{eq:CSUInnerProduct}
    \SP{\boldsymbol{\phi}, \boldsymbol{\psi}}_{\boldsymbol{\mathcal{H}}}
    \cc{\SP{\boldsymbol{\phi}, \boldsymbol{\psi}}}_{\boldsymbol{\mathcal{H}}}
    \le
    \SP{\boldsymbol{\phi}, \boldsymbol{\phi}}_{\boldsymbol{\mathcal{H}}}
    \SP{\boldsymbol{\psi}, \boldsymbol{\psi}}_{\boldsymbol{\mathcal{H}}},
\end{equation}
which holds in general thanks to Lemma~\ref{lemma:Babylonian}, we see
that $\boldsymbol{\mathcal{H}}_0$ is a $\ring{C}[[\lambda]]$-submodule
of $\boldsymbol{\mathcal{H}}$. Thus we can define
\begin{equation}
    \label{eq:clH}
    \cl(\boldsymbol{\mathcal{H}}) 
    = \boldsymbol{\mathcal{H}} \big/ \boldsymbol{\mathcal{H}}_0
\end{equation}
as a $\ring{C}[[\lambda]]$-module. The canonical projection will be
called \emph{classical limit map} and is denoted by $\cl$. Since
$\boldsymbol{\mathcal{H}}_0$ contains
$\lambda\boldsymbol{\mathcal{H}}$ the $\ring{C}[[\lambda]]$-module
structure of $\cl(\boldsymbol{\mathcal{H}})$ is rather trivial:
$\lambda$ always acts as zero. Hence we can forget this
$\ring{C}[[\lambda]]$-module structure and restrict it to the
$\ring{C}$-module structure. Then we define
\begin{equation}
    \label{eq:SPcl}
    \SP{\cl(\boldsymbol{\phi}),
      \cl(\boldsymbol{\psi})}_{\cl(\boldsymbol{\mathcal{H}})}
    =
    \SP{\boldsymbol{\phi},\boldsymbol{\psi}}_{\boldsymbol{\mathcal{H}}} 
    \big|_{\lambda=0} \in \ring{C}.
\end{equation}
It is easy to see that this gives indeed a well-defined and positive
definite inner product for $\cl(\boldsymbol{\mathcal{H}})$ which
thereby becomes a pre-Hilbert space over $\ring{C}$.

Next we consider an adjointable map $\boldsymbol{A} \in
\mathfrak{B}(\boldsymbol{\mathcal{H}}_1, \boldsymbol{\mathcal{H}}_2)$
and define a $\ring{C}$-linear map $\cl(\boldsymbol{A}):
\cl(\boldsymbol{\mathcal{H}}_1) \longrightarrow
\cl(\boldsymbol{\mathcal{H}}_2)$ by
\begin{equation}
    \label{eq:clADef}
    \cl(\boldsymbol{A}) \cl(\boldsymbol{\phi}) 
    = \cl(\boldsymbol{A}\boldsymbol{\phi}).
\end{equation}
Since $\boldsymbol{A}$ is adjointable it turns out that this is
actually well-defined and $\cl(\boldsymbol{A})$ is again adjointable.
Moreover, it is easy to check that $\boldsymbol{A} \mapsto
\cl(\boldsymbol{A})$ is compatible with $\ring{C}$-linear
combinations, taking the adjoints and the compositions of operators to
adjoints and compositions of their classical limits.  Thus we obtain a
\emph{functor} $\cl$ from the category of pre-Hilbert spaces over
$\ring{C}[[\lambda]]$ to the category of pre-Hilbert spaces over
$\ring{C}$ which we shall call the \emph{classical limit}, see
\cite[Sect.~8]{bursztyn.waldmann:2001a}
\begin{equation}
    \label{eq:clPreHilbertFunctor}
    \cl: \mathrm{PreHilbert}(\ring{C}[[\lambda]])
    \longrightarrow
    \mathrm{PreHilbert}(\ring{C}).
\end{equation}
\begin{remark}
    \label{remark:HnullRichtigGross}
    The $\ring{C}[[\lambda]]$-submodule $\boldsymbol{\mathcal{H}}_0$
    contains $\lambda\boldsymbol{\mathcal{H}}$ but in general it is
    much larger as we shall see in the examples. Thus $\cl$ is not
    just the functor `modulo $\lambda$' but takes into account the
    whole pre-Hilbert space structure.
\end{remark}

Now the classical limit functor also induces a classical limit for
$^*$-representations as follows. For a $^*$-representation
$(\boldsymbol{\mathcal{H}}, \boldsymbol{\pi})$ of
$\boldsymbol{\mathcal{A}}$ we define $\cl(\boldsymbol{\pi}):
\mathcal{A} \longrightarrow
\mathfrak{B}(\cl(\boldsymbol{\mathcal{H}}))$ by
\begin{equation}
    \label{eq:clRep}
    \cl(\boldsymbol{\pi})(a) = \cl(\boldsymbol{\pi}(a))
\end{equation}
for $a \in \mathcal{A}$. It is straightforward to check that this
gives indeed a $^*$-representation of the undeformed $^*$-algebra
$\mathcal{A}$ on $\cl(\boldsymbol{\mathcal{H}})$. Moreover, for an
intertwiner $\boldsymbol{T}: (\boldsymbol{\mathcal{H}}_1,
\boldsymbol{\pi}_1) \longrightarrow (\boldsymbol{\mathcal{H}}_2,
\boldsymbol{\pi}_2)$ we use its classical limit $\cl(\boldsymbol{T})$
to obtain an intertwiner between the classical limits of
$(\boldsymbol{\mathcal{H}}_1, \boldsymbol{\pi}_1)$ and
$(\boldsymbol{\mathcal{H}}_2, \boldsymbol{\pi}_2)$. Then it is easy to
check that this gives a functor
\begin{equation}
    \label{eq:clRepFunctor}
    \cl: \rep(\boldsymbol{\mathcal{A}}) 
    \longrightarrow \rep(\mathcal{A}),
\end{equation}
still called the \emph{classical limit}. Thus we can always take the
classical limit of $^*$-representations, even in a canonical way
\cite{bursztyn.waldmann:2001a}.

This immediately raises the question whether the converse is true as
well: Can we always \emph{deform} a given $^*$-representation of
$\mathcal{A}$ into a $^*$-representation of $\boldsymbol{\mathcal{A}}$
such that the above defined classical limit gives back the
$^*$-representation we started with?

In general this is a very difficult question whence we consider the
particular case of \emph{GNS representations}.  Thus let
$\boldsymbol{\omega} : \boldsymbol{\mathcal{A}} \longrightarrow
\ring{C}[[\lambda]]$ be a positive functional of the deformed
$^*$-algebra with corresponding GNS representation
$(\boldsymbol{\mathcal{H}_\omega}, \boldsymbol{\pi_\omega})$.  Then we
have the following result \cite[Thm.~1]{waldmann:2001a}:
\begin{theorem}
    \label{theorem:DeformedGNS}
    The classical limit $\cl(\boldsymbol{\mathcal{H}_\omega},
    \boldsymbol{\pi_\omega})$ is unitarily equivalent to the GNS
    representation $(\mathcal{H}_{\omega_0}, \pi_{\omega_0})$
    corresponding to the classical limit $\omega_0 =
    \cl(\boldsymbol{\omega})$ via the unitary intertwiner
    \begin{equation}
        \label{eq:UnitaryDeformedGNS}
        U: \cl(\boldsymbol{\mathcal{H}_\omega}) \ni
        \cl\left(\psi_{\boldsymbol{a}}\right)
        \; \longmapsto \;
        \psi_{\cl(\boldsymbol{a})} \in \mathcal{H}_{\omega_0},
    \end{equation}
    where $\boldsymbol{a} = a_0 + \lambda a_1 + \cdots \in
    \boldsymbol{\mathcal{A}}$ and $\cl(\boldsymbol{a}) = a_0 \in
    \mathcal{A}$.
\end{theorem}
One just has to check that the map $U$ is actually well-defined. The
remaining features of $U$ are obvious. Note however, that though the
theorem is not very difficult to prove, it is non-trivial in so far as
there is \emph{no} simple relation between the Gel'fand ideal of
$\omega_0$ and $\boldsymbol{\omega}$: the later is usually much
smaller than $\mathcal{J}_{\omega_0}[[\lambda]]$.

An immediate consequence is that for positive deformations we can
always deform GNS representations since we can deform the
corresponding positive functionals:
\begin{corollary}
    \label{corollary:GNSDeformable}
    Let $\boldsymbol{\mathcal{A}}$ be a positive deformation of
    $\mathcal{A}$. Then any direct orthogonal sum of GNS
    representations of $\mathcal{A}$ can be deformed.
\end{corollary}
\begin{remark}[Classical Limit and Deformation of GNS Representations]
    \label{remark:DeformedGNS}
    ~
    \begin{compactenum}
    \item Even in the very nice cases there might be representations
        which are \emph{not} direct sums of GNS representations. In
        the $C^*$-algebraic case, every $^*$-representation is known
        to be a \emph{topological} direct sum of GNS representations.
    \item Thanks to Theorem~\ref{theorem:ComPosDef} the above
        corollary applies for star products.
    \item It is a nice exercise to exemplify the theorem for the three
        GNS representations we have discussed in detail in
        Section~\ref{subsec:deltaSchroedingerTrace}, see also
        \cite{waldmann:2001a}.
    \end{compactenum}
\end{remark}

%
%

\section{General $^*$-representation theory}
\label{sec:GeneralRepTheory}

Given a $^*$-algebra $\mathcal{A}$ over $\ring{C}$ the aim of
representation theory would be (in principle) to understand first the
structure of the convex cone of positive functionals, then the
resulting GNS representations and finally the whole category of
$^*$-representations $\Rep(\mathcal{A})$. Of course, beside for very
simple examples this is rather hopeless from the beginning and we can
not expect to get some `final' answers in this fully general algebraic
approach.

One less ambitious aim would be to \emph{compare} the representation
theories of two given $^*$-algebras in a functorial sense
\begin{equation}
    \label{eq:CompareRep}
    \Rep(\mathcal{A}) \leftrightharpoons \Rep(\mathcal{B}),
\end{equation}
and determine whether they are \emph{equivalent}. This is the
principal question of \emph{Morita theory}. It turns out that even if
we do not understand the representation theories of $\mathcal{A}$ and
$\mathcal{B}$ themselves completely, it might still be possible to
understand whether they are equivalent or not.

The question of finding some relations between the two representation
theories is interesting, even if one does not expect to obtain an
equivalence. The physical situation we have in mind is the following:
Consider a `big' phase space $(M, \pi)$ with some \emph{constraint
  manifold} $\iota: C \hookrightarrow M$, like e.g. the zero level set
of a momentum map or just some coisotropic submanifold (which
corresponds to first class constraints). Then the `physically
interesting' phase space would be the \emph{reduced phase space}
$M_{\mathrm{red}} = C\big/\!\!\!\sim$ endowed with the reduced Poisson
structure $\pi_{\mathrm{red}}$. Any \emph{gauge theory} is an example
of this situation. See e.g. the monography \cite{ortega.ratiu:2004a}
for details and further references on phase space reduction.

Of course we would like to understand the quantum theory of the whole
picture, i.e. the \emph{quantization of phase space reduction}. In
deformation quantization this amounts to find a star product $\star$
for $M$ which induces a star product $\star_{\mathrm{red}}$ on
$M_{\mathrm{red}}$. This has been discussed in various ways in
deformation quantization, see e.g.\cite{fedosov:1998a, fedosov:1994b,
  bordemann.herbig.waldmann:2000a, bordemann:2000a}, culminating
probably in the recent work of Bordemann \cite{bordemann:2004a:pre}.

Having understood the relation between the quantized observable
algebra $\mathcal{A} = (C^\infty(M)[[\lambda]], \star)$ and
$\mathcal{A}_{\mathrm{red}} = (C^\infty(M_{\mathrm{red}})[[\lambda]],
\star_{\mathrm{red}})$ we would like to understand also the relations
between their representation theories
\begin{equation}
    \label{eq:AAredReps}
    \Rep(\mathcal{A}) 
    \leftrightharpoons \Rep(\mathcal{A}_{\mathrm{red}}),
\end{equation}
and now we can not expect to get an equivalence of categories as the
geometrical structure on $M$ may be much richer `far away' from the
constraint surface $C$ whence it is not seen in the reduction process.
Nevertheless, already \emph{some} relation would be helpful.

Motivated by this we give now a rather general procedure to construct
functors $\Rep(\mathcal{A}) \longrightarrow \Rep(\mathcal{B})$.

%
%

\subsection{$^*$-Representation theory on pre-Hilbert modules}
\label{subsec:RepOnPreHilbertModules}

First we have to enlarge the notion of representation in order to get
a more coherent picture: we have to go beyond representations on
pre-Hilbert spaces over $\ring{C}$ but use general \emph{pre-Hilbert
  modules} instead \cite{bursztyn.waldmann:2003a:pre}. We consider an
auxilliary $^*$-algebra $\mathcal{D}$ over $\ring{C}$.
\begin{definition}[Pre-Hilbert Module]
    \label{definition:PreHilbertModule}
    A pre-Hilbert right $\mathcal{D}$-module $\HD$ is a right
    $\mathcal{D}$-module together with a map
    \begin{equation}
        \label{eq:SPDmodule}
        \SPD{\cdot,\cdot}: \mathcal{H} \times \mathcal{H}
        \longrightarrow \mathcal{D}
    \end{equation}
    such that
    \begin{compactenum}
    \item $\SPD{\cdot,\cdot}$ is $\ring{C}$-linear in the second
        argument.
    \item $\SPD{\phi, \psi} = \IP{\psi, \phi}{}{*}{\mathcal{D}}$ for
        $\phi, \psi \in \mathcal{H}$.
    \item $\SPD{\phi, \psi \cdot d} = \SPD{\phi, \psi} d$ for $\phi,
        \psi \in \mathcal{H}$ and $d \in \mathcal{D}$.
    \item $\SPD{\cdot,\cdot}$ is non-degenerate, i.e. $\SPD{\phi,
          \cdot} = 0$ implies $\phi =0$ for $\phi \in \mathcal{H}$.
    \item $\SPD{\cdot,\cdot}$ is completely positive, i.e. for all $n
        \in \mathbb{N}$ and all $\phi_1, \ldots, \phi_n \in
        \mathcal{H}$ we have $(\SPD{\phi_i, \phi_j}) \in
        M_n(\mathcal{D})^+$.
    \end{compactenum}
    In addition, $\SPD{\cdot,\cdot}$ is called \emph{strongly
      non-degenerate} if the map $x \mapsto \SPD{x, \cdot} \in \HD^* =
    \Hom_{\mathcal{D}}(\HD, \mathcal{D})$ is bijective.
\end{definition}

As we will have different inner products with values in different
algebras simultanously, we shall sometimes index the algebra in our
notion to avoid missunderstandings as well as for the modules.
Clearly we have an analogous definition for pre-Hilbert \emph{left}
$\mathcal{D}$-modules where now the inner product is $\ring{C}$-linear
and $\mathcal{D}$-linear to the \emph{left} in the \emph{first}
argument.
\begin{remark}[Pre-Hilbert Modules]
    \label{remark:PreHilbertModules}
    ~
    \begin{compactenum}
    \item This definition generalizes the notion of Hilbert modules
        over $C^*$-algebras, see e.g. \cite{lance:1995a,
          landsman:1998a}. In this case it is well-known that
        positivity of the inner product implies complete positivity,
        see e.g.~\cite[Lem.~4.2]{lance:1995a}.
    \item Pre-Hilbert spaces are obtained for $\mathcal{D} =
        \ring{C}$, the complete positivity of $\SP{\cdot,\cdot}$ for
        pre-Hilbert spaces over $\ring{C}$ can be shown
        \cite[App.~A]{bursztyn.waldmann:2001a} to be a consequence of
        the positivity of the inner product.
    \item We have obvious definitions for $\mathfrak{B}(\HD, \HpD)$,
      $\mathfrak{B}(\HD)$, $\mathfrak{F}(\HD, \HpD)$, and
      $\mathfrak{F}(\HD)$ analogously to pre-Hilbert spaces.
  \item $\HD$ is a $(\mathfrak{B}(\HD), \mathcal{D})$-bimodule since
      adjointable operators are necessarily right
      $\mathcal{D}$-linear. Moreover, $\mathfrak{F}(\HD)$ is a
      $^*$-ideal inside the $^*$-algebra $\mathfrak{B}(\HD)$.
    \end{compactenum}
\end{remark}

The following example shows that such pre-Hilbert modules arise very
naturally in differential geometry:
\begin{example}[Hermitian vector bundles]
    \label{example:HermVecBundle}
    Let $E \longrightarrow M$ be a complex vector bundle with a
    Hermitian fibre metric $h$. Then the right module
    \begin{equation}
        \label{eq:GammaECM}
        \Gamma^\infty(E)_{C^\infty(M)}
    \end{equation}
    with the inner product defined by $\SP{s, s'}(x) = h_x(s(x),
    s'(x))$, where $x \in M$, is a pre-Hilbert right
    $C^\infty(M)$-module. In this case,
    \begin{equation}
        \label{eq:BGammaE}
        \mathfrak{B}\left(\Gamma^\infty(E)_{C^\infty(M)}\right)
        =
        \Gamma^\infty(\End(E)),
    \end{equation}
    with their usual action on $\Gamma^\infty(E)$ and $^*$-involution
    induced by $h$. Moreover, for the finite-rank operators we have
    $\mathfrak{F}(\Gamma^\infty(E)_{C^\infty(M)}) =
    \Gamma^\infty(\End(E))$ as well. This is clear in the case where
    $M$ is compact but it is also true for non-compact $M$ as sections
    of vector bundles are still \emph{finitely generated} modules over
    $C^\infty(M)$. In fact, all these statements are a consequence of
    the Serre-Swan-Theorem \cite{swan:1962a}.
\end{example}

\begin{definition}
    \label{definition:GeneralRep}
    A $^*$-representation $\pi$ of $\mathcal{A}$ on a pre-Hilbert
    $\mathcal{D}$-module $\HD$ is a $^*$-homomorphism
    \begin{equation}
        \label{eq:GenStarRep}
        \pi: \mathcal{A} \longrightarrow \mathfrak{B}(\HD).
    \end{equation}
\end{definition}
Clearly, we can transfer the notions of intertwiners to this framework
as well whence we obtain the category of $^*$-representations of
$\mathcal{A}$ on pre-Hilbert $\mathcal{D}$-modules which we denote by
$\rep[\mathcal{D}](\mathcal{A})$. The strongly non-degenerate ones are
denoted by $\Rep[\mathcal{D}](\mathcal{A})$ where again in the unital
case we require $^*$-representations to fulfill $\pi(\Unit) = \id$.

%
%

\subsection{Tensor products and Rieffel induction}
\label{subsec:TensorRieffel}

The advantage of looking at $\Rep[\mathcal{D}](\mathcal{A})$ for all
$\mathcal{D}$ and not just for $\mathcal{D} = \ring{C}$ is that we now
have a tensor product operation. This will give us functors for
studing $\Rep[\mathcal{D}](\mathcal{A})$ and in particular
$\Rep(\mathcal{A})$.

The construction will be based on Rieffel's internal tensor product of
pre-Hilbert modules. Rieffel proposed this originally for
$C^*$-algebras \cite{rieffel:1974b, rieffel:1974a}, see also
\cite{raeburn.williams:1998a, landsman:1998a, lance:1995a}, but the
essential construction is entirely algebraic whence we obtain a quite
drastic generalization, see \cite{bursztyn.waldmann:2001a,
  bursztyn.waldmann:2003a:pre}.

We consider $\CFB \in \rep[\mathcal{B}](\mathcal{C})$ and $\BEA \in
\rep[\mathcal{A}](\mathcal{B})$. Then we have the algebraic tensor
product $\CFB \tensor[B] \BEA$ which is a $(\mathcal{C},
\mathcal{A})$-bimodule in a natural way, since we started with
bimodules. Out of the given inner products on $\CFB$ and $\BEA$ we
want to construct an inner product $\SPFEA{\cdot,\cdot}$ with values
in $\mathcal{A}$ on this tensor product such that the left
$\mathcal{C}$-module structure becomes a $^*$-representation. This can
actually by done. We define $\SPFEA{\cdot,\cdot}$ for elementary
tensors by
\begin{equation}
    \label{eq:SPFEADef}
    \SPFEA{x \otimes \phi, y \otimes \psi}
    =
    \SPEA{\phi, \SPFB{x, y} \cdot \psi},
\end{equation}
for $x, y \in \mathcal{F}$ and $\phi, \psi \in \mathcal{E}$ and extend
this by $\ring{C}$-sesquilinearity to $\CFB \tensor[B] \BEA$.
\begin{remark}
    \label{remark:TensorProdWelldef}
    One can show by some simple computations that $\SPFEA{\cdot,
      \cdot}$ is indeed well-defined on the $\mathcal{B}$-tensor
    product. Moreover, it has the correct $\mathcal{A}$-linearity
    properties and $\mathcal{C}$ acts by adjointable operators. This
    is all rather straightforward.
\end{remark}
The problem is the non-degeneracy and the complete positivity. Here we
have the following result
\cite[Thm.~4.7]{bursztyn.waldmann:2003a:pre}:
\begin{theorem}
    \label{theorem:ComPosTensor}
    If $\SPEA{\cdot,\cdot}$ and $\SPFB{\cdot,\cdot}$ are completely
    positive then $\SPFEA{\cdot,\cdot}$ is completely positive as well.
\end{theorem}
\begin{proof}
    Let $\Phi^{(1)}, \ldots, \Phi^{(n)} \in \mathcal{F} \tensor[B]
    \mathcal{E}$ be given. Then we must show that the matrix $A =
    \left(\SPFEA{\Phi^{(\alpha)}, \Phi^{(\beta)}}\right) \in
    M_n(\mathcal{A})$ is positive. Thus let $\Phi^{(\alpha)} =
    \sum_{i=1}^N x^{(\alpha)}_i \otimes \phi^{(\alpha)}_i$ with
    $x^{(\alpha)}_i \in \mathcal{F}$ and $\phi^{(\alpha)}_i \in
    \mathcal{E}$ where we can assume without restriction that $N$ is
    the same for all $\alpha = 1, \ldots, n$.  First we claim that the
    map $f: M_{nN}(\mathcal{B}) \longrightarrow M_{nN}(\mathcal{A})$
    defined by
    \begin{equation}
        \label{eq:fMnBPositive}
        f: \left(B^{\alpha\beta}_{ij}\right) \; \longmapsto \;
        \left(
            \SPEA{\phi^{(\alpha)}_i, 
              B^{\alpha\beta}_{ij} \phi^{(\beta)}_j}
        \right)
    \end{equation}
    is positive. Indeed, we have for any $B =
    \left(B^{\alpha\beta}_{ij}\right) \in M_{nN}(\mathcal{B})$
    \[
    f(B^*B)
    = \left(
        \SPEA{\phi^{(\alpha)}_i, 
          (B^*B)^{\alpha\beta}_{ij} \cdot \phi^{(\beta)}_j}
    \right) 
    =
    \sum_{\gamma=1}^n \sum_{k=1}^N
    \left(
        \SPEA{B^{\gamma\alpha}_{ki} \cdot \phi^{(\alpha)}_i,
          B^{\gamma\beta}_{kj} \cdot \phi^{(\beta)}_j}
    \right),
    \]
    and each term $\SPEA{B^{\gamma\alpha}_{ki} \cdot
      \phi^{(\alpha)}_i, B^{\gamma\beta}_{kj} \cdot \phi^{(\beta)}_j}$
    is a positive matrix in $M_{nN}(\mathcal{A})$ since
    $\SPEA{\cdot,\cdot}$ is completely positive. Thus $f(B^*B) \in
    M_{nN}(\mathcal{A})^+$ whence by
    Remark~\ref{remark:PositiveStuff}, \textit{iii.)} we conclude that
    $f$ is a positive map.  Since $\SPFB{\cdot,\cdot}$ is completely
    positive, the matrix $B = \SPFB{x^{(\alpha)}_i, x^{(\beta)}_j}$ is
    positive. Thus
    \[
    f(B)
    =
    \left(
        \SPEA{\phi^{(\alpha)}_i, 
          \SPFB{x^{(\alpha)}_i, x^{(\beta)}_j} \cdot \phi^{(\beta)}_j
        }
    \right)
    \]
    is a positive matrix $f(B) \in M_{nN}(\mathcal{A})^+$.  Finally,
    the summation over $i, j$ is precisely the completely positive map
    $\tau: M_{nN}(\mathcal{A}) \longrightarrow M_n(\mathcal{A})$.
    Hence
    \[
    \begin{split}
        \tau(f(B))
        &=
        \left(
            \sum_{i, j=1}^N
            \SPEA{\phi^{(\alpha)}_i, 
              \SPFB{x^{(\alpha)}_i, x^{(\beta)}_j} \cdot \phi^{(\beta)}_j
            }
        \right)\\
        &=
        \left(
            \SPEA{
              \sum_{i=1}^N x^{(\alpha)}_i \otimes \phi^{(\alpha)}_i,
              \sum_{j=1}^N x^{(\beta)}_j \otimes \phi^{(\beta)}_j
            }
        \right)\\
        &=
        \left(\SPFEA{\Phi^{(\alpha)}, \Phi^{(\beta)}}\right)
        \in
        M_n(\mathcal{A})^+
    \end{split}
    \]
    and thus the theorem is shown.
\end{proof}

A final check shows that the degeneracy space $(\mathcal{F} \tensor[B]
\mathcal{E})^\bot$ of $\SPFEA{\cdot,\cdot}$ is preserved under the
$(\mathcal{C}, \mathcal{A})$-bimodule structure. Thus we can divide by
this degeneracy space and obtain a pre-Hilbert $\mathcal{A}$-module,
again together with a $^*$-representation of $\mathcal{C}$. We denote
this new $^*$-representation by
\begin{equation}
    \label{eq:FhattensorBE}
    \CFB \tensM[B] \BEA 
    = (\mathcal{F} \tensor[B] \mathcal{E})
    \big/(\mathcal{F} \tensor[B] \mathcal{E})^\bot 
    \in \rep[\mathcal{A}](\mathcal{C}).
\end{equation}
The whole procedure is canonical, i.e. compatible with intertwiners at
all stages. So finally we have a \emph{functor}
\begin{equation}
    \label{eq:hatTensorFun}
    \tensM[B]:
    \rep[\mathcal{B}](\mathcal{C})
    \times
    \rep[\mathcal{A}](\mathcal{B})
    \longrightarrow
    \rep[\mathcal{A}](\mathcal{C})
\end{equation}
which restricts to strongly non-degenerate $^*$-representations in the
following way
\begin{equation}
    \label{eq:hatTensSND}
    \tensM[B]:
    \Rep[\mathcal{B}](\mathcal{C})
    \times
    \rep[\mathcal{A}](\mathcal{B})
    \longrightarrow
    \Rep[\mathcal{A}](\mathcal{C}).
\end{equation}
By fixing one of the two arguments of this tensor product we obtain
the following two particular cases:
\begin{compactenum}
\item \emph{Rieffel induction:} We fix $\BEA$. Then the functor
    \begin{equation}
        \label{eq:RieffelE}
        \mathsf{R}_{\mathcal{E}} = \BEA \tensM[A] \cdot:
        \rep[\mathcal{D}](\mathcal{A}) 
        \longrightarrow
        \rep[\mathcal{D}](\mathcal{B}) 
    \end{equation}
    is called Rieffel induction.
\item \emph{Change of base ring:} We fix $\DGDp$ then the functor
    \begin{equation}
        \label{eq:ChangeOfBase}
        \mathsf{S}_{\mathcal{G}} = \cdot \tensM[D] \DGDp:
        \rep[\mathcal{D}](\mathcal{A}) 
        \longrightarrow
        \rep[\mathcal{D}'](\mathcal{A}) 
    \end{equation}
    is called the change of base ring functor.
\end{compactenum}
We clearly have the following commutative diagram
\begin{equation}
    \label{eq:ReSFCommute}
    \begin{CD}
        \rep[\mathcal{D}](\mathcal{A})
        @>\mathsf{S}_{\mathcal{G}}>>
        \rep\nolimits_{\mathcal{D}'}(\mathcal{A}) \\
        @V\mathsf{R}_{\mathcal{E}}VV
        @VV\mathsf{R}_{\mathcal{E}}V \\
        \rep[\mathcal{D}](\mathcal{B})
        @>\mathsf{S}_{\mathcal{G}}>>
        \rep\nolimits_{\mathcal{D}'}(\mathcal{B}),
    \end{CD}
\end{equation}
which commutes in the sense of functors, i.e. up to natural
transformations. This is due to the simple fact that the tensor
product $\tensM$ is associative up to the usual natural
transformations.

%
%

\subsection{A non-trivial example: Dirac's monopole}
\label{subsec:DiracMonopole}

The following example is a particular case of the results of
\cite[Sect.~4.2]{bursztyn.waldmann:2002a} and can be understood better
in the context of Morita equivalence. Nevertheless we mention the
example already here.

We consider the configuration space $Q = \mathbb{R}^3\setminus\{0\}$
of an electrically charged particle within the external field of a
magnetic monopole, which sits in the origin. Thus the magnetic field
is described by a closed two-form $B \in \Gamma^\infty(\Lambda^2
T^*Q)$ which is not exact due to the presence of a `magnetic charge'
at $0 \in \mathbb{R}^3$. We assume furthermore that $\frac{1}{2\pi} B$
is an integral two-form, i.e. the magnetic charge satisfies Dirac's
quantization condition. Mathematically this means that
$\frac{1}{2\pi}[B] \in \HdR^2(M, \mathbb{Z})$ is in the
\emph{integral} deRham cohomology.

Consider $T^*Q$ with the Weyl-Moyal star product $\starweyl$ and
replace now the canonical symplectic form $\omega_0$ by the formal
symplectic form $\omega_B = \omega_0 - \lambda \pi^*B$.  This is the
`minimal coupling' corresponding to switching on the magnetic field.
One can now construct by this minimal coupling a star product
$\star_B$ out of $\starweyl$ by essentially replacing \emph{locally}
the momentum variables $p_i$ by $p_i - \lambda A_i$ where $A_i$ are
the components of a \emph{local} potential $A \in \Gamma^\infty(T^*Q)$
of $B$, i.e. $\D A = B$. It turns out that $\star_B$ is actually
globally defined, i.e. independent of the choices of $A$ but only
depending on $B$. The characteristic class of $\star_B$ is given by
\begin{equation}
    \label{eq:classStarB}
    c(\star_B) = \I [\pi^*B].
\end{equation}

Since $B$ is integral it defines a (non-trivial) line bundle
\begin{equation}
    \label{eq:LtoQ}
    \ell \longrightarrow Q,
\end{equation}
whose Chern class is given by the class $\frac{1}{2\pi}[B]$. This line
bundle is unique up to isomorphism and up to tensoring with a flat
line bundle.  On $\ell$ we choose a Hermitian fibre metric $h$. Thus
we also have the pull-back bundle $L = \pi^*\ell$ with Chern class
$\frac{1}{2\pi}[\pi^*B]$ together with the corresponding pull-back
fibre metric $H = \pi^*h$.

Then it is a fact that on $\Gamma^\infty(L)[[\lambda]]$ there exists a
$(\star_B, \starweyl)$-bimodule structure deforming the classical
bimodule structure of $\Gamma^\infty(L)$ viewed as a
$C^\infty(T^*Q)$-bimodule. Moreover, there exists a deformation of the
Hermitian fibre metric $H$ into a $(C^\infty(T^*Q)[[\lambda]],
\starweyl)$-valued positive inner product $\boldsymbol{H}$. This way,
the sections of $L$ become a $^*$-representation of
$(C^\infty(T^*Q)[[\lambda]], \star_B)$ on a pre-Hilbert module over
$(C^\infty(T^*Q)[[\lambda]], \starweyl)$,
\begin{equation}
    \label{eq:GammaLDeformed}
    \left(\Gamma^\infty(L)[[\lambda]], \boldsymbol{H}\right)
    \in
    \Rep[\protect{(C^\infty(T^*Q)[[\lambda]], \starweyl)}]
    \left(C^\infty(T^*Q)[[\lambda]], \star_B\right).
\end{equation}
This construction can be made very precise using Fedosov's approach to
the construction of symplectic star products, for detail we refer to
\cite{bordemann.neumaier.pflaum.waldmann:2003a,
  bursztyn.waldmann:2002a, waldmann:2002b}.

Having such a bimodule we can use it for Rieffel induction. Since for
$\starweyl$ we have a representation which is of particular interest,
we apply the Rieffel induction functor to the Schrödinger
representation $(C^\infty_0(Q)[[\lambda]], \weylrep)$.  Then it is
another fact that the resulting $^*$-representation of $\star_B$ is
precisely the usual `Dirac-type' representation on the pre-Hilbert
space of sections $\Gamma^\infty_0(\ell)[[\lambda]]$ of $\ell$ endowed
with the inner product
\begin{equation}
    \label{eq:SectionLInnerProduct}
    \SP{s, s'} = \int_{Q} h(s, s') \D^n q.
\end{equation}
The representation is given as follows: The configuration space
variables act by multiplication operators while the corresponding
canonical conjugate momenta act by covariant derivatives using a
connection on $\ell$ whose curvature is given by $B$. This is exactly
the \emph{minimal coupling} expected for quantization in presence of a
magnetic field.
\begin{remark}
    \label{remark:MoritaDirac}
    The above `ad hoc' construction (observation) finds its deeper
    explanation in Morita theory stating that the above bimodule is
    actually an equivalence bimodule, see
    Section~\ref{subsec:MoritaEquivStarProducts}. Moreover, an
    arbitrary star product $\star'$ on $T^*Q$ turns out to be Morita
    equivalent to $\starweyl$ if and only if the characteristic class
    of $\star'$ is integral. This is the Dirac's quantization
    condition for magnetic monopoles in the light of Morita theory and
    Rieffel induction applied to the usual Schrödinger representation.
\end{remark}
\begin{remark}
    \label{remark:DiracCondOnTstarQ}
    Note also, that the whole construction works for \emph{any}
    cotangent bundle $T^*Q$. One has very explicit formulas for the
    star products as well as the representations on the sections of
    the involved line bundles, see \cite{bursztyn.waldmann:2002a,
      bordemann.neumaier.pflaum.waldmann:2003a, waldmann:2002b}.
\end{remark}

%
%

\section{Strong Morita equivalence and the Picard groupoid}
\label{sec:StrongMoritaEquivalence}

We shall now give an introduction to Morita theory of $^*$-algebras
over $\ring{C}$ based on the crucial notion of the Picard groupoid.

%
%

\subsection{Morita equivalence in the ring-theoretic setting}
\label{subsec:MoritaRingtheoretic}

As warming-up we start recalling the ring-theoretic situation. Thus
let $\mathcal{A}$, $\mathcal{B}$ be two $\mathbb{k}$-algebras, where
we consider only the \emph{unital} case for simplicity.
By $\mathcal{A}\Mod$ we denote the category of left
$\mathcal{A}$-modules where we always assume that $\Unit_{\mathcal{A}}
\cdot m = m$ for all $m \in \AM$ where $\AM \in \mathcal{A}\Mod$. The
morphisms of this category are just the usual left
$\mathcal{A}$-module morphisms.

Given a $(\mathcal{B}, \mathcal{A})$-bimodule $\BEA$ one obtains a
functor by tensoring
\begin{equation}
    \label{eq:BEATensorFunctor}
    \BEA \tensor[A] \cdot: 
    \mathcal{A}\Mod \longrightarrow \mathcal{B}\Mod.
\end{equation}
In particular, the canonical bimodule $\AAA$ gives a functor naturally
equivalent to the identity functor $\id_{\mathcal{A}\Mod}$.

This motivates the following definition of an equivalence bimodule in
this ring-theoretic framework: $\BEA$ is called a \emph{Morita
  equivalence bimodule} if it is `invertible' in the sense that there
exit bimodules $\AEpB$ and $\AEppB$ such that
\begin{equation}
    \label{eq:InvertibleBimodules}
    \BEA \tensor[A] \AEpB \cong \BBB
    \quad
    \textrm{and}
    \quad
    \AEppB \tensor[B] \BEA \cong \AAA
\end{equation}
as bimodules.

In this case, it is easy to see that the functor
\eqref{eq:BEATensorFunctor} is an equivalence of categories. This is
essentially the associativity of the tensor product up to a natural
transformation. Moreover, $\AEpB \cong \AEppB$ as bimodules and 
\begin{equation}
    \label{eq:InvBEAHom}
    \AEpB \cong \Hom_{\mathcal{A}} (\BEA, \mathcal{A})
\end{equation}
as bimodules. In addition, $\BEA$ is finitely generated and projective
over $\mathcal{A}$, i.e. of the form
\begin{equation}
    \label{eq:BEAfgpm}
    \BEA \cong e \mathcal{A}^n
    \quad
    \textrm{with}
    \quad
    e = e^2 \in M_n(\mathcal{A})
\end{equation}
and $\mathcal{B}$ is determined up to isomorphism by
\begin{equation}
    \label{eq:BisEndeA}
    \mathcal{B} \cong \End_{\mathcal{A}}(\EA)
    \cong e M_n(\mathcal{A}) e.
\end{equation}
Finally, the idempotent $e$ is \emph{full} in the sense that the ideal
in $\mathcal{A}$
generated by all the matrix coefficients $e_{ij}$ of $e = (e_{ij}) \in
M_n(\mathcal{A})$ coincides with the whole algebra $\mathcal{A}$.

The converse statement is true as well: Given a full idempotent $e \in
M_n(\mathcal{A})$ the $(eM_n(\mathcal{A})e, \mathcal{A})$-bimodule
$e\mathcal{A}^n$ is invertible in the above sense and gives an
equivalence of the categories of modules over the algebras
$\mathcal{B} = e M_n(\mathcal{A})e$ and $\mathcal{A}$. These are the
statements of the classical Morita theory, see
e.g.~\cite{morita:1958a, lam:1999a}.

%
%

\subsection{Strong Morita equivalence}
\label{subsec:StrongMorita}

We want to specialize the notion of Morita equivalence to the case of
$^*$-algebras over $\ring{C}$ such that the specialized Morita
equivalence implies the equivalence of the categories
$\Rep[\mathcal{D}](\cdot)$ for all auxilliary $^*$-algebras
$\mathcal{D}$. In fact, it will be an algebraic generalization of
Rieffel's notion of \emph{strong Morita equivalence} of
$C^*$-algebras, hence the name. We state the following definition
\cite{bursztyn.waldmann:2003a:pre}.
\begin{definition}[Strong Morita equivalence]
    \label{definition:SMEBimodule}
    Let $\mathcal{A}$, $\mathcal{B}$ be $^*$-algebras over $\ring{C}$.
    A $(\mathcal{B}, \mathcal{A})$-bimodule $\BEA$ with inner products
    $\SPEA{\cdot,\cdot}$ and $\BSPE{\cdot,\cdot}$ is called strong
    Morita equivalence bimodule if the following conditions are
    satisfied:
    \begin{compactenum}
    \item Both inner products are non-degenerate and completely
        positive.
    \item $\mathcal{B} \cdot \BEA = \BEA = \BEA \cdot \mathcal{A}$.
    \item Both inner products are full, i.e.
        \begin{equation}
            \label{eq:FullInnerProductA}
            \ring{C}\textrm{-}\mathrm{span}
            \{\SPEA{x,y} \; | \; x, y \in \BEA \} = \mathcal{A}
        \end{equation}
        \begin{equation}
            \label{eq:FullInnerProductB}
            \ring{C}\textrm{-}\mathrm{span}
            \{\BSPE{x,y} \; | \; x, y \in \BEA \} = \mathcal{B}
        \end{equation}
    \item We have the compatibility conditions
        \begin{equation}
            \label{eq:Bstaracts}
            \SPEA{x, b \cdot y} = \SPEA{b^* \cdot x, y}
        \end{equation}
        \begin{equation}
            \label{eq:Astaracts}
            \BSPE{x, y \cdot a} = \BSPE{x \cdot a^*, y}
        \end{equation}
        \begin{equation}
            \label{eq:SPAssoc}
            \BSPE{x, y} \cdot z = x \cdot \SPEA{y, z}
        \end{equation}
        for all $b \in \mathcal{B}$, $a \in \mathcal{A}$, and $x, y, z
        \in \BEA$.
    \end{compactenum}
    If such a bimodule exists then $\mathcal{A}$ and $\mathcal{B}$ are
    called strongly Morita equivalent.
\end{definition}
\begin{remark}[$^*$-Morita equivalence]
    \label{remark:StarMoritaEquiv}
    Without the above complete positivity requirements this notion is
    called $^*$-Morita equivalence and the bimodules are called
    $^*$-Morita equivalence bimodules, see Ara's works
    \cite{ara:1999b,ara:1999a}.
\end{remark}
\begin{remark}
    \label{remark:BSPDetermined}
    It is easy to see that the $\mathcal{B}$-valued inner product is
    completely determined by \eqref{eq:SPAssoc} since this simply
    means that $\BSPE{x, y}$ acts as $\Theta_{x, y}$ or, in Dirac's
    bra-ket notation, as $|x\rangle\langle y|$.
\end{remark}

From now on we shall assume that all $^*$-algebras are
\emph{non-degenerate} in the sense that $a \cdot \mathcal{A} = 0$
implies $a = 0$ and \emph{idempotent} in the sense that $a = \sum_i
b_i c_i$ for any $a \in \mathcal{A}$ with some $b_i, c_i \in
\mathcal{A}$. In particular, \emph{unital} $^*$-algebras are
non-degenerate and idempotent.  This restriction is reasonable
according to the following standard example:
\begin{example}
    \label{example:MnAAnA}
    Consider the $(M_n(\mathcal{A}), \mathcal{A})$-bimodule
    $\mathcal{A}^n$ for $n \ge 1$ with the canonical inner product
    \begin{equation}
        \label{eq:AnInnerProd}
        \SPA{x, y} = \sum_{i=1}^n x_i^*y_i
    \end{equation}
    and $\MnASP{\cdot,\cdot}$ is determined uniquely by the
    requirement \eqref{eq:SPAssoc}. Then one can show that both inner
    products are indeed completely positive, see
    \cite[Ex.~5.11]{bursztyn.waldmann:2003a:pre}. Moreover,
    $\SPA{\cdot,\cdot}$ is a non-degenerate inner product if and only
    if $\mathcal{A}$ is non-degenerate and it is full if and only if
    $\mathcal{A}$ is idempotent. Thus, under the above assumption on
    the class of $^*$-algebras we are interested in, $\mathcal{A}$ is
    strongly Morita equivalent to $M_n(\mathcal{A})$ via
    $\mathcal{A}^n$.
\end{example}
\begin{example}
    \label{example:IsoImpliesME}
    Strong Morita equivalence is implied by $^*$-isomorphism. Indeed,
    let $\Phi: \mathcal{A} \longrightarrow \mathcal{B}$ be a
    $^*$-isomorphism. Then we take $\mathcal{B}$ as a left
    $\mathcal{B}$-module in the canonical way and endow it with a
    right $\mathcal{A}$-module structure by setting $x
    \mathbin{\cdot_\Phi} a = x\Phi(a)$ for $x \in \mathcal{B}$ and $a
    \in \mathcal{A}$. For the inner products we take the canonical one
    with values in $\mathcal{B}$
    \begin{equation}
        \label{eq:BSP}
        \BSP{x, y} = xy^*
    \end{equation}
    and
    \begin{equation}
        \label{eq:SPPhiA}
        \SPA{x,y} = \Phi^{-1}(x^*y)
    \end{equation}
    for the $\mathcal{A}$-valued one. A simple verification shows that
    this gives indeed a strong Morita equivalence bimodule. Hence
    $^*$-isomorphic $^*$-algebras are strongly Morita equivalent.
\end{example}
\begin{example}
    \label{example:ccSMEBimodule}
    Let $\BEA$ be a strong Morita equivalence bimodule. Then we
    consider the complex-conjugate bimodule $\cc{\mathcal{E}}$: as
    $\ring{R}$-module it is just $\mathcal{E}$ but $\ring{C}$ acts now
    as
    \begin{equation}
        \label{eq:ccBimodule}
        \alpha \cc{x} = \cc{\cc{\alpha} x}
    \end{equation}
    where $x \mapsto \cc{x}$ is the identity map of the underlying
    $\ring{R}$-modules. Then $\cc{\mathcal{E}}$ becomes a
    $(\mathcal{A}, \mathcal{B})$-bimodule by the definitions
    \begin{equation}
        \label{eq:AccEB}
        a \cdot \cc{x} = \cc{x \cdot a^*}
        \quad
        \textrm{and}
        \quad
        \cc{x} \cdot b = \cc{b^* \cdot x}.
    \end{equation}
    Moreover, we can take the `old' inner products of $\BEA$ and
    define
    \begin{equation}
        \label{eq:ccInnerProducts}
        \ASPccE{\cc{x}, \cc{y}} = \SPEA{x,y}
        \quad
        \textrm{and}
        \quad
        \SPccEB{\cc{x}, \cc{y}} = \BSPE{x, y}.
    \end{equation}
    Then a simple conputation shows that $\AccEB$ with these inner
    products gives indeed a strong Morita equivalence $(\mathcal{A},
    \mathcal{B})$-bimodule.
\end{example}
Using Example~\ref{example:MnAAnA} for $n=1$ gives that strong Morita
equivalence is a reflexive relation, while
Example~\ref{example:ccSMEBimodule} gives symmetry. For transitivity,
we have to use again the tensor product operation $\tensM$ which can
be shown to be compatible with the other inner product as well. Thus
we finally arrive at the following statement
\cite[Thm.~5.9]{bursztyn.waldmann:2003a:pre}:
\begin{theorem}
    \label{theorem:SMEisE}
    Within the class of non-degenerate and idempotent $^*$-algebras
    strong Morita equivalence is an equivalence relation.
\end{theorem}
We shall denote the tensor product of strong equivalence bimodules by
$\tensB$ to emphasize that in this case we have to take care of
\emph{two} inner products instead of one as for $\tensM$.

Our original motivation of finding conditions for the equivalence of
the representation theories of $^*$-algebras finds now its
satisfactory answer:
\begin{theorem}
    \label{theorem:SMEEquivReps}
    If $\BEA$ is a strong Morita equivalence bimodule then the
    corresponding Rieffel induction functors
    \begin{equation}
        \label{eq:RieffelEquivalence}
        \mathsf{R}_{\mathcal{E}}: 
        \Rep[\mathcal{D}](\mathcal{A})
        \longrightarrow 
        \Rep[\mathcal{D}](\mathcal{B})
    \end{equation}
    and
    \begin{equation}
        \label{eq:RieffelEquivalenceBack}
        \mathsf{R}_{\cc{\mathcal{E}}}: 
        \Rep[\mathcal{D}](\mathcal{B})
        \longrightarrow 
        \Rep[\mathcal{D}](\mathcal{D})
    \end{equation}
    give an equivalence of categories of strongly non-degenerate
    $^*$-representations for all auxilliary $^*$-algebras
    $\mathcal{D}$.
\end{theorem}

Since it will turn out to be much easier to determine the strongly
Morita equivalent $^*$-algebras to a given $^*$-algebra $\mathcal{A}$
than understanding $\Rep[\mathcal{D}](\mathcal{A})$ itself we are now
interested in finding \emph{invariants} of strong Morita equivalence
like $\Rep[\mathcal{D}](\cdot)$.

For a detailed comparison of strong Morita equivalence with the
original definition of Rieffel \cite{rieffel:1974a, rieffel:1974b},
which contains also additional completeness requirements, we refer to
\cite{bursztyn.waldmann:2001b, bursztyn.waldmann:2003a:pre}. It turns
out that the strong Morita theory of $C^*$-algebras in Rieffel's sense
is already completely determined by the above algebraic version. Thus
it is indeed a generalization extending Rieffel's definition.

%
%

\subsection{The strong Picard Groupoid}
\label{subsec:StrongPicardGroupoid}

In order to understand strong Morita equivalence and its invariants
better, it is usefull to consider not only the question of whether
there \emph{is} a strong Morita equivalence bimodule between
$\mathcal{A}$ and $\mathcal{B}$ at all but also \emph{how many} there
may be.
\begin{definition}
    \label{definition:Pic}
    For $^*$-algebras $\mathcal{A}$, $\mathcal{B}$ we define
    $\StrPic(\mathcal{B}, \mathcal{A})$ to be the class of isometric
    isomorphism classes of strong Morita equivalence $(\mathcal{B},
    \mathcal{A})$-bimodules. We set $\StrPic(\mathcal{A}) =
    \StrPic(\mathcal{A}, \mathcal{A})$.
\end{definition}
Here isometric isomorphism classes mean isomorphic as $(\mathcal{B},
\mathcal{A})$-bimodules and isometric with respect to both inner
products.

If $\mathcal{A}$ and $\mathcal{B}$ are unital we already know that the
strong Morita equivalence bimodules are (particular) finitely
generated projective modules. Thus the class $\StrPic(\mathcal{B},
\mathcal{A})$ is a set. In the following we shall ignore the possible
logical subtleties which may arise for non-unital $^*$-algebras for
which we do not know a priori if $\StrPic(\mathcal{B}, \mathcal{A})$
is a set at all.

There are analogous definitions using $^*$-Morita equivalence or
ring-theoretic Morita equivalence yielding $\starPic(\cdot,\cdot)$ and
$\Pic(\cdot,\cdot)$, see e.g.~\cite{bass:1968a, benabou:1967a} for the
ring-theoretic version.

We have now the following structure for the collection of all
$\StrPic(\cdot,\cdot)$, see
\cite[Sect.~6.1]{bursztyn.waldmann:2003a:pre} and
\cite{waldmann:2003c:pre, waldmann:2003b:pre}:
\begin{theorem}[Strong Picard Groupoid]
    \label{theorem:PicGroupoid}
    $\StrPic(\cdot,\cdot)$ is a large Groupoid, called the strong
    Picard Groupoid, with the $^*$-algebras as units , the complex
    conjugate bimodules $[\AccEB]$ as inverses and the tensor product
    $[\CFB] [\BEA] = [\CFB \tensB[B] \BEA]$ as product.
\end{theorem}
The proof consists in showing the groupoid requirements up to
isomorphisms for the bimodules directly. A large Groupoid means that
the collection of objects is not necessarily a set. Here it is the
class of $^*$-algebras over $\ring{C}$ which are non-degenerate and
idempotent.
\begin{corollary}[Strong Picard Group]
    \label{corollary:PicGroup}
    $\StrPic(\mathcal{A})$ is a group, called the strong Picard group
    of $\mathcal{A}$. It corresponds to the isotropy group of the
    strong Picard Groupoid at the unit $\mathcal{A}$.
\end{corollary}
\begin{corollary}
    \label{corollary:PicOrbit}
    A $^*$-algebra $\mathcal{B}$ is strongly Morita equivalent to
    $\mathcal{A}$ if and only if they lie in the same orbit of
    $\StrPic$. In this case $\StrPic(\mathcal{B}) \cong
    \StrPic(\mathcal{A})$ as groups and $\StrPic(\mathcal{A})$ acts
    freely and transitively on $\StrPic(\mathcal{B}, \mathcal{A})$.
\end{corollary}
\begin{corollary}
    \label{corollary:PicMorphisms}
    There are canonical `forgetful' Groupoid morphisms
    \begin{equation}
        \label{eq:CanonicalGroupoidMorphs}
        \bfig
        \Vtriangle[\StrPic(\cdot,\cdot)`\starPic(\cdot,\cdot)`\Pic(\cdot,\cdot);(a)`(c)`(b)]
        \efig,
    \end{equation}
    such that this diagram commutes. 
\end{corollary}
\begin{remark}
    \label{remark:PicMorphsNotBij}
    In general the Groupoid morphism $(a)$ is not surjective as there
    may be more inner products (with different `signatures') on a
    $^*$-equivalence bimodule than only the positive ones. For the
    same reason, $(b)$ is not injective in general. However, even
    $(c)$ shows a non-trivial and rich behaviour: it is neither
    surjective nor injective in general.  For $C^*$-algebras it turns
    out to be always injective but not necessarily surjective. Thus we
    obtain interesting information about $\mathcal{A}$ by considering
    these Groupoid morphisms.
\end{remark}

%
%

\subsection{Actions and invariants}
\label{subsec:Actions}

The idea we want to discuss now is that strong Morita invariants can
arise from groupoid actions of $\StrPic$ on `something'. Then
`something' is preserved along the orbits of the Groupoid $\StrPic$,
i.e. the strong Morita equivalence classes of $^*$-algebras.
This is of course more a philosophical statement than a theorem and we
do not want to make any attempt to make this proposal
precise. However, we can illustrate this principle by several examples
following \cite{waldmann:2003c:pre}:
\begin{example}[Strong Picard groups]
    \label{example:PicInvariant}
    The strong Picard group $\StrPic(\mathcal{A})$ is a strong Morita
    invariant. Indeed, $\StrPic$ acts on itself so the isotropy groups
    are all isomorphic along an orbit. Any element in
    $\StrPic(\mathcal{B}, \mathcal{A})$ provides then a group
    isomorphism between $\StrPic(\mathcal{A})$ and
    $\StrPic(\mathcal{B})$. This is in some sense the most fundamental
    Morita invariant.
\end{example}
\begin{example}[Hermitian $K_0$-groups]
    \label{example:KnullInvariant}
    Recall that the \emph{Hermitian $K_0$-group} $K_0^H(\mathcal{A})$
    of a unital $^*$-algebra $\mathcal{A}$ is defined as follows: one
    considers finitely generated projective modules with strongly
    non-degenerate and completely positive inner products
    $\SPA{\cdot,\cdot}$. We can take direct orthogonal sums without
    loosing these properties so taking isometric isomorphism classes
    gives us an (abelian) semigroup with respect to $\oplus$. Then
    $K^H_0(\mathcal{A})$ is defined as the Grothendieck group of this
    semigroup.

    Now if $\FB$ is such a finitely generated projective module and
    $\BEA$ is a strong Morita equivalence bimodule then $\FB \tensM[B]
    \BEA$ is again finitely generate projective and the
    $\mathcal{A}$-valued inner product is still strongly
    non-degenerate. Moreover, $\tensM$ is clearly compatible with
    direct orthogonal sums. Thus, by passing to isometric isomorphism
    classes, one obtains an action of $\StrPic$ on $K^H_0$
    \begin{equation}
        \label{eq:PicOnK}
        K_0^H(\mathcal{B}) \times \StrPic(\mathcal{B}, \mathcal{A})
        \longrightarrow K_0^H(\mathcal{A})
    \end{equation}
    by group isomorphisms. This has the following consequences: First,
    the strong Picard group $\StrPic(\mathcal{A})$ acts by group
    isomorphisms on the abelian group $K^H_0(\mathcal{A})$. Second,
    $K^H_0(\mathcal{A})$ is a strong Morita invariant, even as a
    $\StrPic(\mathcal{A})$-space.
\end{example}
\begin{example}[Representation theories]
    \label{example:RepInvariant}
    The strong Picard Groupoid acts on the representation theories
    $\Rep[\mathcal{D}](\cdot)$ by Rieffel induction
    \begin{equation}
        \label{eq:PicActsOnRep}
        \StrPic(\mathcal{B}, \mathcal{A}) \times
        \Rep[\mathcal{D}](\mathcal{A})
        \longrightarrow
        \Rep[\mathcal{D}](\mathcal{B}).
    \end{equation}
    However, this is not an honest action as the action properties are
    only fulfilled up to unitary equivalences of representations. Thus
    this should better be interpreted as an `action' of the strong
    Picard \emph{bigroupoid} on the categories
    $\Rep[\mathcal{D}](\cdot)$, where the strong Picard bigroupoid
    consists of all equivalence bimodules \emph{without} identifying
    them up to isometric isomorphisms. Since it would require 10
    additional pages of commutative diagrams to give a definition of
    what the action of a bigroupoid should be, we do not want to make
    this more precise but leave it as a heuristic example to challenge
    the imagination of the reader, see also \cite{benabou:1967a}.
    Another option is to consider the unitary equivalence classes of
    $^*$-representations instead of $\Rep[\mathcal{D}](\cdot)$: Then
    the Picard groupoid acts by Rieffel induction in a well-defined
    way.
\end{example}

There are many more examples of strong Morita invariants like the
centers of $^*$-algebras or their lattices of closed $^*$-ideals in
the sense of \cite{bursztyn.waldmann:2001b}. Thus it is interesting to
see whether one can view all strong Morita invariants as arising of an
appropriate action of the strong Picard Groupoid:
\begin{question}
    \label{question:AllInvariants}
    Can one view any strong Morita invariant as coming from an action
    of the strong Picard Groupoid?
\end{question}
Probably it becomes tautological if one formulates this in the
appropriate context. Nevertheless, a consequence of a positive answer
would be that any strong Morita invariant carries an action of the
Picard group which is invariant itself.

It is clear that also in the ring-theoretic framework as well as for
$^*$-Morita equivalence one can pose the same question. In fact, some
of the above strong Morita invariants have their immediate and
well-known analogs for these coarser notions of Morita equivalence.

%
%

\subsection{Strong vs. ring-theoretic Morita equivalence}
\label{subsec:SMEvsME}

Let us now discuss the relation between strong Morita equivalence and
ring-theoretic Morita equivalence more closely. For simplicity, we
shall focus on \emph{unital} $^*$-algebras throughout this section.
Then it is clear that strong Morita equivalence implies Morita
equivalence since we have even a groupoid morphism
\begin{equation}
    \label{eq:StrPicPic}
    \StrPic \longrightarrow \Pic.
\end{equation}
Thus if $\StrPic(\mathcal{B}, \mathcal{A})$ is non-empty then the
image under \eqref{eq:StrPicPic} is non-empty as well. To understand
the relation between strong and ring-theoretic Morita equivalence
better we want to understand the kernel and the image of the groupoid
morphism \eqref{eq:StrPicPic}.

Thus we first have to determine the structure of strong equivalence
bimodules as precise as possible. The following proposition gives a
simple proof of the well-known fact that equivalence bimodules are
finitely generated and projective using the inner products of a strong
equivalence bimodule:
\begin{proposition}
    \label{proposition:FGPM}
    Let $\BEA$ be a strong Morita equivalence bimodule. Then there
    exist $\xi_i, \eta_i \in \BEA$, $i = 1, \ldots, n$, such that
    \begin{equation}
        \label{eq:xisx}
        x = \sum_{i=1}^n \xi_i \cdot \SPEA{\eta_i, x}.
    \end{equation}
    It follows that $\BEA$ is finitely generated and projective as
    right $\mathcal{A}$-module and by symmetry the same statement
    holds for $\mathcal{B}$.
\end{proposition}
\begin{proof}
    Indeed, let $\Unit_{\mathcal{B}} = \sum_{i=1}^n \BSPE{\xi_i,
      \eta_i}$ by fullness. Then the compatibility of the inner
    products gives
    \[
    x = \Unit_{\mathcal{B}} \cdot x = \sum_{i=1}^n \BSPE{\xi_i,
      \eta_i} \cdot x
    = \sum_{i=1}^n \xi_i \cdot \SPEA{\eta_i, x}.
    \]
    Since $\SPEA{\cdot, \cdot}$ is $\mathcal{A}$-linear to the right
    in the second argument, it follows that the $\xi_i$ together with
    the functionals $\SPEA{\eta_i, \cdot}$ form a finite dual basis.
    By the dual basis lemma, see e.g. \cite[Lemma~(2.9)]{lam:1999a},
    this is equivalent to the fact that $\BEA$ is finitely generatey
    (by the generators $\xi_i$), and projective.
\end{proof}

We shall call the $\{\xi_i, \eta_i\}_{i=1, \ldots, n}$ with the above
property a \emph{Hermitian dual basis}. Thus we have
\begin{equation}
    \label{eq:EAFGPM}
    \EA \cong e\mathcal{A}^n
\end{equation}
for some idempotent $e = e^2 \in M_n(\mathcal{A})$. In fact, $e$ can
be expressed in terms of the Hermitian dual basis explicitly by
\begin{equation}
    \label{eq:eisxieta}
    e = (e_{ij})
    \quad
    \textrm{with}
    \quad
    e_{ij} = \SPEA{\eta_i, \xi_i},
\end{equation}
and the isomorphism \eqref{eq:EAFGPM} is simply given by
\begin{equation}
    \label{eq:EaIsoeAn}
    \EA \ni x
    \mapsto
    \left(\SPEA{\eta_i, x}\right)_{i = 1, \ldots, n} 
    \in e \mathcal{A}^n \subseteq \mathcal{A}^n.
\end{equation}
In particular, it follows that the inner products on a strong Morita
equivalence bimodule are always strongly non-degenerate, see
Definition~\ref{definition:PreHilbertModule}, in the case of unital
$^*$-algebras.

Note however, that we can \emph{not} conclude that $e$ can be chosen
to be a Hermitian idempotent, i.e. a projection. Thus the question how
the inner product $\SPEA{\cdot,\cdot}$ looks like when we identify
$\BEA$ with $e\mathcal{A}^n$ is difficult to answer: How many
completely positive, full and non-degenerate $\mathcal{A}$-valued
inner products can we have on $e\mathcal{A}^n$ up to isometries?  In
order to be able to say something one has to assume additional
properties of the $^*$-algebras in question. Motivated by the case of
$C^*$-algebras we state the following conditions:
\begin{description}
\item[(I)] For all $n \in \mathbb{N}$ and $A \in M_n(\mathcal{A})$ the
    element $\Unit + A^*A$ is invertible.
\end{description}
In particular, since we require this condition for all $n$ we also
have the invertibility of $\Unit + A_1^*A_1 + \cdots + A_k^*A_k$ for
$A_1, \ldots, A_k \in M_n(\mathcal{A})$.  The relevance of this
condition \cond{I} is classical, see Kaplansky's book
\cite[Thm.~26]{kaplansky:1968a}:
\begin{lemma}
    \label{lemma:Kaplansky}
    Assume that $\mathcal{A}$ satisfies \cond{I}. Then for any
    idempotent $e = e^2 \in M_n(\mathcal{A})$ there exists a
    projection $P = P^2 = P^* \in M_n(\mathcal{A})$ and $u, v \in
    M_n(\mathcal{A})$ with
    \begin{equation}
        \label{eq:Pequive}
        e = uv
        \quad
        \textrm{and}
        \quad
        P = vu,
    \end{equation}
    whence the projective modules $e\mathcal{A}^n$ and
    $P\mathcal{A}^n$ are isomorphic via $v$ and $u$.
\end{lemma}
Thus having the property \cond{I} we can assume for any finitely
generated projective module $e \mathcal{A}^n \cong P\mathcal{A}^n$
with some projection instead of a general idempotent. On
$P\mathcal{A}^n$ there is the restriction of the canonical inner
product $\SP{\cdot,\cdot}$ of $\mathcal{A}^n$ such that
\begin{equation}
    \label{eq:BoundedPAnIsPMnP}
    \mathfrak{B}(P\mathcal{A}^n, \SP{\cdot,\cdot})
    \cong P M_n(\mathcal{A})P
\end{equation}
as $^*$-algebras, since $P = P^*$.

The next question is how many other inner products of interest does
one have on $P\mathcal{A}^n$? The following technical condition will
guarantee that there is only one up to isometric isomorphisms. Again,
$C^*$-algebras are the motivation for this condition:
\begin{description}
\item[(II)] Let $P_\alpha \in M_n(\mathcal{A})$ be finitely many
    pairwise orthgonal projections $P_\alpha P_\beta =
    \delta_{\alpha\beta} P_\alpha = \delta_{\alpha\beta} P_\alpha^*$
    such that $\sum_\alpha P_\alpha = \Unit$ and let $H \in
    M_n(\mathcal{A})^+$ be invertible. If $[H, P_\alpha] = 0$ then
    there exists an invertible $U$ (depending on the $P_\alpha$ and on
    $H$) such that $H = U^*U$ and $[U, P_\alpha] = 0$.
\end{description}
This mimicks in some sense the spectral calculus for matrices and for
$C^*$-algebras this is obviously fulfilled since here for any positive
$H$ one even has a unique positive square root $\sqrt{H}$ which
commutes with all elements commuting with $H$.

Assume that $\mathcal{A}$ satisfies \cond{II} and let $h:
P\mathcal{A}^n \times P\mathcal{A}^n \longrightarrow \mathcal{A}$ be a
completely positive and strongly non-degenerate inner product. Then we
can extend $h$ to $\mathcal{A}^n$ by using e.g. the restriction of the
canonical inner product $\SP{\cdot,\cdot}$ on $(\Unit -
P)\mathcal{A}^n$. The result is a completely positive and strongly
non-degenerate inner product on the free module $\mathcal{A}^n$ which
we denote by $\hat{h}(\cdot,\cdot)$.  Then we define the matrix $H \in
M_n(\mathcal{A})$ by
\begin{equation}
    \label{eq:HfromhandSPDef}
    H_{ij} = \hat{h}(e_i, e_j) = h(Pe_i, Pe_j) 
    + \SP{(\Unit - P)e_i, (\Unit - P)e_j},
\end{equation}
whence
\begin{equation}
    \label{eq:hByH}
    \hat{h}(x, y) = \SP{x, Hy}
\end{equation}
for all $x, y \in \mathcal{A}^n$.  Since $\hat{h}$ is completely
positive $H$ is a positive matrix and since $\hat{h}$ is strongly
non-degenerate one finds that $H$ is invertible. Moreover, it is clear
that $[P, H] = 0$.  Thus we can apply \cond{II} and find an invertible
$U \in M_n(\mathcal{A})$ with $H = U^*U$ and $[P, U] = 0$. Thus
\begin{equation}
    \label{eq:hathIsometricSP}
    \hat{h}(x, y) = \SP{x, Hy} = \SP{Ux, Uy}
\end{equation}
whence $\hat{h}$ is isometric to the canonical inner product
$\SP{\cdot,\cdot}$. Since $[P, U] = 0$ the isometry $U$ restricts
to the projective module $P\mathcal{A}^n$ and gives an isometric
isomorphism between $h$ and the restriction of $\SP{\cdot,\cdot}$.
\begin{lemma}
    \label{lemma:NiceUniqueInnerProducts}
    Assume $\mathcal{A}$ satisfies \cond{II} and let $P = P^* = P^2 \in
    M_n(\mathcal{A})$. Then any two completely positive and strongly
    non-degenerate inner product on $P\mathcal{A}^n$ are isometric.
\end{lemma}
Combining both properties leads to the following characterizations of
strong Morita equivalence bimodules:
\begin{theorem}
    \label{theorem:SMEBimodulesForIandII}
    Let $\mathcal{A}$ and $\mathcal{B}$ be unital $^*$-algebras and
    assume $\mathcal{A}$ satisfies \cond{I} and \cond{II}. If $\BEA$
    is a $^*$-Morita equivalence bimodule with completely positive
    inner product $\SPEA{\cdot,\cdot}$ then we have:
    \begin{compactenum}
    \item There exists a full projection $P \in M_n(\mathcal{A})$ such
        that $\EA \cong P\mathcal{A}^n$ are isometrically isomorphic.
    \item The left action of $\mathcal{B}$ on $\BEA$ and the above
        isomorphism induce a $^*$-isomorphism $B \cong P
        M_n(\mathcal{A})P$ and under this isomorphism
        $\BSPE{\cdot,\cdot}$ becomes the canonical
        $PM_n(\mathcal{A})P$-valued inner product on $P\mathcal{A}^n$.
    \item $\BSPE{\cdot,\cdot}$ is necessarily completely positive,
        too, whence $\BEA$ is already a strong Morita equivalence
        bimodule.
    \end{compactenum}
    Conversely, any full projection $P \in M_n(\mathcal{A})$ gives a
    strong Morita equivalence bimodule $P\mathcal{A}^n$ between
    $\mathcal{A}$ and $PM_n(\mathcal{A})P$.
\end{theorem}
The fullness of the projection $P$ is equivalent to the statement that
the canonical inner product on $P\mathcal{A}^n$ is full.

One easily obtains the following consequences of this theorem:
\begin{theorem}
    \label{theorem:CondIIISME}
    The conditions \cond{I} and \cond{II} together are strongly Morita
    invariant.
\end{theorem}
To see this, we only have to check it for $PM_n(\mathcal{A})P$ by hand
which is straightforward.
\begin{theorem}
    \label{theorem:PisStrPicInj}
    Within the class of unital $^*$-algebras satisfying \cond{I} and
    \cond{II} the groupoid morphism
    \begin{equation}
        \label{eq:PicStrToPic}
        \StrPic \longrightarrow \Pic
    \end{equation}
    is injective (though not necessarily surjective).
\end{theorem}
This is also clear since on a Morita equivalence bimodule we can have
at most one inner product up to isometric isomorphisms according to
Theorem~\ref{theorem:SMEBimodulesForIandII}. This also implies the
following result for general finitely generated projective modules:
\begin{theorem}
    \label{theorem:HermKisK}
    For a unital $^*$-algebra $\mathcal{A}$ satisfying \cond{I} and
    \cond{II} we have canonically
    \begin{equation}
        \label{eq:KHnullIsKnull}
        K_0^H(\mathcal{A}) \cong K_0(\mathcal{A}).
    \end{equation}
\end{theorem}

The question of surjectivity of \eqref{eq:PicStrToPic} is actually
more subtle.  Here we have to impose first another condition on the
$^*$-algebras we consider. The condition is not on a single
$^*$-algebra but on a whole family of $^*$-algebras under
considerations:
\begin{description}
\item[(III)] Let $\mathcal{A}$ and $\mathcal{B}$ be unital
    $^*$-algebras and let $P \in M_n(\mathcal{A})$ be a projection and
    consider the $^*$-algebra $P M_n(\mathcal{A})P$. If $\mathcal{B}$
    and $P M_n(\mathcal{A}) P$ are isomorphic as unital algebras then
    they are also $^*$-isomorphic.
\end{description}
In fact, for unital $C^*$-algebras this is always fulfilled as in this
case $PM_n(\mathcal{A})P$ is a $C^*$-algebra again and thus the
$^*$-involution is uniquely determined, see
\cite[Thm.~4.1.20]{sakai:1971a}.  Another class of $^*$-algebras
satisfying this condition are the Hermitian star products. In fact, if
$\star$ is a star product having a $^*$-involution of the form $f
\mapsto \cc{f} + o(\lambda)$ then it is $^*$-equivalent to a Hermitian
star product, see \cite[Lem.~5]{bursztyn.waldmann:2002a}.

Now consider the automorphism group $\Aut(\mathcal{B})$ of
$\mathcal{B}$ then $\Aut(\mathcal{B})$ acts from the left on the set
$\Pic(\mathcal{B}, \mathcal{A})$ in the following way, see also
Example~\ref{example:IsoImpliesME}. The left $\mathcal{B}$-module
structure of $\BEA$ is twisted by $\Phi$ as
\begin{equation}
    \label{eq:PhiTwistBimodule}
    b \cdot_\Phi x = \Phi^{-1}(b) \cdot x
\end{equation}
while the right $\mathcal{A}$-module structure is untouched. This
gives again an equivalence bimodule, now denoted by
$\Bimod{\Phi}{}{\mathcal{E}}{}{}$. It can be checked easily that this
descends to a group action of $\Aut(\mathcal{B})$ on
$\Pic(\mathcal{B}, \mathcal{A})$, see e.g. the discussion in
\cite{bass:1968a, bursztyn.waldmann:2004a}. The problem with the
surjectivity is that for a given ring-theoretic equivalence bimodule
we may obtain the `wrong' $^*$-involution induced for $\mathcal{B}$:
\begin{proposition}
    \label{proposition:StrPicPicSurjectivity}
    Let $\mathcal{A}$, $\mathcal{B}$ satisfy condition \cond{III} and
    let $\mathcal{A}$ satisfy \cond{I} and \cond{II}. Then the
    map
    \begin{equation}
        \label{eq:StrPicPicOnto}
        \StrPic(\mathcal{B}, \mathcal{A}) \longrightarrow
        \Pic(\mathcal{B},\mathcal{A}) \big/ \Aut(\mathcal{B})
    \end{equation}
    is onto.
\end{proposition}
There are some immediate consequences when we apply this to the
previous examples like $C^*$-algebras and star products:
\begin{corollary}
    \label{corollary:BeerResult}
    Within a class of unital $^*$-algebras satisfying \cond{I},
    \cond{II} and \cond{III} ring-theoretic Morita equivalence implies
    strong Morita equivalence.
\end{corollary}
In fact, for $C^*$-algebras this is Beer's theorem \cite{beer:1982a}
while for star products this was obtained in
\cite{bursztyn.waldmann:2002a}, see
Corollary~\ref{corollary:InjectivePic}.

The obstruction whether \eqref{eq:PicStrToPic} is onto and not only
onto up to automorphisms can be encoded in a particular property of
the automorphism group of the algebras. We state the last condition:
\begin{description}
\item[(IV)] For any $\Phi \in \Aut(\mathcal{A})$ there is an
    invertible $U \in \mathcal{A}$ such that $\Phi^* \Phi^{-1} =
    \Ad(U^*U)$ where $\Phi^*(a) = \Phi(a^*)^*$.
\end{description}
In particular, if $\Phi$ is even a $^*$-automorphism then $\Phi^* =
\Phi$ whence the condition is trivially fulfilled for those. So the
condition says that those automorphisms which are \emph{not}
$^*$-automorphisms have to be `essentially inner'.
\begin{theorem}
    \label{theorem:SurjectivityofStrPicPic}
    Consider $^*$-algebras $\mathcal{A}$, $\mathcal{B}$ in a class of
    unital $^*$-algebras satisfying \cond{I}, \cond{II} and
    \cond{III}.
    \begin{compactenum}
    \item $\StrPic(\mathcal{B}, \mathcal{A}) \longrightarrow
        \Pic(\mathcal{B}, \mathcal{A})$ is surjective if and only if
        $\mathcal{B}$ satisfies \cond{IV}.
    \item Within this class, \cond{IV} is strongly Morita invariant.
    \end{compactenum}
\end{theorem}

The condition \cond{IV} captures very interesting properties of the
automorphism group. One can find explicit examples of $C^*$-algebras
where \cond{IV} is not satisfied. Moreover, for $^*$-products, the
automorphisms in question can be written as exponentials of
derivations, see \cite[Prop.~8.8]{bursztyn.waldmann:2003a:pre} whence
one arrives at the question whether certain derivations are inner or
not, see\cite[Thm.~8.9]{bursztyn.waldmann:2003a:pre}:
\begin{example}[Aharonov-Bohm Effect]
    \label{example:CondIVForStarproducts}
    Let $\star$, $\star'$ be strongly Morita equivalent star products
    on $M$. Then $\StrPic(\star', \star) \longrightarrow \Pic(\star',
    \star)$ is bijective if and only if all derivations of $\star$ are
    quasi-inner, i.e. of the form $D = \frac{\I}{\lambda} \ad(H)$ with
    some $H \in C^\infty(M)[[\lambda]]$. In particular, if $M$ is
    symplectic, then this is the case if and only if $\HdR^1(M,
    \mathbb{C}) = \{0\}$. On the other hand, as argued in
    \cite{bordemann.neumaier.pflaum.waldmann:2003a} for the case of
    cotangent bundles and more generally in
    \cite{bursztyn.waldmann:2004a}, the first deRham cohomology is
    responsible for Aharonov-Bohm like effects in deformation
    quantization. Thus the question of surjectivity of
    \eqref{eq:PicStrToPic} gets the physical interpretation of the
    question whether there are Aharonov-Bohm effects possible or not.
\end{example}

%
%

\section{(Strong) Morita equivalence of star products}
\label{sec:SMEStarProducts}

We shall now apply the results of the last section to Hermitian
deformations of $^*$-algebras in order to investigate their strong
Morita theory. First we have to discuss how the conditions \cond{I}
and \cond{II} behave under deformations, in particular for the case of
star products. Again, in this section all $^*$-algebras are assumed to
be unital.

%
%

\subsection{Deformed $^*$-algebras}
\label{subsec:DeformedStarAlgebras}

We consider a Hermitian deformation $\boldsymbol{\mathcal{A}} =
(\mathcal{A}[[\lambda]], \star)$ of a unital $^*$-algebra
$\mathcal{A}$. Then the following observation is trivial:
\begin{lemma}
    \label{lemma:boldAIiffAI}
    The $^*$-algebra $\boldsymbol{\mathcal{A}}$ satisfies \cond{I} if
    and only if $\mathcal{A}$ satisfies \cond{I}.
\end{lemma}
Thus the condition \cond{I} is \emph{rigid} under Hermitian
deformations. More subtle and surprising is the following rigidity
statement:
\begin{theorem}
    \label{theorem:PropIIRigid}
    The condition \cond{II} is rigid under completely positive
    deformations.
\end{theorem}
The idea of the proof is to consider an invertible $\boldsymbol{H} \in
M_n(\boldsymbol{\mathcal{A}})^+$ whence it's classical limit $H_0 \in
M_n(\mathcal{A})^+$ is still invertible and positive according to
Corollary~\ref{corollary:ClassLimPosEl}, adapted to this more general
formulation. Then if $\boldsymbol{P}_\alpha$ are pairwise
$\star$-commuting projections commuting with $\boldsymbol{H}$ then
their classical limits $P_\alpha$ are pairwise commuting projections
commuting with $H_0$ whence we can apply \cond{II} for the classical
limit $\mathcal{A}$ and find a $U_0$. Then the idea is to lift $U_0$
in an appropriate way to find $\boldsymbol{U}$ with $\boldsymbol{H} =
\boldsymbol{U}^*\boldsymbol{U}$ and $[\boldsymbol{P}_\alpha,
\boldsymbol{U}] = 0$.

This rigidity is very nice since star products are completely positive
deformations according to Theorem~\ref{theorem:ComPosDef}. Thus we
have the following consequences:
\begin{corollary}
    \label{corollary:StarsIandII}
    Hermitian star products satisfy \cond{I} and \cond{II}.
\end{corollary}
\begin{corollary}
    \label{corollary:InjectivePic}
    Hermitian star products are strongly Morita equivalent if and only
    if they are Morita equivalent. Moreover, the groupoid morphism
    \begin{equation}
        \label{eq:PicStrStarInj}
        \StrPic(\star', \star) \longrightarrow  \Pic(\star', \star)
    \end{equation}
    is injective.
\end{corollary}
Thus we only have to understand the ring-theoretic Morita equivalence
of star products to get the strong Morita equivalence for free. Note
however, that the (non-)surjectivity of \eqref{eq:PicStrStarInj}
depends very much on the star products under consideration.

%
%

\subsection{Deformed projections}
\label{subsec:DeformedProjections}

We shall now simplify our discussion to the ring-theoretic Morita
equivalence as for star products this is sufficient thanks to
Corollary~\ref{corollary:InjectivePic}.

For a given (Hermitian) deformation $\boldsymbol{\mathcal{A}} =
(\mathcal{A}[[\lambda]], \star)$ we have to find the full idempotents
$\boldsymbol{P} \in M_n(\boldsymbol{\mathcal{A}})$ in order to find
all other algebras $\boldsymbol{\mathcal{B}}$ which are Morita
equivalent to $\boldsymbol{\mathcal{A}}$ since then
\begin{equation}
    \label{eq:BequivBoldMnA}
    \boldsymbol{\mathcal{B}} 
    = \boldsymbol{P} \star 
    M_n(\boldsymbol{\mathcal{A}}) \star \boldsymbol{P}
\end{equation}
gives all Morita equivalent algebras up to isomorphism.  Thus we have
to investigate the idempotents in $M_n(\boldsymbol{\mathcal{A}})$.
From the defining equation $\boldsymbol{P} \star \boldsymbol{P} =
\boldsymbol{P}$ we find in zeroth order
\begin{equation}
    \label{eq:clBoldPisP}
    P_0 \cdot P_0 = P_0,
\end{equation}
where $\boldsymbol{P} = \sum_{r=0}^\infty \lambda^r P_r$. Thus the
classical limit of an idempotent $\boldsymbol{P}$ is an idempotent
$P_0$ for the undeformed product.
\begin{lemma}
    \label{lemma:PFullIffPnullFull}
    $\boldsymbol{P}$ is full if and only if $P_0$ is full.
\end{lemma}
This is straightforward to show. The next observation is more
non-trivial and can be found implicitly in e.g.
\cite{gerstenhaber.schack:1988a, emmrich.weinstein:1996a} while the
explicit formula can be found in \cite[Eq.~(6.1.4)]{fedosov:1996a}: If
$P_0 \in M_n(\mathcal{A})$ is an idempotent with respect to the
undeformed product then
\begin{equation}
    \label{eq:FedosovProjektion}
    \boldsymbol{P} = \frac{1}{2} + \left(P_0 - \frac{1}{2}\right)
    \star \frac{1}{\sqrt[\star]{\Unit + 4 (P_0 \star P_0 - P_0)}}
\end{equation}
defines an idempotent $\boldsymbol{P} \in
M_n(\boldsymbol{\mathcal{A}})$ with respect to $\star$. Here we have
to assume $\mathbb{Q} \subseteq \ring{R}$ in order to make the series
well-defined. Note that as a formal power series in $\lambda$ the
right hand side in \eqref{eq:FedosovProjektion} is well-defined since
in zeroth order $P_0 \star P_0 - P_0$ vanishes. Moreover, the
classical limit of this $\boldsymbol{P}$ coincides with $P_0$ and if
$P_0$ is already an idempotent with respect to $\star$, then the
formula reproduces $P_0$. Finally, if $\star$ is a Hermitian
deformation of a $^*$-algebra and if $P_0$ is even a projection,
i.e. $P_0^* = P_0$, then $\boldsymbol{P}$ is also a projection.

The next statement concerns about the uniqueness of the deformation
$\boldsymbol{P}$ of a given projection $P_0$. First recall that two
idempotents $P$ and $Q$ are called equivalent if (after embedding into
some big matrix algebra $M_n(\mathcal{A})$) there exist $U$, $V$ such
that $P = UV$ and $Q = VU$. This is the case if and only if the
corresponding projective modules $P \mathcal{A}^n$ and
$Q\mathcal{A}^n$ are isomorphic as right $\mathcal{A}$-modules. In
fact, $U$ and $V$ provide such mutually inverse module isomorphisms
when restricted to $P \mathcal{A}^n$ and $Q\mathcal{A}^n$, see also
Lemma~\ref{lemma:Kaplansky}. For projections one has a slightly
refined notion, namely one demands $P = U^*U$ and $Q = UU^*$. For the
deformations we have now the following statement:
\begin{lemma}
    \label{lemma:EquivalentProjections}
    Two deformed idempotents $\boldsymbol{P}$ and $\boldsymbol{Q}$ are
    equivalent if and only if their classical limits $P_0$ and $Q_0$
    are equivalent.
\end{lemma}
As a consequence one immediately obtains the \emph{rigidity} of the
$K_0$-group under formal deformations, i.e. the classical limit map
induces an isomorphism
\begin{equation}
    \label{eq:RigidityKnull}
    \cl_*: K_0 (\boldsymbol{\mathcal{A}}) 
    \stackrel{\cong}{\longrightarrow} 
    K_0 (\mathcal{A}),
\end{equation}
see \cite{rosenberg:1996a:pre}. One can also show that as
$\ring{C}[[\lambda]]$-modules we have
\begin{equation}
    \label{eq:EndProjModDeformed}
    (P_0 M_n(\mathcal{A}) P_0)[[\lambda]] \cong 
    \boldsymbol{P} 
    \star M_n(\boldsymbol{\mathcal{A}}) 
    \star \boldsymbol{P}
\end{equation}
by an isomorphism with the identity as classical limit, when we view
both spaces as submodules of $M_n(\mathcal{A})[[\lambda]]$.

With these results we have the following picture: Given a Morita
equivalence bimodule $\BEA \cong P_0\mathcal{A}^n$ we know $\mathcal{B}
\cong P_0 M_n(\mathcal{A}) P_0$. Moreover, let a deformation $\star$
of $\mathcal{A}$ be given. Then any choice of an
$\ring{C}[[\lambda]]$-isomorphism as in \eqref{eq:EndProjModDeformed}
gives an isomorphism
\begin{equation}
    \label{eq:BPMNAP}
    \mathcal{B}[[\lambda]] \cong 
    \boldsymbol{P} 
    \star M_n(\boldsymbol{\mathcal{A}}) 
    \star \boldsymbol{P}
\end{equation}
as $\ring{C}[[\lambda]]$-submodules of
$M_n(\boldsymbol{A})[[\lambda]]$ inducing the identity in zeroth
order. Since the right hand side carries a algebra structure this
induces a new associative multiplication $\star'$ for
$\mathcal{B}[[\lambda]]$ which turns out to be a deformation of
$\mathcal{B}$. Since everything is unique up to isomorphisms and since
the isomorphisms can be adapted in such a way that they induce the
identity in zeroth order we find the following:
\begin{compactenum}
\item $\star'$ is unique up to equivalence.
\item $(\mathcal{B}[[\lambda]], \star')$ is Morita equivalent to
    $(\mathcal{A}[[\lambda]], \star)$.
\end{compactenum}
Furthermore, everything depends only on the isomorphism class of the
equivalence bimodule and behaves nicely under tensor products.
Denoting by
\begin{equation}
    \label{eq:DefADef}
    \Def(\mathcal{A}) = \{\textrm{equivalence classes of formal
      associative deformations of}\; \mathcal{A}\}
\end{equation}
the \emph{deformation theory} of $\mathcal{A}$ we have an action
\begin{equation}
    \label{eq:PicActsOfDef}
    \Pic(\mathcal{B}, \mathcal{A}) \times \Def(\mathcal{A})
    \ni ([\mathcal{E}], [\star]) 
    \mapsto
    [\star'] \in \Def(\mathcal{B})
\end{equation}
of the Picard groupoid of the undeformed algebras on their deformation
theories $\Def(\cdot)$ such that $\star'$ gives a Morita equivalent
deformation to $\star$ if and only if $[\star']$ and $[\star]$ lie in
the same orbit of the groupoid action \eqref{eq:PicActsOfDef}, see
\cite{bursztyn:2002a}.
\begin{remark}
    \label{remark:YetAnotherMInvariant}
    With the deformation theories we just found another Morita
    invariant coming from an action of the Picard groupoid, here in
    the ring-theoretic framework, see also
    \cite{gerstenhaber.schack:1988a}.
\end{remark}

%
%

\subsection{Morita equivalent star products}
\label{subsec:MoritaEquivStarProducts}

Now we want to apply these general results to star product
algebras. Thus we first have to identify the classical Picard groupoid
and then determine its action on the deformation theories by examining
the projective modules.

The first task has a well-known solution. The Picard groupoid for
algebras of smooth functions $C^\infty(M)$ is given as follows:
\begin{equation}
    \label{eq:PicMMprime}
    \Pic(C^\infty(M), C^\infty(M')) = \emptyset
    \quad
    \textrm{for}
    \quad
    M \not\cong M'
\end{equation}
\begin{equation}
    \label{eq:PicCinftyM}
    \Pic(C^\infty(M)) 
    = \group{Diff}(M) \ltimes \mathrm{H}^2(M, \mathbb{Z})
\end{equation}
Note that this determines $\Pic(\cdot,\cdot)$ completely. Here
$\group{Diff}(M) = \Aut(C^\infty(M))$ is the diffeomorphism group of
$M$ which twists the bimodules in the usual way.  The equivalence
bimodules where $C^\infty(M)$ acts the same way from left and right
are just the sections $\Gamma^\infty(L)$ of complex line bundles $L
\longrightarrow M$. They are well-known to be classified by the second
integer cohomology of $M$.

For the second step we have to deform the sections $\Gamma^\infty(L)$
such that $(C^\infty(M)[[\lambda]], \star)$ acts from the right on
$\Gamma^\infty(M)[[\lambda]]$ while $(C^\infty(M)[[\lambda]], \star')$
acts from the left and both actions commute. We have to compute the
characteristic class of $\star'$ in terms of the equivalence class of
$\star$ and $[L] \in \mathrm{H}^2(M, \mathbb{Z})$. The result for the
symplectic case has a very appealing formulation using the
characteristic classes of star products, see
\cite{bursztyn.waldmann:2002a}:
\begin{theorem}
    \label{theorem:MEofStarproducts}
    In the symplectic case we have
    \begin{equation}
        \label{eq:cstarprimecstarL}
        c(\star') = c(\star) + 2 \pi \I c_1(L),
    \end{equation}
    where
    \begin{equation}
        \label{eq:charclass}
        c(\star) \in \frac{[\omega]}{\I\lambda} 
        + \HdR^2(M, \mathbb{C})[[\lambda]]
    \end{equation}
    is the characteristic class of $\star$ and $c_1(L) \in \HdR^2(M,
    \mathbb{C})$ is the Chern class of $L$.
\end{theorem}
\begin{remark}[Morita Equivalent Star Products]
    \label{remark:MoritaEquvialenceStarProducts}
    From this one immediately obtains the full classification of
    Morita equivalent star products in the symplectic case as we only
    have to re-implement the automorphisms from $\group{Diff}(M)$. The
    final answer is therefor that $\star$ and $\star'$ are Morita
    equivalent symplectic star products on $(M, \omega)$ if and only
    if there exists a symplectomorphism $\psi$ such that
    \begin{equation}
        \label{eq:starstarME}
        \psi^* c(\star') - c(\star) 
        \in 2 \pi \I \HdR^2(M, \mathbb{Z}).
    \end{equation}
\end{remark}
\begin{remark}
    \label{remark:MEPoissonCase}
    Similar results hold in the Poisson case where
    \eqref{eq:cstarprimecstarL} has to be formally inverted to give
    the formal deformations of the Poisson tensor which classifies the
    star products according to the formality theorem,
    see \cite{jurco.schupp.wess:2002a} for a discussion.
\end{remark}
\begin{remark}[Dirac's Monopole]
    \label{remark:DiracMonopoleAgain}
    From this point of view the results on the Dirac monopole as in
    Section~\ref{subsec:DiracMonopole} are much more transparent:
    Dirac's quantization condition for the monopole charge of a
    magnetic monopole is precisely the integrality condition for the
    two-form $B$ to define a line bundle. But this is the condition
    for the Morita equivalence of the two quantizations. The Rieffel
    induction functor is then just the induction by the equivalence
    bimodule, see also \cite{bursztyn.waldmann:2002a}.
\end{remark}

The proof of the theorem consist in constructing local (even
bidifferential) bimodule multiplications $\bullet$ and $\bullet'$
which allow to construct a deformed version of the transition
functions of the line bundle. Now these deformed transition functions
obey a cocycle identity with respect to the star product $\star$. From
this one can reconstruct the characteristic classes, or better, their
difference $t(\star', \star) = c(\star') - c(\star)$ by using
techniques from \cite{gutt.rawnsley:1999a}.

The other and easier option is to use a Fedosov like construction not
only for the star product $\star$ but also for the whole bimodule
structure $\bullet$, $\bullet'$ and $\star'$ directly, by using a
connection $\nabla^L$ for $L$ in addition to the symplectic
connection, see \cite{waldmann:2002b}. Then the characteristic classes
can be trivially determined in the construction. Since every star
product is equivalent to a Fedosov star product this is sufficient to
deduce \eqref{eq:cstarprimecstarL} in general.

\begin{remark}[Deformed Vector Bundles]
    \label{remark:DeformedVectorBundles}
    The other projective modules are precisely the sections of higher
    rank vector bundles, this is just the statement of the Serre-Swan
    Theorem, see e.g.~\cite{swan:1962a}. From the previous section we
    already see that vector bundles can always be deformed, see
    also\cite{bursztyn.waldmann:2000b}.  Deforming vector bundles into
    bimodules in general is useful and interesting for its own.  In
    physics this gives the playing ground for a geometric formulation
    of the so-called \emph{non-commutative field theories}, see
    e.g.~\cite{bursztyn.waldmann:2000b} as well as
    \cite{bahns.doplicher.fredenhagen.piacitelli:2003a:pre,
      jurco.schupp.wess:2000a, jurco.schupp.wess:2001a,
      jurco.et.al.:2001a, seiberg.witten:1999a,
      connes.douglas.schwarz:1998a} and references therein.
\end{remark}

%
%

\section{Outlook: What comes next?}
\label{sec:Outlook}

To conclude this overview let us just mention a few open questions,
further developments and future projects arising from this discussion.
Some of them are work in progress.
\begin{compactenum}
\item To what extend can one apply these techniques to field theories,
    as e.g. the notions of strong Morita equivalence? In particular,
    it would be interesting to compare formal and strict deformation
    quantizations, see e.g. \cite{dito:2002a, dito:1993a, dito:1992a,
      dito:1990a, duetsch.fredenhagen:2001a,
      duetsch.fredenhagen:2001b, duetsch.fredenhagen:2003a} for
    approaches to quantum field theories using star products.
\item The state space of formal star products is in many aspects still
    not physically satisfying: it is much too big in order to allow
    physical interpretations for \emph{all} positive functionals. Thus
    one should find criteria for positive functionals to describe
    physically relevant situations. In particular, how can one
    classify deformations of classical states? Which are the `minimal'
    ones? What is the relevance of mixed and pure states from the
    deformation point of view? Can one extend the baby versions of the
    Tomita-Takesaki theorems from\cite{waldmann:2000a}?
\item Deformed line bundles and more generally deformed vector bundles
    are the starting point for any geometric description of
    non-commutative field theories. Here one has still many open
    questions concerning e.g. the global aspects of these theories,
    the convergence of star products and bimodule structures or the
    formulation of gauge transformations.
\item Since symmetries play a fundamental role in physics one has to
    investigate the invariant states and their GNS representations in
    more detail. First steps in that direction can be found in
    \cite{bordemann.neumaier.waldmann:1999a}. Here one would
    like to understand also the role of coherent states and
    eigenstates. It seems that on the purely algebraic level of formal
    star products one can not get very far but one needs some
    convergence conditions. Thus the relation between formal GNS
    construction and their convergent counterparts has to be explored.
    Since a $C^*$-algebraic is usually very difficult to obtain in a
    first step the whole situation is probably better formulated for
    some locally convex algebras. Here one can rely on the
    general results on $O^*$-algebras \cite{schmuedgen:1990a} but
    these techniques have still to be adapted to star products.
\item What is the relevance of (strong) Morita equivalence from the
    physical point of view? In particular, can one interprete the
    Morita invariants in a more physical way, like this was done for
    the Dirac monopole?
\end{compactenum}

%
%

\begin{footnotesize}
    \renewcommand{\arraystretch}{0.5} 

\end{footnotesize}

\end{document}